\documentclass[11pt, reqno, a4paper]{amsart}

\usepackage[utf8]{inputenc}
\usepackage{a4wide}

\usepackage{dsfont}
\usepackage{amssymb, amsthm, amsmath, enumerate, physics, mathrsfs, amsfonts, mathtools}
\usepackage[colorlinks=true, linkcolor=black, citecolor=black, urlcolor=black]{hyperref}
 \usepackage{accents}
\usepackage{xcolor}
\usepackage{fullpage}
\usepackage{graphicx}
\usepackage{subcaption}
\usepackage[edges]{forest}
\usepackage[normalem]{ulem}
\usepackage{soul}

\usepackage{tikz, pgfplots}
\usetikzlibrary{decorations.pathreplacing, fit}
\usepackage{tkz-euclide}

\usepackage{natbib}

\usepackage[foot]{amsaddr}

\allowdisplaybreaks
\numberwithin{equation}{section}

\newcommand{\R}{\mathbb{R}}
\newcommand{\Lc}{\mathcal{L}}
\newcommand{\C}{\mathcal{C}}
\newcommand{\E}{E}
\newcommand{\N}{\mathbb{N}}
\newcommand{\G}{\mathcal{G}}
\newcommand{\D}{\mathscr{D}}
\newcommand{\X}{\mathcal{X}}
\newcommand{\Y}{\mathcal{Y}}
\newcommand{\Pn}{\mathscr{P}}
\newcommand{\Pac}{\mathscr{P}_{\mathrm{ac}}}
\newcommand{\Pc}{\mathscr{P}_{\mathrm{c}}}
\newcommand{\T}{\mathcal{T}_c}
\newcommand{\f}{f_{\mathrm{cor}}}
\renewcommand{\L}{L^i_c\qty(\vb*{w}, \vb{z})}
\newcommand{\Lopt}{L^i_c\qty(\vb*{w}_*(\vb{z}), \vb{z})}

\usepackage{scalerel}[2014/03/10]
\usepackage[usestackEOL]{stackengine}
\def\intavg{\,\ThisStyle{\ensurestackMath{%
    \stackinset{c}{0\LMpt}{c}{0\LMpt}{\SavedStyle-}{\SavedStyle\phantom{\int}}}%
    \setbox0=\hbox{$\SavedStyle\int\,$}\kern-\wd0}\int} 

\DeclareMathOperator*{\argmax}{arg\,max}
\DeclareMathOperator*{\argmin}{arg\,min}
\newcommand{\measurerestr}{%
  \,\raisebox{-.127ex}{\reflectbox{\rotatebox[origin=br]{-90}{$\lnot$}}}\,%
}

\theoremstyle{plain}
\newtheorem{thm}{Theorem}[section]
\newtheorem{prop}[thm]{Proposition}
\newtheorem{lemma}[thm]{Lemma}
\newtheorem{cor}[thm]{Corollary}

\theoremstyle{definition}
\newtheorem{defn}[thm]{Definition}
\newtheorem{example}[thm]{Example}
\newtheorem{rmk}[thm]{Remark}
\newtheorem*{ass*}{Assumption}

\usepackage{todonotes}
\usepackage{comment} 

\newcommand{\db}[1]{\textcolor{black}{#1}}
\newcommand{\bp}[1]{\textcolor{black}{#1}}

\begin{document}

\title{Semi-discrete optimal transport techniques for the compressible semi-geostrophic equations}
\author[1]{David P.~Bourne$^*$}
\author[2]{Charlie Egan$^\dagger$} \author[3]{Th\'eo Lavier$^\#$} \author[1]{Beatrice Pelloni$^*$}

\address[*]{Maxwell Institute for Mathematical Sciences and Department of Mathematics, Heriot-Watt University, Edinburgh, UK}

\address[$\dagger$]{Department of Meteorology, University of Reading, Reading, UK}

\address[$\#$]{Laboratoire IMATH, Universit\'e de Toulon, Toulon, France}

\date{}

\begin{abstract}
We prove existence of weak solutions of the 3D compressible semi-geostrophic (SG) equations with compactly supported  measure-valued initial data.
These equations model large-scale atmospheric flows.
Our proof uses a particle discretisation 
and semi-discrete optimal transport techniques.  We show that, if the initial data is a discrete measure, then the compressible SG equations
admit a unique, energy-conserving and global-in-time solution, generated by particle trajectories that are twice continuously differentiable. In general, by discretising the initial measure by particles and sending the number of particles to infinity, we show that for any compactly supported initial measure there exists a global-in-time solution of the compressible SG equations that is Lipschitz in time. This significantly generalises the original results due to Cullen and Maroofi (2003), and it provides the theoretical foundation for the design of numerical schemes using semi-discrete optimal transport to solve the 3D compressible SG equations. 
\end{abstract}

\maketitle

\section{Introduction}

The semi-geostrophic (SG) equations for a compressible fluid are a simplified model of the formation and evolution of atmospheric fronts, which are important for accurate weather prediction on large scales. The model, 
incorporating
the effect of compressibility,  
was introduced and analysed in \cite{Cullen:2003}.

In this paper we study the system formulated in terms of the so-called \emph{geostrophic variables}. In these variables,  the governing evolution equation is a continuity equation for a time-dependent probability measure
$\alpha_t$, 
usually called the \emph{potential vorticity}.\footnote{In physical terms, $\alpha_t$ is the inverse of the physical potential vorticity; see \cite{Cullen:2021}.} The equation for $\alpha_t$ is given by
\begin{align}\label{eq:PDE}
    \partial_t\alpha_t + \nabla \cdot (\alpha_t\mathcal{W}[\alpha_t] ) = 0.
\end{align}
The non-local velocity field $\mathcal{W}[\alpha_t]$ (defined in \eqref{eq:W} below) is known as the \emph{geostrophic wind}. It is defined at each time $t$ by minimising the \emph{geostrophic energy} given the potential vorticity $\alpha_t$.
Here, and in what follows, the subscript $t$ on a dependent variable denotes evaluation at time $t$.
The PDE \eqref{eq:PDE} is derived in \cite{Cullen:2003} by starting from Euler's equations for a compressible fluid, assuming a shallow atmosphere in hydrostatic balance, a rotation-dominated flow, and making the geostrophic approximation, which is valid for large-scale flows, then transforming the resulting equations to geostrophic coordinates.

The main goal of this paper is to construct global-in-time weak solutions of \eqref{eq:PDE} as the limit of spatially discrete approximations. This generalises and extends the semi-discrete optimal transport strategy presented in \cite{Bourne:2022} for the incompressible setting, a generalisation that requires substantial changes and poses new technical challenges. 

There are several benefits of using semi-discrete optimal transport: this approach highlights an intuitive connection between flows in geostrophic coordinates and corresponding flows in the fluid domain, it provides the means to construct explicit solutions, and it serves as the theoretical foundation for adapting Cullen and Purser's groundbreaking geometric method \cite{Cullen:1984, Egan:2022} to solve the compressible SG equations numerically.

\subsection{Summary of results}
We denote by $\mathcal{\X}\subset\R^3$ the compact set representing the physical fluid domain, and by $\Y\subset \R^3$ the  \emph{geostrophic domain} where the equations are defined after the change to geostrophic variables. 
Throughout this paper we assume that $\X$ and $\Y$ satisfy the assumptions stated in 
Section \ref{sec:not_and_assump}.

To define the geostrophic energy of the system,  consider an absolutely continuous probability measure $\sigma \in \Pac(\X)$. The measure $\sigma$ is a physical variable, representing the product of the fluid density and potential temperature. For a given probability measure $\alpha_t \in \Pn(\Y)$ representing the potential vorticity, let $\T({\sigma},\alpha_t)$ be the \emph{optimal transport cost} between ${\sigma}$ and $\alpha_t$, given below by \eqref{eq:opTr}, for the cost function $c:\X\times\Y\to\R$ given by equation \eqref{eq:cost}. The corresponding \emph{geostrophic energy} $\E({\sigma},\alpha_t)$ is given by
\begin{align}\label{firstE}
    \E(\sigma, \alpha_t)=\T(\sigma,\alpha_t)+\kappa\int_{\X}\sigma^{\gamma} \,\dd \vb{x},
\end{align}
where $\gamma\in\qty(1,2)$ and $\kappa>0$ are physical constants\footnote{Interpreted physically, $\gamma$ is the ratio $\frac{c_{\mathrm{p}}}{c_{\mathrm{v}}}$, where $c_{\mathrm{p}}$ is the specific heat at constant pressure and $c_{\mathrm{v}}$ is the specific heat at constant volume, and $\kappa=c_{v}\qty(\frac{R_d}{p_0})^{\gamma-1}$, where $R_d$ is the specific gas constant for dry air, and $p_0$ is the reference pressure.}.
The geostrophic energy is the sum of the kinetic energy of the fluid, its gravitational potential energy, and its internal energy.

The geostrophic wind $\mathcal{W}$ that appears in \eqref{eq:PDE} is defined by minimising the geostrophic energy over all $\sigma$. We show that this minimisation problem has a unique solution
\begin{align}
\label{eq:sourcesimple}
    \sigma_*[\alpha_t]=\argmin_{\sigma \in \Pac(\mathcal{X})}\E\qty(\sigma,\alpha_t),
\end{align}
which we call the \emph{optimal source measure}. Let $T_{\alpha_t}$ denote the \emph{optimal transport map} from $\sigma_*[\alpha_t]$ to $\alpha_t$ for the cost $c$. Then the velocity field $\mathcal{W}[\alpha_t]:\Y \to \mathbb{R}^3$ in \eqref{eq:PDE} can be formally defined by
\begin{align}
\label{eq:W}
\mathcal{W}[\alpha_t](\vb{y}) = J \left( \vb{y}-T^{-1}_{\alpha_t}(\vb{y}) \right),\qquad  J= \f
\begin{pmatrix}
0 & -1 & 0 
\\ 1 & 0 & 0 
\\ 0 & 0 & 0
\end{pmatrix},
\end{align}
where 
$f_{\textrm{cor}} > 0$ is the Coriolis parameter, which we take to be constant.
The optimal transport map $T_{\alpha_t}$ exists and is unique because the cost function $c$ given by \eqref{eq:cost} satisfies the classical \emph{twist} condition. In general, the transport map $T_{\alpha_t}$ is not  invertible, so equation \eqref{eq:W} is only formal; we give a rigorous meaning to the PDE \eqref{eq:PDE} in Definition \ref{def:weak}. 

Equations \eqref{eq:PDE}-\eqref{eq:W} are obtained from the system of equations in physical variables (see the background section below) by the change to {\em geostrophic variables} given by
\[
T_{\alpha_t}(\vb{x})= 
\left(x_1 + \f^{-2} v^g_2(\vb{x},t),x_2 - \f^{-2} v^g_1(\vb{x},t),\theta(\vb{x},t)\right),
\]
where $\vb{v}^g=(v_1^g,v_2^g,0)$ is the geostrophic velocity of the fluid and $\theta$ is its potential temperature \cite[equation (5)]{Cullen:2003}. The physical system conserves geostrophic energy, which in these variables has the expression \eqref{firstE}. The physical requirement that solutions minimise this energy at each fixed time (known as {\em Cullen's convexity principle}) translates mathematically to the requirement that $T_{\alpha_t}$ be the optimal transport map from $\sigma_*[\alpha_t]$ to $\alpha_t$. The time evolution in geostrophic variables is governed by equation \eqref{eq:PDE}, which is coupled to the optimal transport problem through the definition of the velocity $\mathcal{W}[\alpha_t]$.
The knowledge of the map $T_{\alpha_t}$ also yields a way to recover the physical quantities $\vb{v}^g$ and $\theta$ from the components of $T_{\alpha_t}$.
This formulation is analogous to the formulation given for the incompressible system in \cite{Benamou:1998}, whereby one obtains an evolution equation such as \eqref{eq:PDE} coupled with a Monge-Amp\`ere equation for the convex function $P$ such that $T_{\alpha_t}=\nabla P$. The Monge-Amp\`ere connection is specific to the case of the quadratic cost featuring in the incompressible case. For a detailed derivation of equations \eqref{eq:PDE}--\eqref{eq:W} see \cite{Cullen:2003} and \cite{Lavier:2025}.

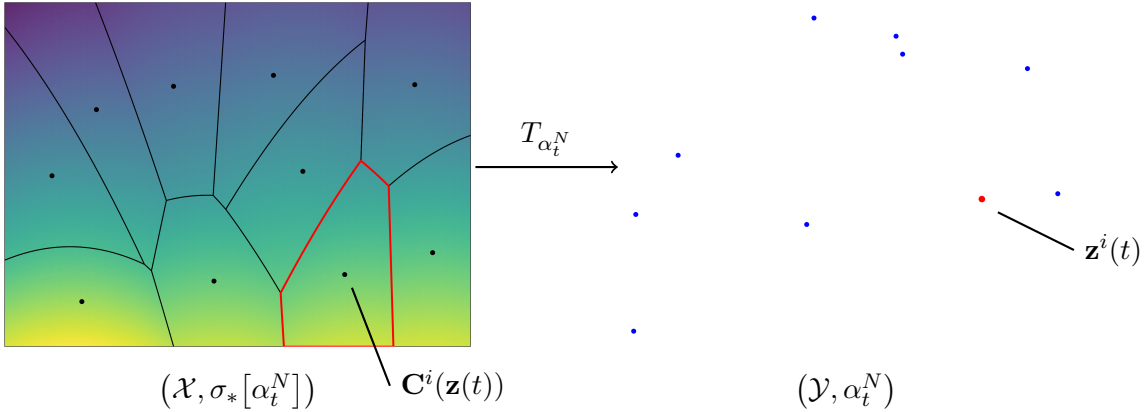
\begin{figure}[!ht]
    \centering
    \begin{tikzpicture}
        \node at (0,0) {\includegraphics[width=0.4\textwidth]{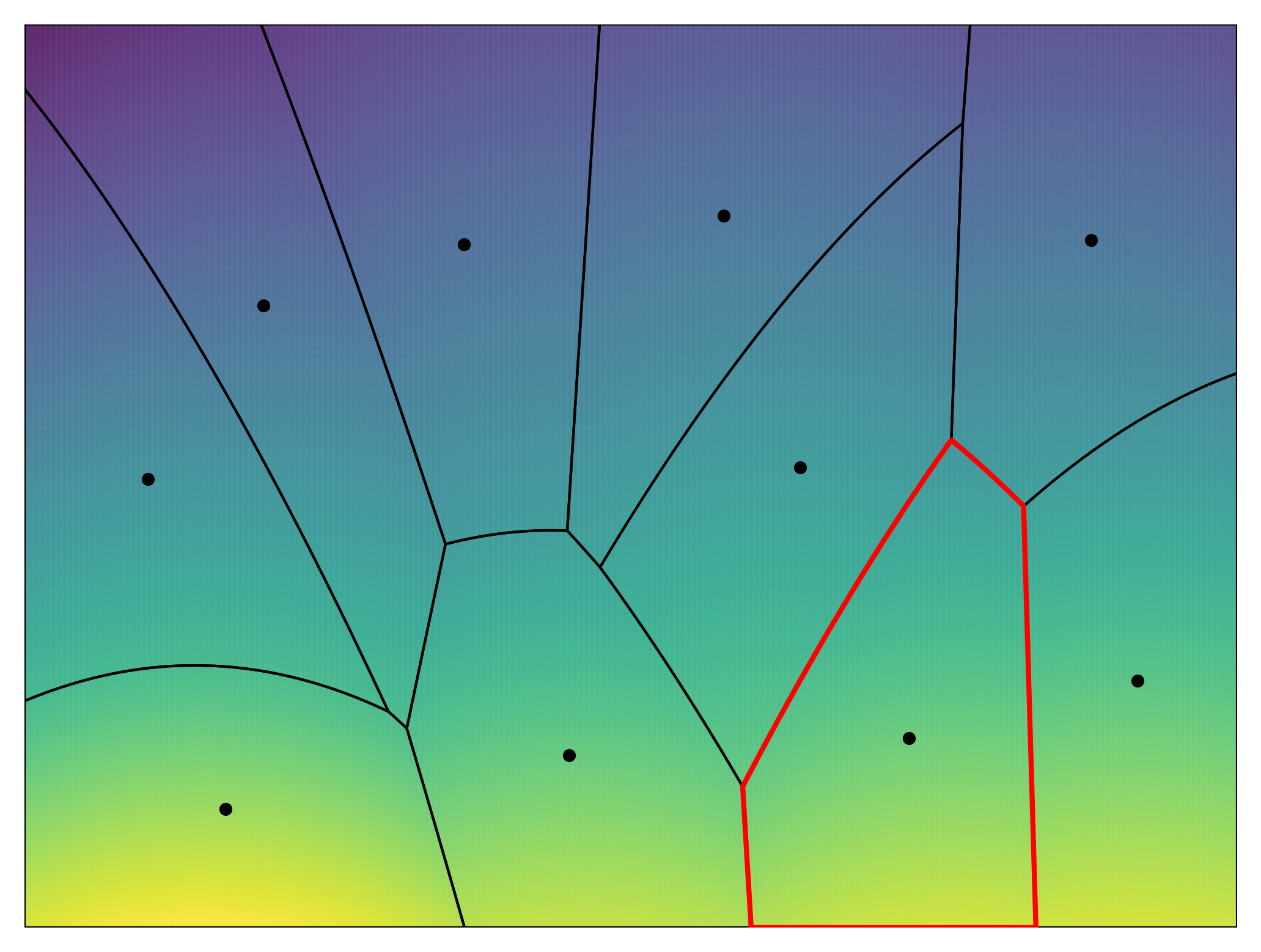}};
        \node[below] at (0,-2.5) {$\qty(\X,\sigma_*\qty[\alpha^N_t])$};
        
        \node at (8,0) {\includegraphics[width=0.4\textwidth]{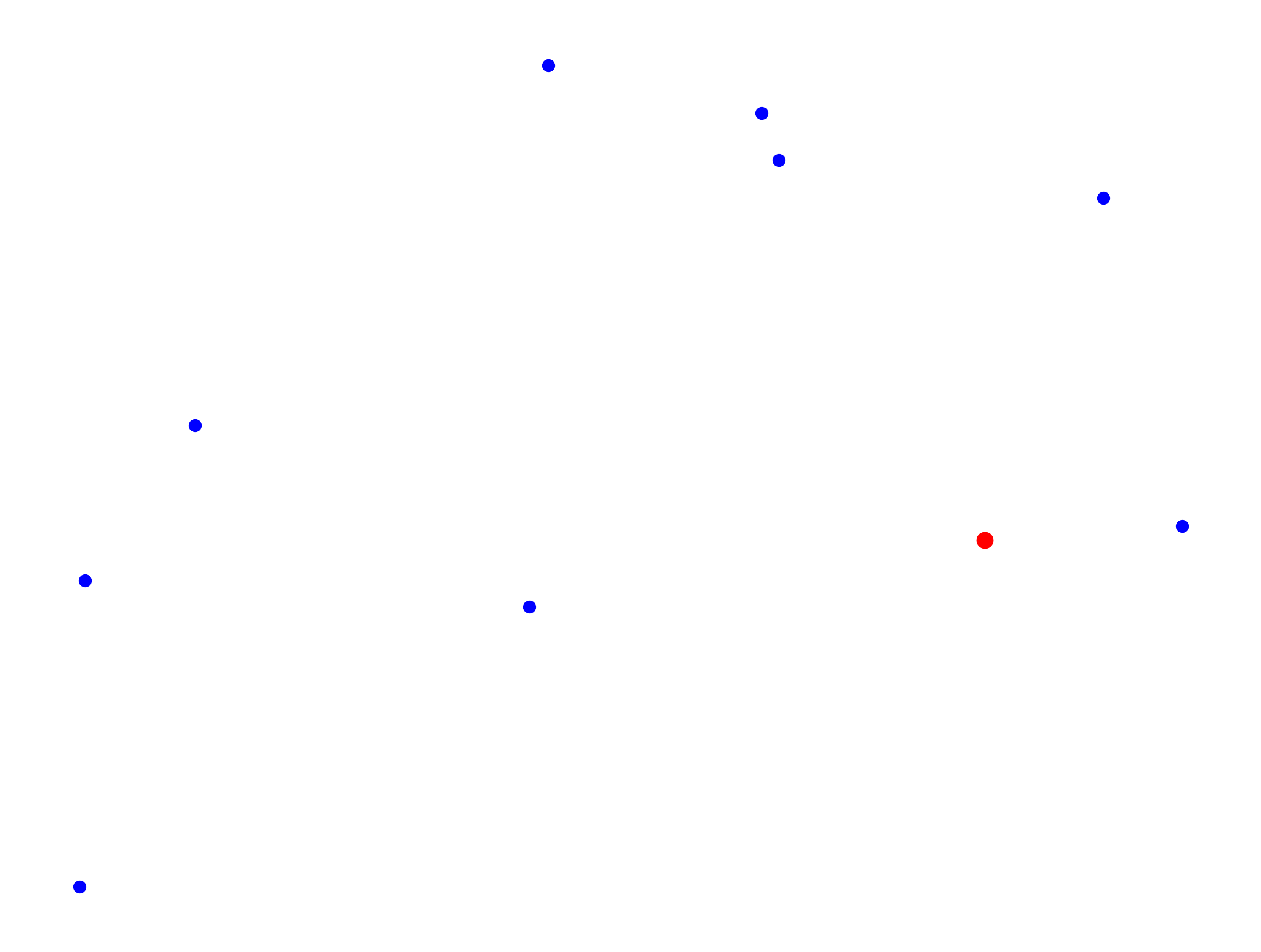}};
        \node[below] at (8, -2.5) {$\qty(\Y,\alpha^N_t)$};
        
        \draw[->, thick] (3.125,0.1) -- (5,0.1) node[above, midway] {$T_{\alpha_t^N}$};

        \draw[-, thick] (1.5, -1.5) -- (2, -2.8) node[right] {$\vb{C}^i(\vb{z}(t))$};
        \draw[-, thick] (10, -0.5) -- (11, -1) node[right] {$\vb{z}^i(t)$};
    \end{tikzpicture}
    \caption{On the left is the source space $\X$ coloured by the density of the optimal source measure $\sigma_*[\alpha_t^N]$. The boundaries and centroids $\vb{C}^j(\vb{z}(t))$ of the regions $T^{-1}_{\alpha^N_t}\{\vb{z}^i(t)\}$ are plotted in black. The boundary of the $i$-th region is highlighted in red. On the right is the target space $\Y$ with the seeds $\vb{z}^j(t)$ in blue. The $i$-th seed, $\vb{z}^i(t)$, corresponding to the $i$-th region is highlighted in red. The union of the seeds is the support of the target measure $\alpha_t^N$.}
    \label{fig1}
\end{figure}

To state the first of our results, we formulate the discrete analogue of the PDE \eqref{eq:PDE} by making the discrete ansatz
\begin{equation}
\label{eq:da}
\alpha^N_t=\sum_{i=1}^N m^i\delta_{\vb{z}^i(t)},
\end{equation}
where $\vb{z}^i(t)\in \Y$.
Let $\vb{z}(t)=\qty(\vb{z}^1(t),\ldots,\vb{z}^N(t)) \in \Y^N$.
We show in Proposition \ref{lem:Correspondence} that
this ansatz gives rise to the ODE
\begin{equation}\label{eq:ODESystem}
    \dot{\vb{z}}=J^N\qty(\vb{z}-\vb{C}(\vb{z})),
\end{equation}
where $J^N\in\R^{3N\times 3N}$ is the block diagonal matrix $ J^N:=\mathrm{diag}(J,\ldots,J)$,
and $\vb{C}$ is the \emph{centroid map}, which is defined in Definition \ref{def:centroid}.
Heuristically, the velocity of the particle $\vb{z}^i$ at time $t$ is the evaluation of the velocity field $\mathcal{W}[\alpha^N_t]$ at $\vb{z}^i(t)$. However, since $\alpha^N_t$ is discrete, the transport map $T_{\alpha^N_t}$ is not invertible and \eqref{eq:W} is not well defined. Instead, for each $i$, a whole region in $\X$ is transported to the single point $\vb{z}^i(t)$ and the term involving $T^{-1}_{\alpha^N_t}$ in \eqref{eq:W} is replaced by the centroid of the region $T^{-1}_{\alpha^N_t}\{\vb{z}^i(t)\}$ resulting in the ODE \eqref{eq:ODESystem}.  See Figure \ref{fig1}. (The sets $T^{-1}_{\alpha^N_t}\{\vb{z}^i(t)\}$ are in fact \emph{c-Laguerre cells}, which we denote by $L^i_c$. See Theorem \ref{thm: well known SDOT}.)

We now state our two main results. The notation we use is listed in Section 2. 
In particular, by $\C^{0,1}([0,\tau];\Pn(\Y))$ we mean Lipschitz continuity with respect to the $W^1$ Wasserstein metric on $\Pn(\Y)$ and the standard metric on $[0,\tau]$.

 Our first theorem asserts the existence and uniqueness of solutions of the ODE \eqref{eq:ODESystem} and correspondingly of  the PDE \eqref{eq:PDE} with discrete initial data.

\begin{thm}[Existence of discrete geostrophic solutions]
\label{thm:discretesolutions}
Let the assumptions in Section \ref{sec:not_and_assump} hold.
Fix $N\in\N$, 
an arbitrary final time $\tau>0$, and an initial discrete probability measure $\overline{\alpha}^N
= \sum_{i=1}^N m^i \delta_{\overline{\vb{z}}^i}
\in\Pn^N(\Y)$ that is well-prepared in the sense of \eqref{wellprep}.
Then there exists a unique solution 
$\vb{z}\in \C^2([0,\tau]; \Y^N)$ of the ODE \eqref{eq:ODESystem} with initial condition 
$\vb{z}(0)=(\overline{\vb{z}}^1,\ldots,\overline{\vb{z}}^N)$. Define
 $\alpha^N\in \C^{0,1}([0,\tau];\Pn^N(\Y))$ 
by
\begin{align*}   \alpha^N_t=\sum_{i=1}^N m^i\delta_{\vb{z}^i(t)}.
\end{align*}
Then $\alpha^N$ 
is the unique weak solution of the compressible SG equations (in the sense of Definition \ref{def:weak}) with discrete initial data $\alpha^N_0=\overline{\alpha}^N$.
Moreover, $\alpha^N$ is energy-conserving in the sense that
\[
\E(\sigma_*[\alpha^N_t],\alpha^N_t)=\E(\sigma_*[\alpha^N_0],\alpha^N_0) \quad \forall \; t \in [0,\tau].
\]
\end{thm}

Our second theorem states
the existence of global-in-time weak solutions of the compressible SG equations \eqref{eq:PDE} for arbitrary, compactly supported initial data.

\begin{thm}[Existence of weak geostrophic solutions]
\label{thm:solutions}
Let the assumptions in Section \ref{sec:not_and_assump} hold.
Let $\tau>0$ be an arbitrary final time, and let $\overline{\alpha}\in \Pc(\Y)$ be an initial compactly-supported probability measure. 
Then there exists $\alpha \in \C^{0,1}([0,\tau];\Pc(\Y))$ 
and $\alpha^N\in \C^{0,1}([0,\tau];\Pn^N(\Y))$, $N \in \mathbb{N}$, 
such that $\alpha$ is a weak solution of the compressible SG equations with initial measure $\overline{\alpha}$ (in the sense of Definition \ref{def:weak}), $\alpha^N$ is a discrete weak solution  
for each $N \in \mathbb{N}$, and the sequence $(\alpha^N)$ converges uniformly to $\alpha$ in the following sense:
\begin{equation}
\label{eq: uniform}
\lim_{N\to\infty}\sup_{t\in[0,\tau]}W_1(\alpha_t^N,\alpha_t)=0.
\end{equation}
\end{thm}

It is not known whether weak solutions of the compressible SG equations are unique. In particular, Theorem \ref{thm:solutions} only asserts that the weak solution that we construct can be uniformly approximated by discrete weak solutions; it does not assert that \emph{any} weak solution can be uniformly approximated in this way.

We stress that while the incompressible SG system  considered in \cite{Bourne:2022} involves the standard quadratic cost and a fixed source measure, the compressible SG system involves a non-standard cost $c$ and a time-dependent source measure,  defined by the minimisation problem \eqref{eq:sourcesimple}. This presents new, significant challenges regarding the regularity of the centroid map $\vb{C}$.

Our method of proof gives rise naturally to a numerical method. Indeed, to show that the ODE \eqref{eq:ODESystem} has a unique solution, we first derive a concave dual formulation of the energy minimisation problem \eqref{eq:sourcesimple} in the case where $\alpha_t$ is a discrete measure; see Section \ref{sec:dual_energy}. We show that this dual problem is solvable and that its solution can be used to define the optimal source measure $\sigma_*[\alpha_t]$. While we do not compute numerical solutions here, the dual problem is an unconstrained, finite-dimensional, concave maximisation problem, which is numerically tractable. (The analogous problem with the quadratic cost is solved numerically for example in \cite{Sarrazin:2022}.) Hence the velocity field in the ODE \eqref{eq:ODESystem} can, in principle, be computed numerically. Importantly, solutions of the ODE \eqref{eq:ODESystem} give rise to discrete, energy-conserving solutions of the PDE \eqref{eq:PDE}. It has been observed in the incompressible setting that numerical solutions inherit this conservation property \cite{Egan:2022,Lavier:2024}. Implementing a numerical scheme based on these observations would be a natural extension to the compressible setting of Cullen and Purser's geometric method \cite{Cullen:1984, Egan:2022}. 

\subsection{Background}
The Navier-Stokes equations are the fundamental mathematical model used by meteorologists to simulate the dynamics of the atmosphere and oceans.
Since the viscosity of air is low, the atmosphere is often modelled as inviscid, so that the diffusion coefficient in the Navier-Stokes equations can be set to zero to obtain the Euler equations.
At length scales on the order of tens of kilometres, such flows develop singularities known as atmospheric fronts. These can be modelled using the SG approximation, which is a good approximation for shallow, rotationally-dominated flows (small Rossby numbers). This approximation was first formulated in 1949 by Eliassen \cite{Eliassen:1949} and rediscovered in the 1970s by Hoskins \cite{Hoskins:1975}.
Since then it has been studied extensively (see for example the book \cite{Cullen:2021} and the bibliography therein) and the SG equations have been used as a diagnostic tool at the UK Met Office \cite{Cullen:2018}. 

The SG approximation is based on two balances that dictate the dynamics of the system. The first is \emph{hydrostatic balance}, namely the assumption that the vertical component of the pressure gradient is balanced by the gravitational force.
The second is \emph{geostrophic balance}, which assumes that the first two components of the pressure gradient are balanced by the Coriolis force \cite{Eliassen:1949}. The corresponding \emph{geostrophic wind} (in physical variables) is a large-scale, two-dimensional approximation of the fluid velocity and is directed along isobars.

 From a physical  point of view, the model described by the SG equations conserves the \emph{geostrophic energy}: the sum of kinetic, potential and internal energy.  However, mathematically, the equations do not appear to be in a conservative form, as is the case for the well-known {\em quasi-geostrophic} reduction.  This  inconsistency was considered by Hoskins \cite{Hoskins:1975}, who clarified the mathematical structure of the incompressible SG equations by an ingenious change of variables to the so-called {\em geostrophic} coordinates, sometimes called {\em dual coordinates}.

\subsection{Previous results} 
The change of coordinates to geostrophic variables was made mathematically rigorous in the groundbreaking work of Benamou and Brenier \cite{Benamou:1998}, where the incompressible SG equations are formulated in geostrophic variables as a transport equation coupled to a Monge-Amp\`ere equation. Brenier's solution of the optimal transport problem for the quadratic cost in \cite{Brenier:1991} was used to prove the existence of global-in-time weak solutions of the incompressible SG equations in geostrophic variables. This seminal work provided the foundation for much of the subsequent analysis of the SG equations. Indeed, similar reformulations are established in \cite{Cullen:2001} for the SG shallow water equations and in \cite{Cullen:2003} for the compressible SG equations, which we study in this paper, where the template of the proof given in \cite{Benamou:1998} is used to prove global existence of weak solutions.

In \cite{Cullen:2003} the existence of global-in-time solutions of the compressible SG equations is established for compactly-supported initial potential vorticity $\overline{\alpha} \in L^r(\Y)$ for $r\in(1,\infty)$.  The case $r=1$ was added by Faria in \cite{Faria:2013}.
We extend this result to $\overline{\alpha} \in \Pc(\Y)$. However, 
in \cite{Cullen:2003}
the term in the cost function $c$ representing the gravitational potential energy is taken to be any $C^2$ function whose derivative with respect to the vertical variable is non-zero, which is more general than the case we treat. The same result is proved in \cite{Cullen:2014} using the general theory for Hamiltonian ODEs in Wasserstein space developed by Ambrosio and Gangbo in \cite{Ambrosio:2008}. Additionally,  \cite{Cullen:2014} presents a proof of the existence of weak Lagrangian solutions in physical space, analogous to the results of Cullen, Feldman \& Tudorascu for the incompressible SG equations \cite{Cullen:2006, Feldman:2013, Feldman:2016, Feldman:2017}.

In \cite{Bourne:2022} a \emph{semi-discrete} optimal transport approach is used to prove the global existence of weak solutions for the \emph{incompressible} SG equations for measure-valued initial data. Exact solutions of the incompressible SG equations are constructed starting from spatially discrete initial data. These discrete solutions are then used to approximate solutions for arbitrary initial data. Similar existence results were already established in increasing levels of generality in \cite{Benamou:1998, Filho:2002, Loeper:2006, Feldman:2015}. The approach of \cite{Bourne:2022} is distinct from those of earlier papers in its use of a spatial discretisation based on semi-discrete optimal transport, rather than a combination of time discretisation and spatial mollification as used in \cite{Benamou:1998}. It is also distinct from the fully discrete approach used to study a variant of the incompressible SG equations in \cite{Cullen:2007}, where both physical and geostrophic variables are spatially discretised.
Recently entropy-regularised, fully discrete optimal transport was used to prove global existence of weak solutions for the incompressible SG equations and convergence of a numerical method \cite{Carlier:2024}.
While we cannot give here an exhaustive reference list to what is now a very substantial body of work, other important contributions to the analysis of the SG equations include the works of Ambrosio et al. \cite{Ambrosio:2012, Ambrosio:2014} and recent work on the incompressible SG equations with non-constant Coriolis parameter \cite{Silini:2023}. 
The incompressible SG equations with non-constant Coriolis parameter
was first considered in \cite{ChengCullenFeldman}.
In this paper we assume that the Coriolis parameter is constant.

Numerical methods for the SG equations date back to the development of the geometric method of Cullen \& Purser \cite{Cullen:1984}. More recently, major developments in the field of numerical optimal transport have given rise to efficient numerical methods for solving the SG equations. Semi-discrete optimal transport has been used to obtain solutions of the 2-dimensional SG Eady slice equations \cite{Egan:2022} and the 3D incompressible SG equations \cite{Lavier:2024}, while entropy-regularised discrete optimal transport 
has been used to solve the Eady slice equations 
\cite{Benamou:2023, Carlier:2024} and the semigeostrophic shallow water equations \cite{Benamou:2025}. 
Note that $\mathcal{W}[\alpha_t]$ can be formally rewritten as $\mathcal{W}[\alpha_t](\vb{y}) = - \f^{-2} \, y_3 \, J \nabla \psi_t(\vb{y})$, where $(\varphi_t,\psi_t)$ is an optimal Kantorovich potential pair for transporting $\sigma_*[\alpha_t]$ to $\alpha_t$ with cost $c$ \cite{Lavier:2025}. This expression could be convenient for applying the techniques from \cite{Benamou:2025,Benamou:2023,Carlier:2024} to the compressible SG equations. Other applications of semi-discrete optimal transport in fluid dynamics are studied for example in \cite{Gallouet:2018, Gallouet:2022}.
The finite element method has recently been used to solve a 2D version of the compressible SG equations \cite{Cotter:2025}.

\subsection{Outline of the paper}

Section \ref{sec:prelim} includes notation, background material, and the definition of a weak solution of \eqref{eq:PDE}. In Section \ref{sec:dual_energy} we derive a dual formulation of the energy minimisation problem \eqref{eq:sourcesimple} and in Section \ref{sec:discrete_solutions} we prove Theorem~\ref{thm:discretesolutions}. Our main theorem, Theorem \ref{thm:solutions}, is proved in Section \ref{sec:exist}. In Section \ref{sec:examples} we give two explicit solutions of the compressible SG equations, namely an absolutely continuous steady state solution and a discrete solution with a single particle. 
In Appendix \ref{sec:domain} we discuss a technical assumption on the fluid domain. In Appendix \ref{Sec: Proof of Lemma 3.8} we prove a technical lemma, Lemma \ref{lem:dualdifferentiability},
which is important for establishing the regularity of the centroid map $\vb{C}$.

\section{Preliminaries}
\label{sec:prelim}

\subsection{Assumptions and notation}
\label{sec:not_and_assump}
Throughout this paper we make the following assumptions:
\begin{itemize}
\item $g>0$ is the acceleration due to gravity and $f_{\textrm{cor}}> 0$ is the Coriolis parameter, which we assume to be constant.
    \item $\X\subset \R^3$ denotes the fluid domain.
    We assume that $\X$ is nonempty, connected, compact, and that it coincides with the closure of its interior. Moreover, we assume that the set $\Phi^{-1}(\X)$ is convex, where $\Phi:\R^3 \to \R^3$ is defined by
    \begin{equation}
    \label{eq: Department of Education}
    \Phi(\vb{x}) = \left( \f^{-2} \, x_1,  \,\f^{-2} \, x_2, \, g^{-1} \left( x_3 - \tfrac 12 \f^{-2} \left( x_1^2 + x_2^2 \right) \right) \right).
    \end{equation}
    In Appendix \ref{sec:domain} we give examples of domains satisfying this assumption and show that it is equivalent to assuming that $\X$ is $c$-\emph{convex} \cite[Definition 1.2]{Kitagawa:2019}, which is an important notion in regularity theory for semi-discrete optimal transport. We require this technical assumption to prove Lemma \ref{lem:dualdifferentiability}. 
    \item $\Y = \R^2\times (\delta,\tfrac{1}{\delta})$ denotes the geostrophic domain, where $\delta \in (0,1)$. This is equivalent to assuming that the initial potential temperature of the fluid is bounded away from zero and infinity.
    \end{itemize}
    We also use the following standard notation:
    \begin{itemize}
    \item $\D(\Y\times\R)$ denotes the space of real-valued infinitely-differentiable and compactly supported test functions on $\Y \times\R$.
    
    \item $\mathcal{H}^{d}$ denotes the $d$-dimensional Hausdorff measure.
    \item For a Borel set $A\subset \X$, $\mathds{1}_{A}:\X\to\{0,1\}$ denotes the characteristic function of the set $A$. That is
    $$
    \mathds{1}_{A}(\vb{x}) =
    \begin{cases}
        1 \quad \text{if}\quad \vb{x}\in A,\\
        0 \quad \text{otherwise.}
    \end{cases}
    $$
    \item $\Pn(A)$ denotes the set of Borel probability measures on $A\subseteq\R^d$.

    \item $\Pac(A)$ denotes the set of  Borel probability measures on $A\subseteq\R^d$ that are absolutely continuous with respect to the Lebesgue measure. Throughout this paper we identify absolutely continuous measures with their density functions.

\item Given a Borel map $S : A \to B$ and a measure $\mu \in \Pn(A)$, the pushforward of $\mu$ by $S$ is the measure
    $S_{\#}\mu \in \Pn(B)$ defined by $S_{\#}\mu(U)=\mu(S^{-1}(U))$ for all Borel sets $U \subseteq B$.
    
    \item $\Pc(\Y)$ denotes the set of Borel probability measures on $\Y$ with compact support. In this paper we always equip $\Pc(\Y)$ with the $W_1$ Wasserstein metric, which is defined by $W_1:\Pc(\Y) \times \Pc(\Y) \to [0,\infty)$, 
    \[
    \qquad \quad
    W_1(\alpha,\beta)=\inf \left\{ \int _{\Y\times \Y} \| \vb{y}-\tilde{\vb{y}} \| \,\dd \gamma (\vb{y},\tilde{\vb{y}}) \, : \, \gamma \in \Pn (\Y\times\Y), \, \pi_{1 \#}\gamma=\alpha, \;\pi_{2 \#}\gamma=\beta \right\}, 
    \]
    where $\pi_{1} : \Y \times \Y \to \Y$ is defined by 
$\pi_1(\vb{y},\tilde{\vb{y}})=\vb{y}$
    and
    $\pi_{2} : \Y \times \Y \to \Y$ is defined by 
$\pi_1(\vb{y},\tilde{\vb{y}})=\tilde{\vb{y}}$.
\item $\Pn^N(\Y)$ denotes the set of discrete probability measures on $\Y$ with support of cardinality $N$; see equation \eqref{eq:PN(Y)}.
\end{itemize}


\subsection{The compressible SG equations}
In this section we define weak solutions of \eqref{eq:PDE}.

\begin{defn}[The cost function]\label{def:costfunction}
The cost function $c:\X\times\Y\to\R$ associated with the compressible SG equations 
is given by
\begin{equation}\label{eq:cost}
c(\vb{x},\vb{y})= \frac{1}{y_3}\left(\frac{\f^2}{2}(x_1-y_1)^2+\frac{\f^2}{2}(x_2-y_2)^2+gx_3\right).
\end{equation}
\end{defn}

It is easy to check that $c$ is \emph{twisted}, namely that $c$ is differentiable in $\vb{x}$ and that the map $\vb{y} \mapsto \nabla_{\vb{x}}c(\vb{x},\vb{y})$ is injective for all $\vb{x} \in \X$
\cite[Definition 1.16]{Santambrogio:2015}.

Note that, while the term $gx_3$ in \eqref{eq:cost} is physically motivated, it could be
replaced by $\phi(\vb{x})$ for any function $\phi \in\C^2(\X)$ such that  $\partial_{x_3}\phi\neq 0$ (so that $c$ is twisted).
In \cite{Cullen:2003} the existence of weak solutions of the compressible SG equations is proved in this more general setting, while the assumptions on the initial data in \cite{Cullen:2003}, as well as the subsequent \cite{Cullen:2014} and \cite{Faria:2013}, are more restrictive than ours. 

\begin{defn}[Optimal transport cost]
The optimal transport cost from $\sigma \in \Pn(\X)$ to $\alpha\in \Pn(\Y)$ for the cost function $c$ is defined by
\begin{equation}
\label{eq:opTr}     \T(\sigma,\alpha):=\inf_{\substack{T:\X \to \Y \\ T_\#\sigma=\alpha}}\int_{\X}c(\vb{x},T(\vb{x})) \,\dd \sigma(\vb{x}).
\end{equation}
If $\sigma \in \Pac(\X)$ and $\alpha \in \Pc(\Y)$, then there exists a unique map $T:\X\to\Y$ achieving the minimum in \eqref{eq:opTr}; see, for example, \cite{Santambrogio:2015}. We call $T$ the \emph{optimal transport map} from $\sigma$ to $\alpha$ for the cost $c$.  
\end{defn}

\begin{defn}[Geostrophic energy]
\label{def:E}
We define the compressible geostrophic energy functional $\E:\Pn(\X)\times\Pc(\Y)\to \R \cup \{ + \infty\}$ by
\begin{equation}\label{eq:energy}
    \E(\sigma, \alpha)=
    \T(\sigma,\alpha) + F(\sigma),
\end{equation}
where $F:\Pn(\X) \to \R \cup \{ + \infty\}$ is defined by
\[
F(\sigma) = 
\begin{cases}    
\displaystyle
\kappa\int_{\X}\sigma(\vb{x})^{\gamma} \,\dd \vb{x} & \textrm{if } \sigma \in \Pac(\X),
\\
+ \infty & \textrm{otherwise},
\end{cases}
\]
where $\gamma\in\qty(1,2)$ and $\kappa>0$ are constants.
\end{defn}


\begin{lemma}[Existence and uniqueness of minimisers of $E$]
\label{lem:E unique minimiser}
Given $\alpha \in \Pc(\Y)$, the functional on $\Pn(\X)$ defined by $\sigma \mapsto E(\sigma,\alpha)$ is strictly convex and has a unique minimiser $\sigma_{*} \in \Pac(\X)$.
\end{lemma}

\begin{proof}
This is proved in \cite[Theorem 4.1]{Cullen:2003} for the special case $\alpha \in \Pac(\db{\Y})$. 
The functional $\Pn(\X) \ni \sigma \mapsto \T(\sigma,\alpha)$ is continuous and convex \cite[Propositions 7.4 \& 7.17]{Santambrogio:2015}, 
and $F$
is lower semi-continuous and strictly convex
\cite[Proposition 7.7 \& Section 7.2]{Santambrogio:2015}.
Therefore $E$ is lower semi-continuous and strictly convex. 
Since $\X$ is compact, $\Pn(\X)$ is compact with respect to weak convergence of measures. Therefore $E(\cdot,\alpha)$ has a minimiser $\sigma_*$, which belongs to $\Pac(\X)$ because $E(\cdot,\alpha)$ takes the value $+\infty$ on $\Pn(\X) \setminus \Pac(\X)$. 
Since $E$ is strictly convex, its minimiser is unique. 
\end{proof}

\begin{defn}[Optimal source measure and transport map]\label{def:optimal_source_general}
Given $\alpha\in\Pc(\Y)$, we define the map $\sigma_{*}:\Pc(\Y)\to\Pac(\X)$ by
\begin{equation}\label{eq:sigma_min_general}
    \sigma_*[\alpha]=\argmin_{\sigma\in\Pac(\X)}\E\qty(\sigma,\alpha).
\end{equation}
We denote by $T_{\alpha}:\X \to \Y$ the optimal transport map from $\sigma_*[\alpha]$ to $\alpha$ for the cost $c$.
\end{defn}

The formulation of the compressible SG equations given in \eqref{eq:PDE} and \eqref{eq:W}  is formal because the inverse of $T_{\alpha_t}$ may not be well defined. 
For example, if $\alpha_t$ is discrete, then $T_{\alpha_t}$ sends a set of positive Lebesgue measure to each point in the support of $\alpha_t$, violating injectivity. 
We therefore introduce a suitable weak formulation of \eqref{eq:PDE}.
 
To derive the weak formulation, we 
assume provisionally that $T_{\alpha_t}^{-1}$ exists and that all functions appearing in \eqref{eq:PDE} are sufficiently regular. Multiplying \eqref{eq:PDE} by $\varphi\in\D(\Y\times\R)$ and integrating by parts over the domain $\Y\times\qty[0,\tau]$ yields
\[
\int_0^{\tau}\int_{\Y}
\qty(\partial_t\varphi_t(\vb{y})+J\qty(\vb{y}-T^{-1}_{\alpha_t}(\vb{y}))\cdot\nabla\varphi_t(\vb{y})) \,\dd \alpha_t(\vb{y}) \dd t 
        = \int_{\Y}\varphi_{\tau}(\vb{y}) \,\dd \alpha_{\tau}(\vb{y}) -\int_{\Y}\varphi_0(\vb{y}) \,\dd \alpha_0(\vb{y}),
\]
where we have used the notation $\varphi_t(\vb{y}):=\varphi(\vb{y},t)$.
By using the push-forward constraint
$(T_{\alpha_t})_{\#}\sigma_*[\alpha_t]=\alpha_t$ to rewrite the term involving  $T_{\alpha_t}^{-1}$, we arrive at the following weak formulation of the compressible SG equations.

\begin{defn}[Weak solution]\label{def:weak}
We say that $\alpha \in \C([0,\tau];\Pc(\mathcal{Y}))$ is a \emph{weak solution} of  
the compressible SG equation \eqref{eq:PDE} with initial condition $\alpha_0=\overline{\alpha} \in \Pc(\Y)$ if for all $\varphi\in\D(\Y\times\R)$,
\begin{multline}
\label{eq:WeakFormPDE}
\int_0^{\tau} \int_{\Y} 
\left[ \partial_t \varphi_t (\vb{y}) + (J\vb{y}) \cdot \nabla \varphi_t(\vb{y}) \right] \,\dd \alpha_t(\vb{y}) 
\, \dd t
- \int_0^{\tau} \int_{\X}(J\vb{x}) \cdot \nabla \varphi_t (T_{\alpha_t}(\vb{x})) \,\dd \sigma_*[\alpha_t](\vb{x}) \, \dd t 
\\
= \int_{\Y} \varphi_{\tau} (\vb{y}) \, \dd \alpha_{\tau}(\vb{y}) -\int_{\Y} \varphi_0 (\vb{y}) \,\dd \overline{\alpha} (\vb{y}).
\end{multline}
\end{defn}

\begin{rmk}[Equivalent weak formulation]
Following \cite{Benamou:2023,Carlier:2024}, the nonlinear term in \eqref{eq:WeakFormPDE} can be rewritten as a term that is \emph{linear} in the optimal transport plan:
\begin{equation}
\label{eq: Greenland}    
\int_0^{\tau} \int_{\X}(J\vb{x}) \cdot \nabla \varphi_t (T_{\alpha_t}(\vb{x})) \,\dd \sigma_*[\alpha_t](\vb{x}) \, \dd t 
= 
\int_0^{\tau} \int_{\X \times \Y}
(J\vb{x}) \cdot \nabla \varphi_t (\vb{y}) \, \dd \gamma[\alpha_t](\vb{x},\vb{y}) \, \dd t,
\end{equation}
where $\gamma[\alpha_t] = (\mathrm{id}_\X \times T_{\alpha_t})_\# \sigma_*[\alpha_t]$ is the optimal plan for transporting $\sigma_*[\alpha_t]$ to $\alpha_t$.
\end{rmk}

\subsection{Semi-discrete optimal transport}
\label{sec:SDOT}
In this section we summarise the definitions and results from semi-discrete optimal transport theory used throughout this paper. Here $N\in\N$ is any natural number.

\begin{itemize}
    \item The set $D^N$ of {\em seed vectors} $\vb{z}$ is given by
        \begin{equation*}
           D^N:=\left\{ \vb{z}= \left(\vb{z}^1,\ldots,\vb{z}^N \right) \in\Y^{N}:\vb{z}^i\neq \vb{z}^j \textrm{ whenever } i\neq j \right\}.
        \end{equation*}
    \item We define
        \begin{equation*}
        \Delta^N:=\left\{ \vb{m}=(m^1,\ldots,m^N) \in (0,1)^N : \sum_{i=1}^N m^i = 1\right\}.
        \end{equation*}
    \item The class of discrete probability measures with $N$ masses is 
        \begin{equation}
        \label{eq:PN(Y)}
            \Pn^N(\Y):=\left\{ \sum_{i=1}^N m^i\delta_{\vb{z}^i}:\vb{z}\in D^N, \;\vb{m}\in\Delta^N \right\}.
        \end{equation}
    \item We say that a discrete probability measure 
        $\displaystyle{\sum_{i=1}^N m^i\delta_{\vb{z}^i} \in \Pn^N(\Y)}$
        is {\em well-prepared} if 
        \begin{equation}\label{wellprep}
            \vb{z}\in D_0^N \textrm{ where } D_0^N \coloneqq \left\{ \vb{z} \in D^N \, : \, 
            z^i_3 \ne z^j_3 \; \forall \; 
            i,j \in \{1,\ldots,N\}, \, i \ne j
            \right\}.
        \end{equation}
        In other words, a discrete measure is well-prepared if the seeds $\vb{z}^i$ lie in distinct horizontal planes. It will be crucial to assume that the initial data is well-prepared to prove existence and uniqueness for the ODE \eqref{eq:ODESystem}; see Proposition \ref{prop:ODESolutions}. Solutions of the ODE \eqref{eq:ODESystem} satisfy $z^i_3(t) \ne z^j_3(t)$ for all $i \ne j$ and all $t$, provided that $z^i_3(0) \ne z^j_3(0)$ (because $\dot{z}^k_3 = 0$ for all $k$). Consequently the particles $\vb{z}^i$ cannot collide, which would preclude global existence of solutions because the vector field in the ODE \eqref{eq:ODESystem} is not defined if the seeds are not distinct (the centroid map is not defined).
        
\end{itemize}

\begin{defn}[$c$-Laguerre tessellation]
\label{def:cLagCells} 
Given $(\vb*{w},\vb{z}) \in \R^N\times D^N$, the \emph{$c$-Laguerre tessellation of $\X$ generated by $(\vb*{w},\vb{z})$} is the collection of \emph{Laguerre cells} $\{ \L \}_{i=1}^N$ defined by
\begin{align*}
    \L:=\qty{\vb{x}\in \X:c(\vb{x},\vb{z}^i)-w^i\leq c(\vb{x},\vb{z}^j)-w^j\; \forall\,j\in\qty{1,\ldots,N}},
    \quad i \in \{1,\ldots,N\}.
\end{align*}
\end{defn}

\begin{figure}[t]
  \centering
  \begin{subfigure}{0.4\textwidth}
    \centering
    \includegraphics[width=\linewidth, angle = 180]{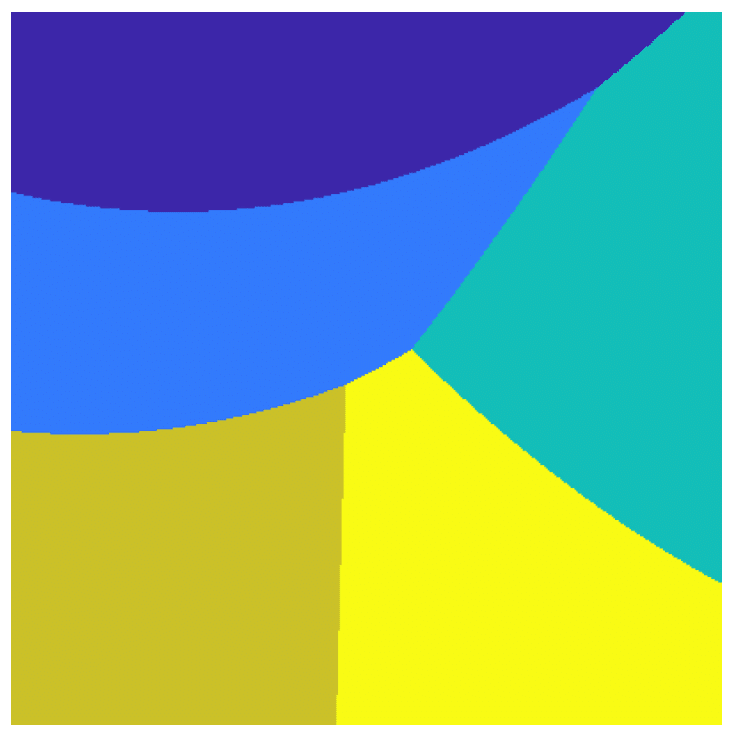}
    \caption{5 Laguerre cells}
    \label{fig:5Cells}
  \end{subfigure}%
  \hfill
  \begin{subfigure}{0.4\textwidth}
    \centering
    \includegraphics[width=\linewidth, angle = 180]{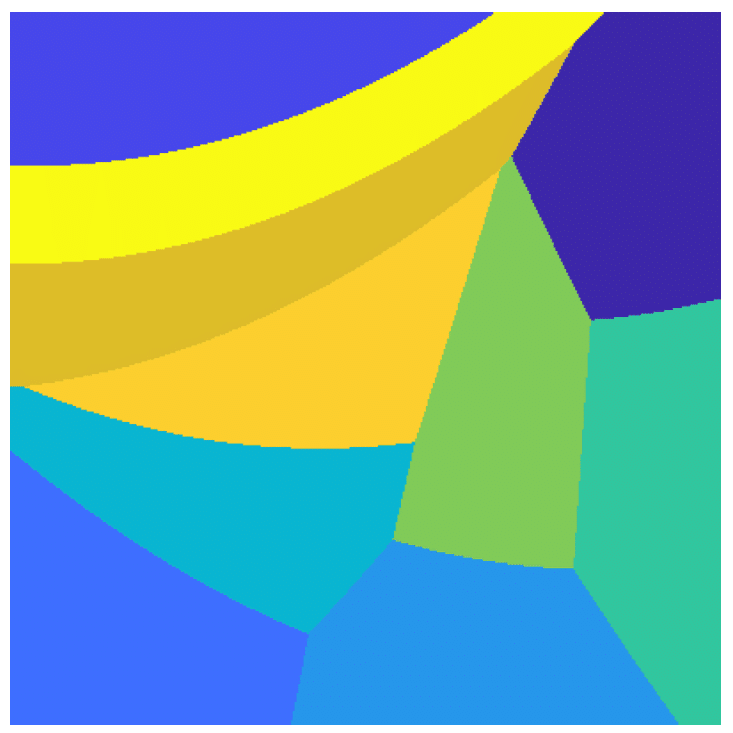}
    \caption{10 Laguerre cells}
    \label{fig:10Cells}
  \end{subfigure}

  \vspace{1em}

  \begin{subfigure}{0.4\textwidth}
    \centering
    \includegraphics[width=\linewidth, angle = 180]{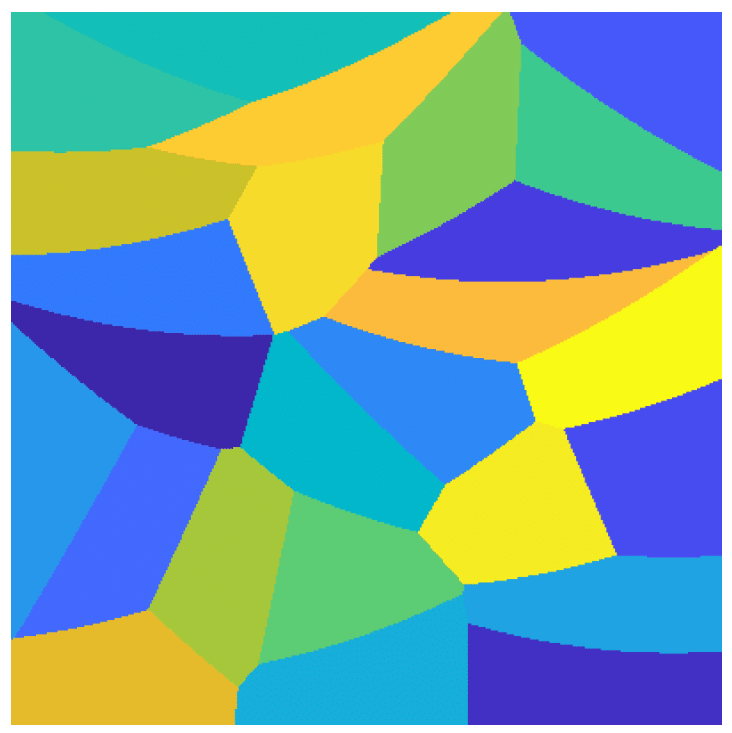}
    \caption{25 Laguerre cells}
    \label{fig:25Cells}
  \end{subfigure}%
  \hfill
  \begin{subfigure}{0.4\textwidth}
    \centering
    \includegraphics[width=\linewidth, angle = 180]{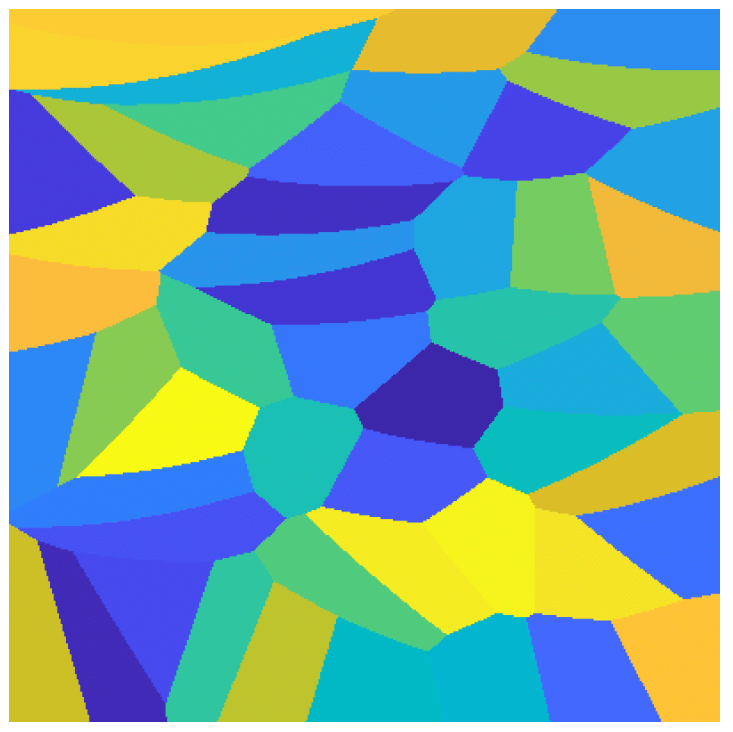}
    \caption{50 Laguerre cells}
    \label{fig:50Cells}
  \end{subfigure}

  \caption{$c$-Laguerre tessellations
   (see Definition \ref{def:cLagCells}) in the $(x_1,x_3)$-plane for the cost function $c=c_{\mathrm{2d}}$
   (see \eqref{eq:c2d}). The colours distinguish the cells.
  For each plot, $\X=[0,1]^2$, $f_{\textrm{cor}}=1$, $g=1$, the seeds $\vb{z}^i$ were sampled uniformly from $\X$, and the weights $w^i$ were chosen so that the cells have equal area (by maximising the dual function as in \eqref{eq:optdiscT}).}
  \label{fig:CompressibleCells}
\end{figure}

The $c$-Laguerre cells 
form a tessellation of $\X$ in the sense that $\bigcup_{i=1}^N L^i_c(\vb*{w},\vb{z}) = \X$ and $\Lc^3 \big( 
L^i_c(\vb*{w},\vb{z}) \cap L^j_c(\vb*{w},\vb{z}) \big) = 0$ if $i \ne j$ (by \cite[Proposition 37]{Merigot:2020}). 
Intuitively, changing the \emph{weight vector} $\vb*{w}=(w^1,\ldots,w^N)$ shifts the boundaries of the Laguerre cells, thus controlling their volumes.
It is easy to check that if $\vb{z} \in D_0^N$, then each Laguerre cell $\L$ is the intersection of $\X$ with $N-1$ paraboloids of the form
\[
\left\{ \vb{x} \in \R^3 : x_3 \le
    -\frac{\f^2}{2g}\left( (x_1-a_1)^2 + (x_2-a_2)^2 \right) + b
\right\},
\]
where $a_1, a_2 \in \R$ depend on $\vb{z}$ and $b \in \R$ depends on $(\vb*{w},\vb{z})$. 
In Figure~\ref{fig:CompressibleCells} we plot $c$-Laguerre tessellations in two dimensions for the following 2D version of the cost \eqref{eq:cost}:
\begin{align}
\label{eq:c2d}
    c_{\mathrm{2d}}((x_1,x_3),(y_1,y_3)) = \frac{1}{y_3}\qty(\frac{f_{\mathrm{cor}
    }^2}{2}(x_1-y_1)^2+gx_3).
\end{align}

We recall the following well-known result in semi-discrete optimal transport theory, which we will use frequently throughout the paper. For a proof see for example \cite[Section 4]{Merigot:2020}.
\begin{thm}[The optimal transport map for semi-discrete optimal transport]
\label{thm: well known SDOT}
Let
$\sigma\in\mathcal{P}_{\mathrm{ac}}(\mathcal{X})$ and $\alpha^N = \sum_{i=1}^N m^i\delta_{\vb{z}^i} \in \Pn^N(\Y)$. 
There there exists  $\vb*{w}\in\R^N$ satisfying
\begin{equation}
\label{eq:massConstraint}
\sigma(\L) = m^i \quad \forall\, i \in \{ 1,\ldots, N \}.
\end{equation}
Given $\vb*{w}$ satisfying \eqref{eq:massConstraint}, define $T:\X \to \Y$ by
\begin{equation}
\label{eq:TSOT}
T = \sum_{i=1}^{N} \vb{z}^i \mathds{1}_{\L}.
\end{equation}
In other words, $T$ is the piecewise constant function given by $T(\vb{x})=\vb{z}^i$ if $\vb{x} \in \L$.
Then $T:\X \to \Y$ is the unique optimal transport map from $\sigma$ to $\alpha^N$ for the cost $c$. Moreover, the following Kantorovich duality holds:
\begin{equation}
\label{eq:optdiscT}
\mathcal{T}_c(\sigma,\alpha^N)
= \max_{\widetilde{\vb*{w}}\in\R^N} \sum_{i=1}^N \left( m^i\widetilde{w}^i+\int_{L^i_c(\widetilde{\vb*{w}},\vb{z})}\qty(c(\vb{x},\vb{z}^i)-\widetilde{w}^i) \,\dd \sigma(\vb{x}) \right).
\end{equation}
Furthermore, $\vb*{w} \in \R^N$ satisfies \eqref{eq:massConstraint} if and only if it maximises the dual function in \eqref{eq:optdiscT}.  
\end{thm}
Note that the solution $\vb*{w} \in \R^N$ of \eqref{eq:massConstraint} is not unique because $L^i_c(\vb*{w}+\vb*{a},\vb{z})=\L$ for any constant vector $\vb*{a}=(a,\ldots,a) \in \R^N$. Also note that
the map $T$ defined in \eqref{eq:TSOT} is well-defined $\Lc^3$-almost everywhere because the intersection of any two $c$-Laguerre cells has zero Lebesgue measure by \cite[Proposition 37] {Merigot:2020} and the twist condition on $c$.

\section{The dual problem}
\label{sec:dual_energy}
In this section we derive a dual formulation of the minimisation problem \eqref{eq:sigma_min_general} for the case where $\alpha$ is a discrete measure. Additionally, we prove necessary and sufficient optimality conditions that characterise the optimal source measure $\sigma_*[\alpha]$ in terms of the solution of the dual problem. The main results are contained in Theorem \ref{thm:dualitytheorem}.

Given $(\vb*{w},\vb{z}) \in \R^N \times D^N$, recall that the $c$-transform of $\vb*{w}$ is the function $\vb*{w}^c (\, \cdot \, ;\vb{z}): \X \to \mathbb{R}$ defined by
\[
\vb*{w}^c(\vb{x};\vb{z}) = \min_{i \in \{1,\ldots,N\}} \qty{ c(\vb{x},\vb{z}^i) - w^i }.
\]
In particular, $\vb*{w}^c(\vb{x};\vb{z}) = c(\vb{x},\vb{z}^i) - w^i$ if and only if $\vb{x} \in L^i_c(\vb*{w},\vb{z})$.

\begin{defn}[Dual functional]
Given $\vb{m} \in \Delta^N$, 
    define the dual functional $\G:\R^N \times D^N \to\R$ by
\[
    \G(\vb*{w},\vb{z})
    = \sum_{i=1}^N m^iw^i - \int_{\X} f^*(-\vb*{w}^c(\vb{x};\vb{z})) \, \dd \vb{x} 
    \\ 
    = \sum_{i=1}^N \left( m^iw^i-\int_{\L}f^*(w^i-c(\vb{x},\vb{z}^i)) \,\dd \vb{x} \right),
\]
where $f^*:\R\to\R$ is the Legendre-Fenchel transform of the function $f:\R\to\R$ given by
\begin{equation*}
    f(s)=\begin{cases}
    \kappa s^{\gamma} &\text{if }  s \ge 0, \\ 
    +\infty &\text{if } s < 0, 
    \end{cases}
\end{equation*}
where $\gamma\in (1,2)$ and $\kappa>0$. 
\end{defn}

\begin{rmk}[Properties of $f$] \label{remf}
The function $f$ is strictly convex on $[0,\infty)$, lower semi-continuous, and its Legendre-Fenchel transform is given by
\begin{align*}
    f^*(t)\coloneqq\sup_{s\in\R}\, (s t - f(s))
    =\begin{cases}
    \frac{1}{\gamma'}
    (\kappa\gamma)^{1-\gamma'}t^{\gamma'} &\text{if } t>0,\\
    0 &\text{if } t\leq 0,
    \end{cases}
\end{align*}
where $\gamma' \in (2,\infty)$ satisfies $\frac{1}{\gamma}+\frac{1}{\gamma'}=1$. Note that $f^* \in \C^2(\R)$.
\end{rmk}

\begin{rmk}[Derivation of the dual functional]
We give a short
derivation of the dual functional before stating a rigorous duality theorem below (see Theorem \ref{thm:dualitytheorem}).
Here we identify vectors $\vb*{w} \in \mathbb{R}^N$ with real-valued functions on the discrete set $\{ \vb{z}^i \}_{i=1}^N$ via $\vb*{w}(\vb{z}^i) = w^i$. 
Applying the standard Kantorovich Duality Theorem from optimal transport theory 
gives
\begin{align*}
\inf_{\sigma\in\Pac(\X)} E(\sigma,\alpha^N) 
& = 
\inf_{\sigma\in\Pac(\X)}   
\Big( {\mathcal{T}_c(\sigma,\alpha^N )} +
F(\sigma)
\Big)
\\
& = \inf_{\sigma\in\Pac(\X)} 
\left(
\sup_{\vb*{w} \in \mathbb{R}^N}
\left( \int_{\X} \vb*{w}^c \, \mathrm{d} \sigma + \int_{\Y}  \vb*{w} \, \mathrm{d} \alpha^N 
\right)
+ F(\sigma)
\right)
\\
& \ge \sup_{\vb*{w} \in \mathbb{R}^N}
 \inf_{\sigma\in\Pac(\X)}
 \left(
\int_{\X} \vb*{w}^c \, \mathrm{d} \sigma + \int_{\Y}  \vb*{w} \, \mathrm{d} \alpha^N
+ F(\sigma)
\right)
\\
& = 
\sup_{\vb*{w} \in \mathbb{R}^N}
\bigg(
 \int_{\Y}  \vb*{w} \, \mathrm{d} \alpha^N
- \sup_{\sigma\in\Pac(\X)}
\left(
\int_{\X} (-\vb*{w}^c) \, \mathrm{d} \sigma
- F(\sigma)
\right)
\bigg)
\\
& = \sup_{\vb*{w} \in \mathbb{R}^N}
\bigg(
 \sum_{i=1}^N m^i w^i
- \int_{\X} f^* ( - \vb*{w}^c)
 \, \mathrm{d} \vb{x}
 \bigg)
  \\
 & = \sup_{\vb*{w} \in \mathbb{R}^N} \G (\vb*{w},\vb{z}).
 \end{align*}
 In Theorem \ref{thm:dualitytheorem} we prove that the $\inf$ on the left-hand side is in fact a $\min$, the $\sup$ on the right-hand side is in fact a $\max$, and the inequality is an equality.
\end{rmk}

\begin{thm}[Duality Theorem]
\label{thm:dualitytheorem} 
Let $\vb{z} \in D^N$, $\vb{m} \in \Delta^N$ and define
\[
\alpha^N:=\sum_{i=1}^Nm^i\delta_{\vb{z}^i}.
\]
Assume that $z^i_3 \ne z^j_3$ for all $i,j \in \{1,\ldots,N\}$ with $i \ne j$.
Then
the map $\sigma \mapsto E(\sigma,\alpha^N)$ is strictly convex and has a unique minimiser, and the map $\vb*{w} \mapsto \G(\vb*{w},\vb{z})$ is concave and has a unique maximiser.
Moreover, 
\begin{equation}
\label{eq:duality}    
    \min_{\sigma\in\Pac(\X)} \E(\sigma,\alpha^N) = \max_{\vb*{w}\in\R^N} \G(\vb*{w},\vb{z}).
\end{equation}
Furthermore, $\sigma \in\Pac(\X)$ minimises $E(\cdot,\alpha^N)$ and 
$\vb*{w} \in\R^N$ maximises $\G(\cdot,\vb{z})$
if and only if
\begin{align} 
\label{eq:gMassConstraint}
    m^i & = \int_{L_c^i(\vb*{w},\vb{z})} (f^*)'(w^i-c(\vb{x},\vb{z}^i)) \,\dd \vb{x} \quad \forall \; i\in\qty{1,\ldots,N},\\
    \label{eq:density}
    \sigma(\vb{x}) 
    & = (f^*)'(-\vb*{w}^c(\vb{x};\vb{z})) \quad \textrm{for } \mathcal{L}^3\textrm{-almost every } \vb{x} \in \X. 
\end{align}
\end{thm}

Theorem \ref{thm:dualitytheorem} is analogous to the duality result given in \cite[Proposition 12]{Sarrazin:2022} but with the cost $c$ rather than the quadratic cost. 
Other closely related results include \cite[Proposition 5.3]{Leclerc:2020} and \cite[Theorems 3.1 and 3.2]{Bourne:2018}. In addition, the Euler-Lagrange equation \eqref{eq:density} can be derived from \cite[Proposition 7.20]{Santambrogio:2015}.
 
\begin{rmk}[Optimal transport map]
If $\sigma \in\Pac(\X)$ minimises $E(\cdot,\alpha^N)$ and 
$\vb*{w} \in\R^N$ maximises $\G(\cdot,\vb{z})$,
then the optimality conditions \eqref{eq:gMassConstraint} and \eqref{eq:density} imply that the mass constraint \eqref{eq:massConstraint} holds. In particular, this means that the map $T$ given by 
\eqref{eq:TSOT}
is the optimal transport map from $\sigma$ to $\alpha^N$ for the cost $c$.
\end{rmk}

We begin by establishing the regularity of $\G$; see Proposition \ref{lem:gderivatvies}. 
To do this we use the following two lemmas, 
which will also be used in Section \ref{sec:discrete_solutions} to prove that the centroid map is continuously differentiable.

\begin{lemma}[Regularity of integrals over Laguerre cells]
\label{lem:dualdifferentiability}
Let $U \subset \R^N \times D_0^N$ be the open set
\[
U = \bigg\{ (\vb*{w},\vb{z}) \in \R^N \times D_0^N \, : \, \mathcal{L}^3 (L^i_c(\vb*{w},\vb{z})) > 0 \; \forall \; i \in \{1,\ldots,N\} \bigg\}.
\]
Define 
$\vb{\Psi}=(\Psi^1,\ldots,\Psi^N):U\to\R^N$
by
\begin{equation*}
    \Psi^i(\vb*{w},\vb{z}):=\int_{\L} \zeta(\vb{x},\vb*{w},\vb{z}) \,\dd \vb{x},
\end{equation*}
where
$\zeta: \X \times \R^N \times D^N \to\R$
satisfies the following:
\begin{itemize}
    \item 
$\zeta(\cdot,\vb*{w}, \vb{z})\in\C(\X)$
for all
$(\vb*{w},\vb{z}) \in \R^N \times D^N$;
\item
for each compact set $K \subset \R^N \times D^N$, 
$\zeta(\vb{x},\cdot,\cdot)\in\C^{0,1}(K)$ 
for all $\vb{x} \in \X$ and 
there exists a constant $L(K)>0$ such that $\sup_{\vb{x} \in \X} \mathrm{Lip}(\zeta(\vb{x},\cdot,\cdot)|_K) \le L(K)$;
\item for each
$(\vb*{w}_0,\vb{z}_0) \in \R^N \times D^N$,
 there exists an open set $S(\vb*{w}_0,\vb{z}_0) \subseteq \X$ with $\mathcal{L}^3(\X \setminus S(\vb*{w}_0,\vb{z}_0))=0$ such that 
the partial derivatives 
$\partial \zeta/\partial \vb*{w}(\vb{x},\cdot,\cdot)$, $\partial \zeta/\partial \vb{z}(\vb{x},\cdot,\cdot)$ exist and are continuous at $(\vb*{w}_0,\vb{z}_0)$ for all $\vb{x} \in S(\vb*{w}_0,\vb{z}_0)$.
\end{itemize}
Then $\vb{\Psi}$
is continuously differentiable. Moreover, for all $i,j \in \{1,\ldots,N\}$, $i \ne j$,
\begin{align*}
\frac{\partial \Psi^i}{\partial w^j}(\vb*{w},\vb{z}) 
& = 
\int_{L^i_c(\vb*{w},\vb{z})} \frac{\partial \zeta}{\partial w^j} (\vb{x},\vb*{w},\vb{z})  \, \dd \vb{x}
- \int_{L^i_c(\vb*{w},\vb{z}) \cap L^j_c(\vb*{w},\vb{z})} 
\frac{\zeta(\vb{x},\vb*{w},\vb{z})}{\| \nabla_{\vb{x}} c(\vb{x},\vb{z}^i) -  \nabla_{\vb{x}} c(\vb{x},\vb{z}^j) \|} \, \dd  \mathcal{H}^2(\vb{x}),
\\
\frac{\partial \Psi^i}{\partial \vb{z}^j}(\vb*{w},\vb{z}) 
& = 
\int_{L^i_c(\vb*{w},\vb{z})} \frac{\partial \zeta}{\partial \vb{z}^j} (\vb{x},\vb*{w},\vb{z})  \, \dd \vb{x}
+
\int_{L^i_c(\vb*{w},\vb{z}) \cap L^j_c(\vb*{w},\vb{z})} 
\frac{\nabla_{\vb{y}}c(\vb{x},\vb{z}^j) \, \zeta(\vb{x},\vb*{w},\vb{z})}{\| \nabla_{\vb{x}} c(\vb{x},\vb{z}^i) -  \nabla_{\vb{x}} c(\vb{x},\vb{z}^j) \|} \, \dd  \mathcal{H}^2(\vb{x}),
\end{align*}
and
\begin{align*}
\frac{\partial \Psi^i}{\partial w^i}(\vb*{w},\vb{z}) 
& = 
\int_{L^i_c(\vb*{w},\vb{z})} \frac{\partial \zeta}{\partial w^i} (\vb{x},\vb*{w},\vb{z})  \, \dd \vb{x}
+ \sum_{j \ne i} \int_{L^i_c(\vb*{w},\vb{z}) \cap L^j_c(\vb*{w},\vb{z})} 
\frac{\zeta(\vb{x},\vb*{w},\vb{z})}{\| \nabla_{\vb{x}} c(\vb{x},\vb{z}^i) -  \nabla_{\vb{x}} c(\vb{x},\vb{z}^j) \|} \, \dd  \mathcal{H}^2(\vb{x}),
\\    
\frac{\partial \Psi^i}{\partial \vb{z}^i}(\vb*{w},\vb{z}) 
& =
\int_{L^i_c(\vb*{w},\vb{z})} \frac{\partial \zeta}{\partial \vb{z}^i} (\vb{x},\vb*{w},\vb{z})  \, \dd \vb{x}
- \sum_{j \ne i} 
\int_{L^i_c(\vb*{w},\vb{z}) \cap L^j_c(\vb*{w},\vb{z})} 
\frac{\nabla_{\vb{y}}c(\vb{x},\vb{z}^j) \, \zeta(\vb{x},\vb*{w},\vb{z})}{\| \nabla_{\vb{x}} c(\vb{x},\vb{z}^i) -  \nabla_{\vb{x}} c(\vb{x},\vb{z}^j) \|} \, \dd  \mathcal{H}^2(\vb{x}).
\end{align*}
\end{lemma}

\begin{proof}
The proof is rather technical and is given in Appendix \ref{Sec: Proof of Lemma 3.8}. 
\end{proof}

Note that the set $U$ in the statement of Lemma \ref{lem:dualdifferentiability} is open because the map $(\vb*{w},\vb{z}) \mapsto \mathcal{L}^3 (L^i_c(\vb*{w},\vb{z}))$ is continuous.
(Continuity of this map with respect to $\vb*{w}$ is proved for example in \cite[Proposition 38(vii)]{Merigot:2020}. The proof of continuity with respect to $(\vb*{w},\vb{z})$ is identical.)
The constraint that the seeds in $U$ satisfy $|z^i_3 - z^j_3| > 0$ is not necessary, but it is sufficient for our purposes and it slightly simplifies the proof; see Remark \ref{remark: B.2}.

\begin{lemma}
\label{lem:fregularity}
Define $\zeta: \X \times \R^N \times D^N \to\R$ by
\begin{align*}
    \zeta(\vb{x},\vb*{w},\vb{z}):=
    (f^*)'(-\vb*{w}^c(\vb{x};\vb{z})).
\end{align*}
Then $\zeta$ satisfies the assumptions of Lemma~\ref{lem:dualdifferentiability}.
\end{lemma}
\begin{proof}
Since $(f^*)'$ is non-decreasing, we have that 
\[
\zeta(\vb{x},\vb*{w},\vb{z})
= (f^*)' \Big( -\min_{i \in \{1,\ldots,N\}} (c(\vb{x},\vb{z}^i) - w^i) \Big)
= \max_{i \in \{1,\ldots,N\}} (f^*)'(w^i-c(\vb{x},\vb{z}^i)).
\]
For each $i \in \{1,\ldots,N\}$, define $\zeta^i :\X \times \R^N \times D^N  \to \mathbb{R}$ by
\begin{equation}
\label{eq: Catherine Zeta i -Jones}
    \zeta^i (\vb{x},\vb*{w},\vb{z}) = (f^*)'(w^i-c(\vb{x},\vb{z}^i)).
\end{equation}
Then $\zeta^i$
is continuously differentiable because $(f^*)'$ and $c$ are continuously differentiable. 

Fix $(\vb*{w},\vb{z}) \in \R^N \times D^N$ and consider the function $\zeta(\cdot,\vb*{w},\vb{z}) = \max_i \zeta^i(\cdot,\vb*{w},\vb{z})$.  For each $i \in \{1,\ldots,N\}$ the function $\zeta^i(\cdot,\vb*{w},\vb{z})$ is globally Lipschitz on $\X$ because $\zeta^i$ is continuously differentiable and $\X$ is compact. Therefore $\zeta(\cdot,\vb*{w},\vb{z})$ is the pointwise maximum of a finite family of Lipschitz functions, hence Lipschitz. In particular, it is continuous.

Now we consider the functions $\zeta(\vb{x},\cdot,\cdot) = \max_i \zeta^i (\vb{x},\cdot,\cdot)$ indexed by $\vb{x} \in \X$. Let $K \subset \R^N \times D^N$ be compact. 
As above, $\zeta(\vb{x},\cdot,\cdot)|_K$ is the pointwise maximum of a finite family of of Lipschitz functions, hence Lipschitz. Moreover, the Lipschitz constant of $\zeta^i(\vb{x},\cdot,\cdot)|_K$ can be bounded from above by a constant $L^i(K)$ independent of $\vb{x}$ since $\partial \zeta^i/\partial \vb*{w}$, $\partial \zeta^i/\partial \vb{z}$ are continuous, and $\X$ and $K$ are compact. Therefore
$\sup_{\vb{x} \in \X} \mathrm{Lip}(\zeta(\vb{x},\cdot,\cdot)|_K) \le \max_i L^i(K) =: L(K)$, as required.

Let $(\vb*{w}_0,\vb{z}_0) \in \R^N \times D^N$ and
define $S(\vb*{w}_0,\vb{z}_0)=\bigcup_{j=1}^N \mathrm{int}(L^j_c(\vb*{w}_0,\vb{z}_0))$.
Since the cost $c$ is twisted,
$\mathcal{L}^3(\X \setminus S(\vb*{w}_0,\vb{z}_0))=0$
by \cite[Proposition 37]{Merigot:2020}.  For all $\vb{x} \in S(\vb*{w}_0,\vb{z}_0)$, we have
\begin{align*}
\frac{\partial \zeta}{\partial w^i}(\vb{x},\vb*{w}_0,\vb{z}_0)
& = 
\begin{cases}
(f^*)''(w^i_0 - c(\vb{x},\vb{z}^i_0)) & \textrm{if } \vb{x} \in \mathrm{int}(L^i_c(\vb*{w}_0,\vb{z}_0)),
\\
0 & \textrm{if } \vb{x} \in \mathrm{int}(L^j_c(\vb*{w}_0,\vb{z}_0)) \textrm{ for } j \ne i,
\end{cases}
\\
\frac{\partial \zeta}{\partial \vb{z}^i}(\vb{x},\vb*{w}_0,\vb{z}_0)
& = 
\begin{cases}
- \nabla_{\vb{y}}c(\vb{x},\vb{z}^i_0) (f^*)''(w^i_0 - c(\vb{x},\vb{z}^i_0)) & \textrm{if } \vb{x} \in \mathrm{int}(L^i_c(\vb*{w}_0,\vb{z}_0)),
\\
0 & \textrm{if } \vb{x} \in \mathrm{int}(L^j_c(\vb*{w}_0,\vb{z}_0)) \textrm{ for } j \ne i.
\end{cases}
\end{align*}
(Note that the partial derivatives do not exist if $\vb{x} \in L^i_c(\vb*{w}_0,\vb{z}_0) \cap L^j_c(\vb*{w}_0,\vb{z}_0)$ for $j \ne i$.)
In particular, for all $\vb{x} \in S(\vb*{w}_0,\vb{z}_0)$, the partial derivatives 
$\partial \zeta/\partial \vb*{w}(\vb{x},\cdot,\cdot)$, $\partial \zeta/\partial \vb{z}(\vb{x},\cdot,\cdot)$ exist at $(\vb*{w}_0,\vb{z}_0)$. 
Now we prove that the partial derivatives are continuous. 
Let $\vb{x}\in S(\vb*{w}_0,\vb{z}_0)$. Then there exists $i\in\{1,\ldots,N\}$ such that $\vb{x}\in \mathrm{int}(L^i_c(\vb*{w}_0,\vb{z}_0))$, so
\[
c(\vb{x},\vb{z}^i_0) - w^i_0
< c(\vb{x},\vb{z}^j_0) - w^j_0
\quad \forall \; j \in  \{1,\ldots,N\} \setminus \{ i \}.
\]
Since the inequality is strict and $c$ is continuous,
$\vb{x} \in \mathrm{int}(L^i_c(\vb*{w},\vb{z}))$
for all $(\vb*{w},\vb{z})$ sufficiently close to
$(\vb*{w}_0,\vb{z}_0)$, and
\begin{multline*}
    \frac{\partial \zeta}{\partial w^i}(\vb{x},\vb*{w},\vb{z})
- 
\frac{\partial \zeta}{\partial w^i}(\vb{x},\vb*{w}_0,\vb{z}_0)
= 
\\
\begin{cases}
(f^*)''(w^i - c(\vb{x},\vb{z}^i))
-
(f^*)''(w^i_{0} - c(\vb{x},\vb{z}^i_{0}))
& \textrm{if } \vb{x} \in \mathrm{int}(L^i_c(\vb*{w}_0,\vb{z}_0)),
\\
0 & \textrm{if } \vb{x} \in \mathrm{int}(L^j_c(\vb*{w}_0,\vb{z}_0)) \textrm{ for } j \ne i.
\end{cases}
\end{multline*}
Since $(f^*)''$ and $c$ are continuous, it follows that
$\partial \zeta/\partial \vb*{w}(\vb{x},\cdot,\cdot)$ is continuous at $(\vb*{w}_0,\vb{z}_0)$. Similarly, it can be shown that
$\partial \zeta/\partial \vb{z}(\vb{x},\cdot,\cdot)$ is continuous at $(\vb*{w}_0,\vb{z}_0)$.
This completes the proof. 
\end{proof}

\begin{prop}[Regularity of $\G$] 
\label{lem:gderivatvies}
Let $\vb{z} \in D^N$ and $\vb*{w} \in \R^N$. Then $\G(\cdot,\vb{z}) \in \C^1(\R^N)$, $\G(\vb*{w},\cdot)\in\C^1(D^N)$, and
\begin{align}
\label{eq:dGdw}
\frac{\partial \G}{\partial w^i} (\vb*{w},\vb{z}) 
& =
m^i-\int_{L^i_c(\vb*{w},\vb{z})}(f^*)'(w^i-c(\vb{x},\vb{z}^i)) \,\dd \vb{x}, 
\\
\label{eq:dGdz}
\frac{\partial\G}{\partial\vb{z}^i} (\vb*{w}, \vb{z}) 
& = \int_{\L}\nabla_{\vb{y}}c(\vb{x},\vb{z}^i) \, (f^*)'(w^i-c(\vb{x},\vb{z}^i)) \, \dd \vb{x}, 
\end{align}
for all $i \in \{1,\ldots,N\}$. 
Moreover, if $\vb{z} \in D^N$ satisfies $z^i_3 \ne z^j_3$ for all $i,j \in \{1,\ldots,N\}$, $i \ne j$, and 
$\vb*{w} \in \R^N$ satisfies 
$\mathcal{L}^3(L^i_c(\vb*{w},\vb{z}))>0$ for all $i \in \{1,\ldots,N\}$, 
then 
$\G$ is
twice continuously differentiable 
at 
$(\vb*{w},\vb{z})$, and the second-order partial derivatives of $\G$ with respect to $\vb*{w}$ are
\begin{align}
\begin{dcases}
\label{eqn:dGdwdw}
\frac{\partial^2 \G}{\partial w^i \partial w^j} (\vb*{w},\vb{z}) 
& = 
\int_{L^i_c(\vb*{w},\vb{z}) \cap L^j_c(\vb*{w},\vb{z})} 
\frac{(f^*)'(w^i - c(\vb{x},\vb{z}^i))}{\| \nabla_{\vb{x}} c(\vb{x},\vb{z}^i) -  \nabla_{\vb{x}} c(\vb{x},\vb{z}^j) \|} \,\dd \mathcal{H}^2(\vb{x}) \quad \text{for}\; j\neq i,
\\
\frac{\partial^2 \G}{\partial w^i \partial w^i} (\vb*{w},\vb{z}) 
& = 
-\int_{L^i_c(\vb*{w},\vb{z})} (f^*)'' (w^i-c (\vb{x},\vb{z}^i)) \, \dd \vb{x} 
-\sum_{j \neq i} \frac{\partial^2 \G}{\partial w^i \partial w^j} (\vb*{w},\vb{z}).    
\end{dcases}
\end{align}
If in addition
\begin{equation}\label{eqn:mass_pos}
\int_{L^i_c(\vb*{w},\vb{z})} (f^*)' (w^i-c (\vb{x},\vb{z}^i)) \, \dd \vb{x} > 0
\end{equation}
for all $i \in \{1,\ldots,N\}$, 
then the Hessian matrix $D^2_{\vb*{w} \vb*{w}} \G (\vb*{w},\vb{z})$ is negative definite.
\end{prop}

\begin{proof}
Fix $\vb{z}_0\in D^N$. We begin by showing that $\G(\cdot,\vb{z}_0)\in \C^1(\R^N)$. Let $\vb*{w}_0\in \R^N$ and, as in the proof of Lemma \ref{lem:fregularity}, define the set
\[
S(\vb*{w}_0,\vb{z}_0)\coloneqq \bigcup_{j=1}^N \mathrm{int}(L^j_c(\vb*{w}_0,\vb{z}_0)).
\]
For any $\vb{x}\in S(\vb*{w}_0,\vb{z}_0)$ there exists $i\in \{1,\ldots,N\}$ and $\varepsilon>0$ such that, for all $\vb*{w}\in B_\varepsilon(\vb*{w}_0)$,
\[
c(\vb{x},\vb{z}_0^i) - w^i< c(\vb{x},\vb{z}_0^j) - w^j\quad \forall \, j\in\{1,\ldots,N\}.
\]
In particular, 
$\vb{x} \in L^i_c(\vb*{w},\vb{z}_0)$ for all $\vb*{w}\in B_\varepsilon(\vb*{w}_0)$, i.e.,
\[
-\vb*{w}^c(\vb{x};\vb{z}) = w^i - c(\vb{x},\vb{z}_0^i) \quad \forall \,\vb*{w}\in B_\varepsilon(\vb*{w}_0).
\]
It follows that
\begin{equation}
\label{eq:Dwf*}
\forall \; \vb{x} \in \mathrm{int}(L_c^i(\vb*{w}_0,\vb{z}_0)), \quad 
\frac{\partial}{\partial w^j}f^*(-\vb*{w}^c(\vb{x};\vb{z}))
\Big|_{\vb*{w}=\vb*{w}_0} = 
\begin{cases}
(f^*)^\prime(w^i_0 - c(\vb{x},\vb{z}_0^i)) & \text{if } j=i,
    \\
    0 & \text{otherwise.}
\end{cases}
\end{equation}
By continuity of $(f^*)^\prime$,
\begin{equation}\label{eqn:deriv_cts}
\lim_{\vb*{w}\to\vb*{w}_0} (f^*)^\prime(w^j - c(\vb{x},\vb{z}_0^j))\mathds{1}_{L^j_c(\vb*{w},\vb{z}_0)}(\vb{x}) 
= (f^*)^\prime(w^j_0 - c(\vb{x},\vb{z}_0^j))\mathds{1}_{L^j_c(\vb*{w}_0,\vb{z}_0)}(\vb{x})
\end{equation}
for all $j\in\{1,\ldots,N\}$.

As in the proof of Lemma \ref{lem:fregularity}, the map $\vb*{w} \mapsto f^*(-\vb*{w}^c(\vb{x};\vb{z}))$ is locally Lipschitz continuous for all $\vb{x} \in \X$, and its Lipschitz constant can be bounded uniformly, independently of $\vb{x} \in \X$. Moreover, for all $\vb{x} \in S(\vb*{w}_0,\vb{z}_0)$ this map is differentiable at $\vb*{w}_0$ with derivative given by \eqref{eq:Dwf*}.
Then, since $\X$ is compact and $\X\setminus S(\vb*{w}_0,\vb{z}_0)$ is Lebesgue negligible, by the Dominated Convergence Theorem
\begin{align*}
\frac{\partial \G}{\partial w^j}(\vb*{w}_0,\vb{z}_0) 
&= \frac{\partial}{\partial w^j}\left.\left(\sum_{i=1}^N m^i w^i - \int_{\X} f^*(-\vb*{w}^c(\vb{x};\vb{z}))\, \dd \vb{x} \right)\right\vert_{\vb*{w}=\vb*{w}_0}\\
&= m^j - \int_{\X} \frac{\partial}{\partial w^j}f^*(-\vb*{w}^c(\vb{x};\vb{z}))
\Big|_{\vb*{w}=\vb*{w}_0}\, \dd \vb{x} \\
&= m^j - \int_{L^j_c(\vb*{w}_0,\vb{z}_0)}(f^*)^\prime(w^j_0 - c(\vb{x},\vb{z}_0^j))\,\dd\vb{x}.
\end{align*}
This establishes \eqref{eq:dGdw}. Similarly, using \eqref{eqn:deriv_cts} and applying the Dominated Convergence Theorem proves continuity of the partial derivatives of $\G(\cdot,\vb{z}_0)$ at $\vb*{w}_0$. Since $\vb*{w}_0$ was arbitrary, it follows that $\G( \cdot ,\vb{z}_0)\in \C^1(\R^N)$.
Analogously, it can be proved that \eqref{eq:dGdz} holds and that
$\G(\vb*{w}_0,\cdot)\in \C^1(D^N)$ for all $\vb*{w}_0 \in \R^N$; we omit the proof for brevity. 

Now suppose that $\vb{z} \in D^N$ satisfies $z^i_3 \ne z^j_3$ for all $i,j \in \{1,\ldots,N\}$, $i \ne j$, and 
$\vb*{w} \in \R^N$ satisfies 
$\mathcal{L}^3(L^i_c(\vb*{w},\vb{z}))>0$ for all $i \in \{1,\ldots,N\}$. By Lemmas \ref{lem:dualdifferentiability} and \ref{lem:fregularity}, and 
equations \eqref{eq:dGdw} and  \eqref{eq:dGdz},
$\G$ is twice continuously differentiable in a neighbourhood of $(\vb*{w}, \vb{z})$ and its second-order partial derivatives with respect to $\vb*{w}$
are given by \eqref{eqn:dGdwdw}.

In addition, suppose that the positivity constraint \eqref{eqn:mass_pos} holds. Then
\[
\int_{L^i_c(\vb*{w},\vb{z})} (f^*)'' (w^i-c (\vb{x},\vb{z}^i)) \, \dd \vb{x} > 0
\]
for all $i \in \{1,\ldots,N\}$.
Using this inequality and the 
non-negativity of $(f^*)^\prime$ and $(f^*)^{\prime\prime}$ 
gives
\[
\left\vert\frac{\partial^2 \G}{\partial w^i \partial w^i} (\vb*{w},\vb{z}) \right\vert
-
\sum_{j \neq i} \left\vert \frac{\partial^2 \G}{\partial w^i \partial w^j} (\vb*{w},\vb{z})\right\vert
= \int_{L^i_c(\vb*{w},\vb{z})} (f^*)'' (w^i-c (\vb{x},\vb{z}^i)) \, \dd \vb{x} 
> 0.
\]
Altogether, the Hessian $D^2_{\vb*{w} \vb*{w}}\G(\vb*{w},\vb{z})$ is strictly diagonally dominant, symmetric, and has negative diagonal entries, which implies that it is negative definite, as claimed.
\end{proof}

Now we prove the duality theorem, Theorem \ref{thm:dualitytheorem}.

\begin{proof}[Proof of Theorem~\ref{thm:dualitytheorem}]
By Lemma \ref{lem:E unique minimiser}, $E(\cdot,\alpha^N)$ is strictly convex and has a unique minimiser over $\Pac(\X)$.

Now we show that the maximum of $\G(\cdot,\vb{z})$ over $\R^N$ is attained. Since $\G(\cdot,\vb{z})$ is continuous, it suffices to show that $\lim_{\| \vb*{w} \| \to \infty} \G(\vb*{w},\vb{z})=-\infty$.
Since $\X$ is compact and $c$ is continuous, there exists a constant $M(\vb{z})>0$ such that
\begin{align*}
\max_{j \in \{ 1,\ldots,N \}}\max_{\vb{x}\in\X}|c(\vb{x},\vb{z}^j)|\leq M(\vb{z}).
\end{align*}
Since $f^*$ is non-decreasing, 
\[
f^*(-\vb*{w}^c(\vb{x};\vb{z}))
= f^* \left( 
\max_{j \in \{ 1,\ldots,N \}} 
\left( w^j - c(\vb{x},\vb{z}^j) \right)
\right)
\ge 
f^*
\left( 
\max_j w^j - M(\vb{z})
\right).
\]
For each $\vb*{w} \in \R^N$, write $\vb*{w} = \vb*{w}_+ - \vb*{w}_-$, where $\vb*{w}_+,\vb*{w}_-$ are the positive and negative parts of $\vb*{w}$, namely, $w_+^i = \max \{w^i,0 \}$, $w_-^i = - \min \{w^i,0\}$ for all $i \in \{1,\ldots,N\}$.
Since $\sum_{i=1}^N m^i = 1$,
\begin{align*}
\G(\vb*{w})
& \le 
\sum_{i=1}^N m^i (w_+^i-w_-^i) - \mathcal{L}^3(\X) 
f^*
\left( 
\max_j w^j - M(\vb{z})
\right)
\\
&
\le \| \vb*{w}_+\|_\infty - \min_i m^i \, \| \vb*{w}_-\|_\infty
- \mathcal{L}^3(\X) 
f^*
\left( 
\max_j w^j - M(\vb{z})
\right)
\\
& = 
- \min_i m^i \, \| \vb*{w}_-\|_\infty
+ 
\begin{cases}
0 & \textrm{if } \vb*{w}_+=0,
\\
\| \vb*{w}_+\|_\infty - \mathcal{L}^3(\X)
f^*
\left( 
\| \vb*{w}_+ \|_\infty - M(\vb{z})
\right)
& \textrm{if } \vb*{w}_+ \ne 0.
\end{cases}
\end{align*}
The first term tends to $-\infty$ as $\| \vb*{w}_-\|_\infty \to \infty$, and the second term tends to $-\infty$
as $\| \vb*{w}_+\|_\infty \to \infty$ because $f^*$ grows superlinearly at $+\infty$. Therefore
$\lim_{\| \vb*{w} \| \to \infty} \G(\vb*{w},\vb{z})=-\infty$, as claimed.

Next we show that $\G(\cdot,\vb{z})$ is concave. Since $f^*$ is non-decreasing, 
\[
\G(\vb*{w},\vb{z}) 
= \sum_{i=1}^N m^iw^i - \int_{\X} f^*(-\vb*{w}^c(\vb{x};\vb{z})) \, \dd \vb{x} 
= \sum_{i=1}^N m^iw^i
- \int_{\X} \max_{j \in \{1,\ldots,N\}} f^*(w^j -   c(\vb{x},\vb{z}^j)) \, \dd \vb{x} .
\]
Note that the map $\vb*{w} \mapsto \max_{j \in \{1,\ldots,N\}} f^*(w^j -   c(\vb{x},\vb{z}^j))$ is convex because 
$f^*$ is convex and the pointwise maximum over a family of convex functions is convex. Therefore it follows by linearity of integration that $\G(\cdot,\vb{z})$ is concave.

Now we prove that the maximiser of $\G(\cdot,\vb{z})$ is unique. 
Let $\vb*{w}$ be a maximiser of $\G(\cdot,\vb{z})$. 
Then $\nabla_{\vb*{w}} \G (\vb*{w},\vb{z}) =0$ and hence \eqref{eqn:mass_pos} holds (by \eqref{eq:dGdw} and because $m^i>0$ for all $i$). Then by Proposition \ref{lem:gderivatvies} the Hessian of $\G(\cdot,\vb{z})$ at $\vb*{w}$ is negative definite, so $\vb*{w}$ is an isolated maximiser. Since $\G(\cdot,\vb{z})$ is concave, its set of maximisers must be convex, hence connected. Therefore, $\vb*{w}$ is the unique maximiser of $\G(\cdot,\vb{z})$.

Next we prove weak duality.
Recall the Fenchel-Young inequality: 
\begin{gather}
\label{eq:FY-1}    
ab \leq f(a) + f^*(b) \quad \forall \; a,b \in \R,
\\
\label{eq:FY-2}
ab = f(a) + f^*(b) 
\quad \Longleftrightarrow \quad 
a\in \partial f^*(b)=\{ (f^*)'(b) \}
\quad \Longleftrightarrow \quad 
b\in\partial f(a).
\end{gather}
By the Kantorovich Duality Theorem  (see \eqref{eq:optdiscT}) and the Fenchel-Young inequality \eqref{eq:FY-1}, for all $\sigma\in \Pac(\X)$, $\vb*{w} \in \R^N$,
\begin{align}
\nonumber
    \E(\sigma,\alpha^N)
    & = \T(\sigma,\alpha^N) + \int_\X f(\sigma(\vb{x}))\,\dd\vb{x}
    \\
    \nonumber
    & = \max_{\widetilde{\vb*{w}}\in\R^N} \left(\sum_{i=1}^N m^i \widetilde{w}^i + \int_{\X} \widetilde{\vb*{w}}^c(\vb{x};\vb{z})\sigma(\vb{x}) \,  \dd \vb{x}\right)
    + \int_\X f(\sigma(\vb{x}))\,\dd\vb{x}
    \\
    \nonumber
    & \ge \sum_{i=1}^N m^i w^i + \int_{\X} \big( \vb*{w}^c(\vb{x};\vb{z})\sigma(\vb{x}) + f(\sigma(\vb{x}))\big)\, \dd \vb{x}
    \\
    \label{eq: Keir Starmer}
    & \geq \sum_{i=1}^N m^i {w}^i - \int_{\X} f^*(-{\vb*{w}}^c(\vb{x};\vb{z}))\, \dd \vb{x} = \G(\vb*{w},\vb{z}).
\end{align}
In particular, we have the weak duality 
\begin{equation}
\label{eq:weak duality}
\min_{\sigma \in \Pac(\X)} E(\sigma,\alpha^N) \ge \max_{\vb*{w} \in \R^N} \G(\vb*{w},\vb{z}).
\end{equation}
It follows that 
if $\E(\sigma,\alpha^N) = \G(\vb*{w},\vb{z})$, then
$\vb*{w}$ is a maximiser of $\G(\cdot,\vb{z})$ and $\sigma$ is a minimiser of $E(\cdot,\alpha^N)$.

Next we show that
\eqref{eq:gMassConstraint} and \eqref{eq:density} are 
sufficient for optimality.
Suppose that $(\sigma, \vb*{w})\in\Pac(\X)\times\R^N$ satisfy
\eqref{eq:gMassConstraint} and \eqref{eq:density}. Then
\[
\int_{L^i_c(\vb*{w},\vb{z})} \sigma(\vb{x}) \, \dd \vb{x} = m^i \quad \forall \; i \in \{1,\ldots,N\}.
\]
By Theorem \ref{thm: well known SDOT}, $\vb*{w}$ is an optimal Kantorovich potential for transporting $\sigma$ to $\alpha^N$ and
\begin{equation}
\label{eq: Elon Musk}    
\T(\sigma,\alpha^N)
= \int_{\X} 
\vb*{w}
^c(\vb{x};\vb{z}) \, \dd  \sigma (\vb{x}) + \sum_{i=1}^N m^ i w^i.
\end{equation}
By \eqref{eq:density},
\[
\sigma(\vb{x})
= (f^*)'(-\vb*{w}^c(\vb{x};\vb{z}))
\in \partial f^*(-\vb*{w}^c(\vb{x};\vb{z})) \quad \forall \; \vb{x}\in \X.
\]
Then by the Fenchel-Young inequality \eqref{eq:FY-2},
\begin{equation}\label{eqn:LF_equality}
-\vb*{w}^c(\vb{x};\vb{z}) \,\sigma(\vb{x}) = 
f(\sigma(\vb{x}))
+ f^*(-\vb*{w}^c(\vb{x};\vb{z}))   \quad \forall \; \vb{x}\in \X.    
\end{equation}
Combining \eqref{eq: Elon Musk} and \eqref{eqn:LF_equality} gives
\[
E(\sigma,\alpha^N)
= \sum_{i=1}^N m^i w^i - \int_{\X} f^*(-\vb*{w}^c(\vb{x};\vb{z}))\, \dd \vb{x} = \G(\vb*{w},\vb{z}).
\]
Therefore, by \eqref{eq:weak duality},  $\vb*{w}$ is a maximiser of $\G(\cdot,\vb{z})$ and $\sigma$ is a minimiser of $E(\cdot,\alpha^N)$.

We now show that \eqref{eq:gMassConstraint} and \eqref{eq:density} are necessary for optimality. Let
$\sigma \in\Pac(\X)$ be the unique minimiser of $E(\cdot,\alpha^N)$ and let $\vb*{w} \in\R^N$ be the unique maximiser of $\G(\cdot,\vb{z})$.
Then 
$\nabla_{\vb*{w}}\G(\vb*{w},\vb{z})=0$ and hence, by \eqref{eq:dGdw},  $\vb*{w}$ satisfies \eqref{eq:gMassConstraint}. Define $\widetilde{\sigma}: \X \to \mathbb{R}$ by
\begin{equation}
\label{eq:sig def}    
\widetilde{\sigma}(\vb{x}) 
 = (f^*)'(-\vb*{w}^c(\vb{x};\vb{z})).
\end{equation}
The pair
$(\widetilde{\sigma}, \vb*{w})\in\Pac(\X)\times\R^N$ satisfies
\eqref{eq:gMassConstraint} and \eqref{eq:density} by construction. It remains to show that $\widetilde{\sigma}$ is a probability measure and
\begin{equation}\label{eqn:sigma_identification}
    \widetilde{\sigma} = \sigma \quad \mathcal{L}^3\text{-almost everywhere.}
\end{equation}
Since $\sum_{i=1}^N m^i = 1$ and \eqref{eq:gMassConstraint} holds, it follows that $\widetilde{\sigma}$ is a probability measure. Evaluating $E(\widetilde{\sigma},\alpha^N)$ as in the previous paragraph reveals that $E(\widetilde{\sigma},\alpha^N)=\G(\vb*{w},\vb{z})$. As argued above, this holds if and only if $\widetilde{\sigma}$ and $\vb*{w}$ are optimal. But since the minimiser of $E(\cdot,\alpha^N)$ is unique, \eqref{eqn:sigma_identification} holds.

To conclude, we have shown that \eqref{eq:gMassConstraint} and \eqref{eq:density} are necessary and sufficient for optimality. The argument above for sufficiency implies that $E(\sigma,\alpha^N)=\G(\vb*{w},\vb{z})$ for the unique pair of optimisers $(\sigma,\vb*{w})$, so strong duality \eqref{eq:duality} holds.
\end{proof}

An immediate consequence of Theorem \ref{thm:dualitytheorem} and  Lemma \ref{lem:fregularity} is the following. 

\begin{cor}[Regularity of the optimal source measure for discrete target measures]
\label{cor: reg source}
Let $\alpha^N \in \Pn^N(\Y)$ be a discrete measure and $\sigma_*[\alpha^N] \in \Pac(\X)$ be the optimal source measure given in Definition \ref{def:optimal_source_general}. Then 
$\sigma_*[\alpha^N]$ is Lipschitz continuous.
\end{cor}

\section{Discrete solutions of the compressible SG equations}
\label{sec:discrete_solutions}

The goal of this section is to prove 
Theorem \ref{thm:discretesolutions}, namely, that
for well-prepared initial data, discrete solutions of the compressible semi-geostrophic equations exist, are unique, and are energy conserving. 


We first define the centroid map $\vb{C}$ that appears in the ODE \eqref{eq:ODESystem} and show that it is continuously differentiable.

\begin{defn}[Optimal weight map]
\label{def:opt_w}
Given $\vb{m} \in \Delta^N$, define the map $\vb*{w}_*:D_0^N\to\R^N$ by
\begin{equation}
\label{eq:dualenergy}
    \vb*{w}_*(\vb{z}):=\argmax_{\vb*{w}\in\R^N} \G(\vb*{w},\vb{z}).
\end{equation}
\end{defn}

Note that $\vb*{w}_*$ is well defined by Theorem \ref{thm:dualitytheorem}. 

\begin{defn}[Centroid map]
\label{def:centroid}
Given $\vb{m} \in \Delta^N$, define the centroid map $\vb{C}:D_0^N\to(\R^3)^N$ by
 $\vb{C}(\vb{z})
 \coloneqq (\vb{C}^1(\vb{z}),\ldots,\vb{C}^N(\vb{z}))$, where
\begin{equation}
\label{eq:centroid0}
     \vb{C}^i(\vb{z})\coloneqq\frac{1}{m^i}\int_{\Lopt}\vb{x} \, \dd \sigma_*[\alpha^N](\vb{x})
     ,
\end{equation}
where
\[
\alpha^N =\sum_{i=1}^Nm^i\delta_{\vb{z}^i}.
\]
That is, $\vb{C}^i(\vb{z})$ is the centroid (or barycentre) of the set of points $T_{\alpha^N}^{-1}(\{ \vb{z}^i \})$  transported to $\vb{z}^i$ in the optimal transport from $\sigma_*[\alpha^N]$ to $\alpha^N$ for the cost $c$.
\end{defn}

To prove the existence of solutions of the ODE \eqref{eq:ODESystem}, we show that the transport velocity $W:D_0^N \to \mathbb{R}^{3N}$ is continuously differentiable, where $W$ is defined by
\[
W(\vb{z})=J^N\qty(\vb{z}-\vb{C}(\vb{z})),
\]
where $J^N$ is the block diagonal matrix $ J^N:=\mathrm{diag}(J,\ldots,J) \in\R^{3N\times 3N}$.

\begin{lemma}[Regularity of $\vb*{w}_*$]
\label{lem:C1weight}
The
optimal weight map $\vb*{w}_*$ is continuously differentiable.
\end{lemma}

\begin{proof}
By definition of $\vb*{w}_*$, $\nabla_{\vb*{w}} \G(\vb*{w}_*(\vb{z}),\vb{z})=0 $ for all $\vb{z}\in D^N$. In order to show that $\vb*{w}_*$ is continuously differentiable, we will apply the Implicit Function Theorem to the function $\nabla_{\vb*{w}} \G$.
Fix any point $\vb{z}_0 \in D_0^N$ and define $\vb*{w}_0 = \vb*{w}_*(\vb{z}_0)$. By 
Proposition \ref{lem:gderivatvies}
and equation \eqref{eq:gMassConstraint}, $\G$ is twice continuously differentiable in a neighbourhood of $(\vb*{w}_0,\vb{z}_0)$ and 
$D^2_{\vb*{w} \vb*{w}} \G (\vb*{w}_0,\vb{z}_0)$ is negative definite, hence invertible.
By the Implicit Function Theorem, there exists an open set $U \subset D_0^{N}$ containing $\vb{z}_0$ and a unique, continuously differentiable function $\widetilde{\vb*{w}} : U \to \R^N$ such that $\widetilde{\vb*{w}}(\vb{z}_0)=\vb*{w}_0$ and $\nabla_{\vb*{w}}\G\qty(\widetilde{\vb*{w}}(\vb{z}),\vb{z})=0$ for all $\vb{z}\in U$. By uniqueness, $\widetilde{\vb*{w}}$ coincides with $\vb*{w}_*$ on $U$. Therefore $\vb*{w}_*$ is continuously differentiable in a neighbourhood of $\vb{z}_0$, as required.
\end{proof}

\begin{lemma}[Regularity of $\vb{C}$]\label{lem:C1centroid}
The centroid map $\vb{C}$
is continuously differentiable.
\end{lemma}

\begin{proof}
Recalling the expression for $\sigma_*[\alpha^N]$ derived in Theorem \ref{thm:dualitytheorem}, we have
\[
\vb{C}^i(\vb{z})
     =
     \frac{1}{m^i}\int_{\Lopt} \vb{x} \,
     (f^*)'(-(\vb*{w}_*(\vb{z}))^c(\vb{x};\vb{z}))
     \, \dd \vb{x}.
\]
Define $\vb{\Psi}=(\Psi^1,\ldots,\Psi^N): \R^N \times D_0^N \to (\R^3)^N$
by
\[
\Psi^i(\vb*{w},\vb{z}) = \frac{1}{m^i} \int_{L^i_c(\vb*{w},\vb{z})} \vb{x} \, (f^*)'(-\vb*{w}^c(\vb{x};\vb{z})) \, \dd \vb{x}.
\]
Then $\vb{C} = \Psi \circ (\vb*{w}_*,\mathrm{id})$.
For $j \in \{1,2,3\}$, define the functions $\zeta_j:\X \times \R^N \times D^N \to \mathbb{R}$ by $\zeta_j(\vb{x},\vb*{w},\vb{z})= x_j \, (f^*)'(-\vb*{w}^c(\vb{x};\vb{z}))$. It follows from Lemma \ref{lem:fregularity} that $\zeta_j$ satisfy the assumptions of Lemma \ref{lem:dualdifferentiability}, so $\Psi$ is continuously differentiable by Lemma \ref{lem:dualdifferentiability}. Therefore $\vb{C}$ is continuously differentiable by Lemma \ref{lem:C1weight} and the Chain Rule.
\end{proof}


We now show that solutions of the ODE \eqref{eq:ODESystem} correspond to weak solutions of the PDE \eqref{eq:PDE} with discrete initial data.The following is essentially a consequence of 
\cite[Lemma 4.2]{Bourne:2022}.

\begin{prop}
\label{lem:Correspondence} 
Let $\overline{z} \in D^N$ and let $\overline{\alpha}^N\in\Pn^N(\Y)$ be the corresponding discrete initial measure given by
\begin{equation}
    \overline{\alpha}^N=\sum_{i=1}^Nm^i\delta_{\overline{\vb{z}}^i}.
\end{equation}
Fix a final time $\tau>0$ and let $\vb{z}\in\C^1([0,\tau];\Y^N)$.
Define $\alpha^N:[0,\tau]\to\Pn^N(\Y)$ by
\begin{equation}\label{eq:PDESolution}
    \alpha^N_t=\sum_{i=1}^Nm^i\delta_{\vb{z}^i(t)}\quad\forall\, t\in\qty[0,\tau].
\end{equation}
Then $\alpha^N$
is a weak solution of the compressible semi-geostrophic equations (in the sense of Definition \ref{def:weak}) with initial data $\alpha_0^N = \overline{\alpha}^N$ if and only if $\vb{z}$ is a solution of the ODE \eqref{eq:ODESystem} with initial data $\vb{z}(0)=\overline{\vb{z}}$.
\end{prop}

\begin{proof}
Suppose that $\vb{z}\in C^1(\qty[0,\tau];\Y^{N})$ is a solution of the ODE \eqref{eq:ODESystem} and let $\alpha^N$ be given by \eqref{eq:PDESolution}.  We show that $\alpha^N$ satisfies \eqref{eq:WeakFormPDE}. 
Let $\varphi\in\D(\Y\times\R)$. Recall that we use the notation
$\varphi_t(y) = \varphi(y,t)$, i.e., the subscript $t$ does not denote differentiation. Then
\begin{align*}
& \int_0^{\tau} \int_{\Y} \left[ \partial_t \varphi_t (\vb{y}) + (J \vb{y}) \cdot \nabla \varphi_t (\vb{y})\right] \, \dd \alpha^N_t (\vb{y}) \, \dd t -\int_0^{\tau} \int_{\X} (J \vb{x} ) \cdot \nabla \varphi_t (T_{\alpha^N_t}(\vb{x})) \, \dd \sigma_*[\alpha^N_t](\vb{x}) \, \dd t 
\\
& = 
\sum_{i=1}^N m^i \int_0^{\tau} \left[ \partial_t \varphi_t (\vb{z}^i(t)) + (J\vb{z}^i(t)) \cdot \nabla \varphi_t (\vb{z}^i(t)) \right] \, \dd t 
\\
& \qquad
- \int_0^{\tau}\sum_{i=1}^N \int_{L^i_c(\vb*{w}_*(\vb{z}(t)),\vb{z}(t))}(J \vb{x}) \cdot \nabla \varphi_t(\vb{z}^i(t)) \,\dd \sigma_*[\alpha^N_t](\vb{x}) \, \dd t 
\\ 
& = 
\sum_{i=1}^N m^i \int_0^{\tau} \left[ \partial_t
\varphi_t ( \vb{z}^i(t)) + J \left( \vb{z}^i(t)-\vb{C}^i(\vb{z}(t)) \right) 
\cdot \nabla \varphi_t (\vb{z}^i(t)) \right]
\, \dd t 
\\
& = 
\sum_{i=1}^N m^i \int_0^{\tau}
\left[ \partial_t \varphi_t (\vb{z}^i(t)) +\dot{\vb{z}}^i (t) \cdot \nabla \varphi_t (\vb{z}^i(t)) \right] \, \dd t 
\\
& = 
\sum_{i=1}^N m^i \int_0^{\tau} 
\frac{d}{dt}
\varphi_t (\vb{z}^i(t)) \,\dd t 
\\
& = 
\sum_{i=1}^N m^i \left( \varphi_{\tau}(\vb{z}^i(\tau)) -\varphi_0 (\vb{z}^i(0)) \right) 
\\
& = \int_{\Y}  \varphi_{\tau} (\vb{y}) \,\dd \alpha^N_{\tau}(\vb{y})-\int_{\Y} \varphi_0 (\vb{y}) \, \dd \overline{\alpha}^N( \vb{y}),
 \end{align*}
as required. From here, the proof follows exactly the proof of Lemma 4.2 given in \cite{Bourne:2022}.
\end{proof}


Next we prove global existence and uniqueness of solutions of the ODE \eqref{eq:ODESystem} using the following \emph{a priori} estimates.

\begin{lemma}[\emph{A priori} estimates]
\label{lem:apriori} 
Let $\vb{z}:[0,\tau] \to D^N$ be continuously differentiable and satisfy the ODE \eqref{eq:ODESystem} on the interval $[0,\tau]$. 
Let $R=R(\X)>0$ satisfy $\X\subset B_R(0)$.
Then, for each $i\in \{ 1, \ldots, N \}$ and $t \in [0,\tau]$,
\begin{align}
\label{eq: a priori 1}
\| \vb{z}^i(t) \| & \leq \| \vb{z}^i(0) \| + \f R \tau,
\\
\label{eq: a priori 2}
\| \dot{\vb{z}}^i(t) \| & \leq \f \left( \| \vb{z}^i(0) \| + R (\f \tau + 1) \right).
\end{align}
\end{lemma}

\begin{proof}
By \eqref{eq:ODESystem}, and since $J$ is skew-symmetric with $\| J \|_2 = \f$,
\[
\frac{d}{dt} \| \vb{z}^i(t) \|
= \frac{\vb{z}^i}{\| \vb{z}^i(t) \|} \cdot
\dot{\vb{z}}^i(t) 
= - \frac{\vb{z}^i(t) \cdot J \, \vb{C}^i(\vb{z}(t))}{\| \vb{z}^i(t) \|}
\le \f \| \vb{C}^i(\vb{z}(t)) \| \le \f R.
\]
Integrating gives \eqref{eq: a priori 1}.
By \eqref{eq:ODESystem} and \eqref{eq: a priori 1},
\[
\| \dot{\vb{z}}^i(t) \|
= \| J (  \vb{z}^i(t) - \vb{C}^i(\vb{z}(t))  ) \| 
\le 
\f \left( \| \vb{z}^i(t) \| + \| \vb{C}^i(\vb{z}(t)) \| \right)
\le \f \left( \| \vb{z}^i(0) \| + \f R \tau + R \right),
\]
as required. 
\end{proof}

\begin{prop}[Existence and uniqueness of solutions of the ODE \eqref{eq:ODESystem}]
\label{prop:ODESolutions}
Let the final time $\tau \in (0,\infty)$ be arbitrary and
let the initial data $\overline{\vb{z}} \in D_0^{N}$.
Then there exists a unique function $\vb{z} \in \C^2([0,\tau];D^N_0)$ such that $\vb{z}$ satisfies the ODE \eqref{eq:ODESystem} on the interval $[0,\tau]$ and $\vb{z}(0) = \overline{\vb{z}}$.
\end{prop}

\begin{proof}
The proof of Proposition \ref{prop:ODESolutions} follows exactly the proof of Proposition 4.4 of \cite{Bourne:2022}. Indeed, local existence and uniqueness follows from the Picard–Lindel\"of Theorem, Lemma \ref{lem:C1centroid}, and the assumption on the initial data $\overline{\vb{z}}$. Global existence then follows from the \emph{a priori} estimates (Lemma~\ref{lem:apriori}) and the fact that the third row of the matrix $J$ is zero, which implies that any solution $\vb{z}$ of the ODE satisfies $z^i_3(t) = z^i_3(0)$ for all $t$ and $i$. In particular, for all $t$, $\vb{z}(t)$ belongs to the subset $D_0^N$ of $D^N$ where the centroid map $\vb{C}$ is continuously differentiable.
\end{proof}

Furthermore, we prove that the discrete solutions conserve energy. 

\begin{prop}[Solutions are energy-conserving]
\label{prop:energy-conserving}
Any solution $\vb{z}\in\C^1([0,\tau];D_0^N)$ of  the ODE \eqref{eq:ODESystem}
is energy-conserving in the sense that
\[
\frac{d}{dt}
E(\sigma_*[\alpha_t^N], \alpha_t^N) = 0
\quad \forall \; t \in [0,\tau],
\]
where
\[    \alpha_t^N = \sum_{i=1}^{N}m^i\delta_{\vb{z}^i(t)}.
\]
\end{prop}

\begin{proof}
By the duality theorem (Theorem \ref{thm:dualitytheorem}),
\begin{align}
\nonumber
\frac{d}{dt}
E(\sigma_*[\alpha_t^N], \alpha_t^N)  
& = \frac{d}{dt} \G (\vb*{w}_*(\vb{z}(t)),\vb{z}(t))
\\
\nonumber
& =
\sum_{i=1}^N 
\left( \frac{\partial \G}{\partial w^i}
(\vb*{w}_*(\vb{z}(t)),\vb{z}(t))
\frac{d}{dt} w^i_*(\vb{z}(t))
+ \frac{\partial \G}{\partial \vb{z}^i}
(\vb*{w}_*(\vb{z}(t)),\vb{z}(t)) \cdot
\dot{\vb{z}}^i(t)
\right)
\\
\label{eq: econ 1}
& =
\sum_{i=1}^N 
\frac{\partial \G}{\partial \vb{z}^i}
(\vb*{w}_*(\vb{z}(t)),\vb{z}(t)) \cdot
\dot{\vb{z}}^i(t)
\end{align}
because
$\vb*{w}_*(\vb{z}(t))$ is the maximiser (and hence a critical point) of $\G(\cdot,\vb{z}(t))$.
For brevity, define
\[
L^i_c(t) := 
L^i_c (\vb*{w}_*(\vb{z}(t)),\vb{z}(t)),
\]
and let 
$\vb{C}^i(t) = \vb{C}^i(\vb{z}(t))$.
By \eqref{eq:density} and \eqref{eq:dGdz},
\begin{align}
\nonumber
\frac{\partial\G}{\partial\vb{z}^i} (\vb*{w}_*(\vb{z}(t)),\vb{z}(t)) 
\cdot \dot{\vb{z}}^i(t)
& = \left(\int_{L^i_c(t)}\nabla_{\vb{y}}c(\vb{x},\vb{z}^i(t))  \, \sigma_*[\alpha_t^N](\vb{x}) \, \dd \vb{x}\right)
\cdot \dot{\vb{z}}^i(t)
\\
\nonumber
& =
\frac{1}{z_3^i(t)}\left(\int_{L^i_c(t)}
\begin{pmatrix}
\f^2(z_1^i(t)-x_1) 
\\ \f^2(z_2^i(t)-x_2) 
\\ -c(\vb{x},\vb{z}^i(t))
\end{pmatrix}
\sigma_*[\alpha_t^N](\vb{x}) \,\dd \vb{x}\right)
 \cdot \dot{\vb{z}}^i(t)
\\
\nonumber
& =
\frac{m^i \f^2}{z_3^i(t)}
(\vb{z}^i - \vb{C}^i(t))\cdot J(\vb{z}^i - \vb{C}^i(t))
\\
\label{eq: econ 2}
& = 0
\end{align}
where the final two equalities hold by \eqref{eq:ODESystem}, skew-symmetry of $J$, and the fact that the third row of $J$ is zero. Combining \eqref{eq: econ 1} and \eqref{eq: econ 2} completes the proof.
\end{proof}

Finally, we prove  
Theorem \ref{thm:discretesolutions}.

\begin{proof}[Proof of Theorem \ref{thm:discretesolutions}]
Let $\overline{\alpha}^N = \sum_{i=1}^N m^i \delta_{\overline{\vb{z}}^i}  \in \Pn^N(\Y)$ be well-prepared in the sense of \eqref{wellprep}.
By Proposition \ref{prop:ODESolutions} the ODE \eqref{eq:ODESystem} has a unique solution $\vb{z} \in \C^2([0,\tau];D_0^N)$ satisfying $\vb{z}(0)=\overline{\vb{z}}$. Define $\alpha^N:[0,\tau] \to \Pn^N(\Y)$ by $\alpha^N_t=\sum_{i=1}^N m^i\delta_{\vb{z}^i(t)}$. By Proposition \ref{lem:Correspondence} $\alpha^N$ is  weak solution of the compressible SG equations (in the sense of Definition \ref{def:weak}) with $\alpha^N_0 = \overline{\alpha}^N$. Furthermore, $\alpha^N$ is energy-conserving by Proposition \ref{prop:energy-conserving}. 
\end{proof}

\section{Existence of weak solutions}
\label{sec:exist}
In this section we prove Theorem \ref{thm:solutions}.
The strategy is to approximate the initial measure $\overline{\alpha} \in \Pc(\Y)$ by a sequence of discrete measures $\overline{\alpha}^N \in \Pn^N(\Y)$, apply  Theorem~\ref{thm:discretesolutions} to obtain a weak solution  of \eqref{eq:PDE} for each initial data $\overline{\alpha}^N$, and then pass to the limit $N \to \infty$ to obtain a weak solution of \eqref{eq:PDE} with initial data $\overline{\alpha}$. This is the same proof strategy as for the \emph{incompressible} system given in \cite{Bourne:2022}. Indeed, some of the preliminary lemmas, namely Lemmas \ref{5.1} and \ref{lem:compact2}, are identical and are stated here without proof. 
The first lemma asserts that any probability measure $\mu \in \Pc(\Y)$ can be approximated in the $W_1$ metric by a discrete measure $\mu^N \in \Pn^N(\Y)$ with seeds in distinct horizontal planes.

\begin{lemma}[Quantization, \protect{\cite[Lemma 5.1]{Bourne:2022}}]\label{5.1}
Let $\mu\in \Pc(\Y)$. There exists a compact set 
$K\subset\Y$ and
a sequence of well-prepared discrete probability measures $\mu^N\in\Pn^N(\Y)$ (in the sense of \eqref{wellprep}) such that 
\begin{gather*}
\mathrm{spt}(\mu^N) \subset K \quad \forall \; N \in \N,
\qquad     \lim_{N\to\infty}W_1(\mu^N,\mu)=0.
\end{gather*}
\end{lemma}

\begin{lemma}[Compactness, \protect{\cite[Lemma 5.2]{Bourne:2022}}] \label{lem:compact2}
Let $\overline{\alpha} \in \Pc(\Y)$ and let $\overline{\alpha}^N \in \Pn^N(\Y)$ given by Lemma \ref{5.1} be a sequence of well-prepared discrete discrete probability measures converging to $\overline{\alpha}$ in $W_1$. 
Fix $\tau \in (0,\infty)$.
Let $\alpha^N\in \C^{0,1}([0,\tau];\Pn^N(\Y))$ given by Theorem \ref{thm:discretesolutions} be a discrete weak solution of \eqref{eq:PDE} with initial data $\overline{\alpha}^N$.
Then the sequence $(\alpha^N)_{N\in\N}$ has a uniformly convergent subsequence in $\C\qty([0,\tau];\Pc(\Y))$. To be precise, there exists a subsequence (which we do not relabel) and a Lipschitz map $\alpha \in \C^{0,1} ([0,\tau];\Pc(\Y))$ such that 
\begin{equation}
\label{eq: uniform convergence}    
\lim_{N \to \infty} \,\sup_{t \in [0,\tau]} W_1 \big( \alpha_t^N,\alpha_t \big) = 0.
\end{equation}
Moreover, there exists a compact set $K \subset \Y$ such that
$\mathrm{spt} ( \alpha_t^N ) \subset K$ and $\mathrm{spt} (\alpha_t ) \subset K$
for all $t\in[0,\tau]$ and $N\in\N$.
\end{lemma}

The final task is to show that $\alpha$ in Lemma \ref{lem:compact2} satisfies the PDE \eqref{eq:PDE}.
One of the ingredients of the proof is the following.

\begin{lemma}[Uniform convergence of source measures]
\label{lem: uniform convergence}
Given a compact set $K \subset \Y$ and a sequence $\beta^N\in\Pn^N(\Y) \cap \Pn(K)$ converging weakly to $\beta\in\Pn(K)$, the sequence of optimal source  measures $\sigma_*[\beta^N]$ converges uniformly to $\sigma_*[\beta]$.
\end{lemma}

\begin{proof}
Recall from Corollary \ref{cor: reg source} that 
$\sigma_*[\beta^N]$ is Lipschitz continuous for all $N \in \N$. Moreover, it follows from the proof of Lemma \ref{lem:fregularity} that its $W^{1,\infty}$-norm can be bounded independently of $N$ (by noting that $\max_{\vb{z} \in D^N \cap K^N} \| \zeta^i(\cdot,\vb*{w}_*(\vb{z}),\vb{z}) \|_{\C^1(\X)}$ is independent of $N$, where $\zeta^i$ was defined in equation \eqref{eq: Catherine Zeta i -Jones}).
Therefore the sequence $\sigma_* [\beta^N]$ is uniformly bounded and  equicontinuous.
By the Ascoli-Arzel\'a Theorem, there exists $\widetilde{\sigma}\in \C(\X) \cap \Pac(\X)$ such that $\sigma_*[\beta^N]$ converges uniformly in $\X$ to $\widetilde{\sigma}$ (up to a subsequence 
that we do not relabel).

Next we show that 
$\widetilde{\sigma} = \sigma_*[\beta]$.
The continuity of $\T$ (see \cite[Theorem 1.51]{Santambrogio:2015}) and the uniform convergence of $\sigma_*[\beta^N]$ imply that
\begin{equation*}
E\qty(\sigma_*[\beta^N], \beta^N) \to E(\widetilde{\sigma}, \beta) 
\qq{and}E(\sigma_*[\beta], \beta^N) \to E(\sigma_*[\beta], \beta) 
\qq{as} 
N \to \infty.
\end{equation*}
Assume for a contradiction that $E(\widetilde{\sigma},\beta)>E(\sigma_*[\beta],\beta)$. Then $E(\sigma_*[\beta^N],\beta^N) > E(\sigma_*[\beta],\beta^N)$ for $N$ sufficiently large, which contradicts the fact that $\sigma_*[\beta^N]$ is the global minimiser of $E(\cdot,\beta^N)$. Hence $E(\widetilde{\sigma},\beta)=E(\sigma_*[\beta],\beta)$ and so $\widetilde{\sigma}=\sigma_*[\beta]$ by the uniqueness of the minimiser of $E(\cdot,\beta)$. 
Finally, note that all convergent subsequences of $\sigma_*[\beta^N]$ converge to $\sigma_*[\beta]$,  
hence the whole sequence must converge.
\end{proof}

For the proof of Theorem \ref{thm:solutions}, weak convergence of $\sigma_*[\alpha_t^N]$ is in fact sufficient. However,  uniform convergence of $\sigma_*[\alpha_t^N]$, along with the fact that $\sigma_*[\alpha^N_t]$ is Lipschitz continuous uniformly in $\alpha_t^N$ (see Lemma \ref{lem:fregularity}), 
implies the following. 

\begin{cor}[Regularity of the optimal source measure for general target measures]
Let $\beta \in \Pc(\Y)$ and $\sigma_*[\beta] \in \Pac(\X)$ be the optimal source measure. Then $\sigma_*[\beta]$ is Lipschitz continuous.
\end{cor}

\begin{proof}
This is an immediate consequence of 
Corollary \ref{cor: reg source} and Lemmas \ref{lem:fregularity}, \ref{5.1} and \ref{lem: uniform convergence}.
\end{proof}

Finally we prove our main theorem.

\begin{proof}[Proof of Theorem \ref{thm:solutions}]
Let $\alpha \in \C^{0,1}([0,\tau];\Pc(\Y))$ be the limit of the sequence $(\alpha^{N})_{N\in\N}$ obtained in Lemma \ref{lem:compact2}, and let $K \subset \Y$ be a compact set such that
$\mathrm{spt} ( \alpha_t^N ) \subset K$ and $\mathrm{spt} (\alpha_t ) \subset K$
for all $t\in[0,\tau]$ and $N\in\N$.
By \eqref{eq:WeakFormPDE} and \eqref{eq: Greenland}, for all $\varphi \in \D (\Y \times \R)$,
\begin{multline}
\label{eq: discrete weak}
\int_0^{\tau} \int_{\Y} 
\left[ \partial_t \varphi_t (\vb{y}) + (J\vb{y}) \cdot \nabla \varphi_t(\vb{y}) \right] \,\dd \alpha^N_t(\vb{y}) 
\, \dd t
- 
\int_0^{\tau} \int_{\X \times \Y}
(J\vb{x}) \cdot \nabla \varphi_t (\vb{y}) \, \dd \gamma[\alpha_t^N](\vb{x},\vb{y}) \, \dd t
\\
= \int_{\Y} \varphi_{\tau} (\vb{y}) \, \dd \alpha^N_{\tau}(\vb{y}) -\int_{\Y} \varphi_0 (\vb{y}) \,\dd \overline{\alpha}^N (\vb{y}).
\end{multline}
Recall the following standard estimate in optimal transport theory (see \cite[Exercise 38]{Santambrogio:2015} or \cite[Theorem 1.14]{Villani:2003}): 
\[
\left| \int_{K} \psi \, \dd (\alpha^N_t-\alpha_t) \right| 
\le \| \nabla \psi \|_{L^{\infty}(K)} W_1(\alpha^N_t,\alpha_t) \quad \forall \; \psi \in W^{1,\infty}(K).
\]
Then
\begin{align}
\nonumber    
& \left|
\int_0^{\tau} \int_{\Y} 
\left[ \partial_t \varphi_t (\vb{y}) + (J\vb{y}) \cdot \nabla \varphi_t(\vb{y}) \right] \, \dd (\alpha^N_t-\alpha_t) \!(\vb{y}) 
\, \dd t
\right|
\\
\nonumber
& \le 
\int_0^{\tau}
\max_{\vb{y} \in K} 
\big\{
\| 
\nabla_{\vb{y}} [\partial_t \varphi_t (\vb{y}) + (J\vb{y}) \cdot \nabla \varphi_t(\vb{y}) ]
\|
\big\}
\, W_1(\alpha^N_t,\alpha_t)
\, \dd t
\\
\label{eq: ingredient 1}
& \le
\tau
\max_{\substack{t \in [0,\tau] \\ \vb{y} \in K}} 
\big\{
\| 
\nabla_{\vb{y}} [\partial_t \varphi_t (\vb{y}) + (J\vb{y}) \cdot \nabla \varphi_t(\vb{y}) ]
\| \,
W_1(\alpha^N_t,\alpha_t)
\big\}
\to 0 \quad \textrm{as} \quad N \to \infty
\end{align}
by uniform convergence of $\alpha^N$ (seee \eqref{eq: uniform convergence}).
Pointwise-in-time convergence of $\alpha^N$ implies that
\[
\lim_{N \to \infty} \int_\Y \varphi_{\tau}(\vb{y}) \,\dd (\alpha_{\tau}^N-\alpha_{\tau}) \! (\vb{y})
= 0
\quad \text{and}\quad
\lim_{N \to \infty} \int_\Y \varphi_0(\vb{y}) \, \dd (\overline{\alpha}^N-\overline{\alpha}) \! (\vb{y})
= 0.
\]

By Lemma \ref{lem: uniform convergence}, for all $t \in [0,\tau]$, $\sigma_*[\alpha^N_t] \rightharpoonup \sigma_*[\alpha_t]$ in $\Pn(\X)$. Also, $\alpha_t^N \rightharpoonup \alpha_t$ in $\Pn(K)$.
Therefore $\gamma[\alpha^N_t] \rightharpoonup \gamma[\alpha_t]$ in $\Pn(\X \times K)$, where $\gamma[\alpha_t]$ is the unique optimal transport plan for transporting $\sigma_*[\alpha_t]$ to $\alpha_t$; see \cite[Theorem 5.20]{Villani:2009} . 
Define $F,F_N \in L^1([0,\tau])$ by
\[
F_N(t) =  \int_{\X \times \Y}
(J\vb{x}) \cdot \nabla \varphi_t (\vb{y}) \, \dd \gamma[\alpha_t^N](\vb{x},\vb{y}),
\quad 
F(t) =  \int_{\X \times \Y}
(J\vb{x}) \cdot \nabla \varphi_t (\vb{y}) \, \dd \gamma[\alpha_t](\vb{x},\vb{y}).
\]
Then $\lim_{N \to \infty} F_N(t) = F(t)$ for all $t \in [0,\tau]$ and
\[
\| F_N \|_{L^\infty([0,\tau])}
\le \f \max_{\vb{x} \in \X} \| \vb{x} \| \max_{t \in [0,\tau]} \| \nabla \varphi_t \|_{L^\infty(K)}. 
\]
By the Lebesgue Dominated Convergence Theorem,
\begin{equation}
\label{eq: last step}
\lim_{N \to \infty}
\int_0^{\tau} \int_{\X \times \Y}
(J\vb{x}) \cdot \nabla \varphi_t (\vb{y}) \, \dd \gamma[\alpha_t^N](\vb{x},\vb{y}) \, \dd t
=
\int_0^{\tau} \int_{\X \times \Y}
(J\vb{x}) \cdot \nabla \varphi_t (\vb{y}) \, \dd \gamma[\alpha_t](\vb{x},\vb{y}) \, \dd t.
\end{equation}
Combining equations \eqref{eq: discrete weak}--\eqref{eq: last step} shows that $\alpha$ is a weak solution of \eqref{eq:PDE}, as required.   
\end{proof}

\begin{rmk}[Convergence of the transport maps]
It can be shown that $T_{\alpha^N_t} \to T_{\alpha_t}$ as $N \to \infty$ in $L^p$ for $p \in (1,\infty)$.
This can be combined with the uniform convergence of the source measures to pass to the limit in the nonlinear term 
\[
\int_0^{\tau} \int_{\X}
(J\vb{x}) \cdot \nabla \varphi_t (T_{\alpha^N_t}(\vb{x})) \,\dd \sigma_*[\alpha^N_t](\vb{x}) \, \dd t
\]
to give an alternative, though more involved,  proof of Theorem \ref{thm:solutions}.
\end{rmk}

\begin{rmk}[Conservation of energy]
Proposition \ref{prop:energy-conserving} states that discrete weak solutions of \eqref{eq:PDE} conserve the geostrophic energy. In fact, the same is true for the weak solutions with arbitrary initial data $\overline{\alpha} \in \Pc(\Y)$ that we construct in the proof of Theorem \ref{thm:solutions}. We give a quick sketch of the proof here. 
Let $\alpha \in \C^{0,1}([0,\tau];\Pc(\Y))$ be a weak solution of \eqref{eq:PDE} with initial measure $\overline{\alpha}$, and 
let $\alpha^N\in \C^{0,1}([0,\tau];\Pn^N(\Y))$, $N \in \mathbb{N}$, be discrete weak solutions that uniformly approximate $\alpha$ in the sense of 
\eqref{eq: uniform}. 
We proved that $\alpha^N$ conserve energy in Proposition \ref{prop:energy-conserving}. We will use this to prove that $\alpha$ also conserves energy. For each $t \in [0,\tau]$ define $H_t,H_t^N:\Pn(\X) \to \mathbb{R} \cup \{ + \infty \}$ by 
\[
H_t(\sigma) = E(\sigma,\alpha_t), \qquad
H^N_t(\sigma) = E(\sigma,\alpha^N_t). 
\]
Then $H_t^N$ $\Gamma$-converges to $H_t$ with respect to weak convergence of measures. To see this, note that $E$ is lower semi-continuous because the transport cost $\mathcal{T}_c$
is continuous  \cite[Theorem 5.20]{Villani:2009} and $F$ is lower semi-continuous \cite[Proposition 7.7]{Santambrogio:2015}. Take any sequence $\sigma^N$ in $\Pn(\X)$ such that $\sigma^N \rightharpoonup \sigma$. Recall that $\alpha_t^N \rightharpoonup \alpha_t$. Then 
\[
H_t(\sigma) = E(\sigma,\alpha_t) \le \liminf_{N \to \infty} E(\sigma^N,\alpha_t^N) = \liminf_{N \to \infty} H^N_t(\sigma^N), 
\]
which is the liminf inequality for $\Gamma$-convergence. Moreover, the constant recovery sequence $\sigma^N=\sigma$ satisfies 
\[
\lim_{N \to \infty} H^N_t(\sigma^N) = \lim_{N \to \infty} E(\sigma,\alpha^N_t) = E(\sigma,\alpha_t) = H_t(\sigma),
\]
which is the limsup inequality for $\Gamma$-convergence.
Therefore $H^N_t \stackrel{\Gamma}{\to} H_t$, as claimed.
By a standard property of $\Gamma$-convergence (convergence of the minimum values of the energy) and Proposition \ref{prop:energy-conserving}, 
\begin{align*}    
E(\sigma_*[\alpha_t],\alpha_t)
&= \min H_t
= \lim_{N \to \infty}
\min H^N_t
= \lim_{N \to \infty} E(\sigma_*[\alpha^N_t],\alpha^N_t)
 = \lim_{N \to \infty} E(\sigma_*[\alpha^N_0],\alpha^N_0)
\\
& = 
\lim_{N \to \infty}  \min H_0^N
=
\min H_0
=
E(\sigma_*[\alpha_0],\alpha_0),
\end{align*}
which proves that $\alpha$ is energy-conserving, as claimed.
Note, however, that we have not proved that the energy is conserved by 
\emph{all} weak solutions of \eqref{eq:PDE}, only by those  
that can be approximated by discrete weak solutions as in Theorem \ref{thm:solutions}.
\end{rmk}

\section{Explicit examples}
\label{sec:examples}
In this section we present two explicit solutions of \eqref{eq:PDE}: a steady state solution and a discrete solution with a single seed.

\subsection{Steady state example}
In this section we make the additional assumption that the fluid domain satisfies $\X \subset \Y$.

\begin{prop}[Time-independent solution]
\label{prop:steady}
    Given $\ell\in\R$, define $\sigma_{\ell}:\X\to\R$ by
    \begin{equation*}
        \sigma_{\ell}(\vb{x}) := (f^*)'(\ell - g\ln x_3).
    \end{equation*}
    Then there exists $\ell_*\in\R$ such that 
    \begin{equation*}
        \int_{\X}\sigma_{\ell_*}(\vb{x}) \,\dd \vb{x} = 1.
    \end{equation*}
  Define $\alpha:[0,\tau] \to\Pc(\Y)$ by 
  \[
  \alpha_t=\sigma_{\ell_*}\Lc^3\measurerestr\X \quad \forall \; t\in [0,\tau].
  \]
  Then $\alpha$ is a weak steady state  solution of the compressible SG equation \eqref{eq:PDE}.
\end{prop}

\begin{rmk}[Formal derivation]
\label{rmk:measurederivation}
    We give a formal derivation of the steady state $\alpha$ before giving a rigorous proof of Proposition~\ref{prop:steady}.
    By Definition \ref{def:weak}, $\alpha \in \C([0,\tau];\Pc(\Y))$ is a steady solution of \eqref{eq:PDE} if $\alpha_t$ is independent of $t$, $\alpha_t = \sigma_*[\alpha_t]$ (in particular, $\mathrm{spt}(\alpha_t) \subseteq \X$), and $T_{\alpha_t}=\mathrm{id}$.
    Then
    \[
        \alpha_t = \argmin_{\sigma \in \Pac(\X)} E(\sigma,\alpha_t)
    \]
    for all $t \in [0,\tau]$. 
    It follows from \cite[Proposition 7.20]{Santambrogio:2015} that $\alpha_t$ satisfies the Euler-Lagrange equations
    \[
    \begin{cases} 
    \displaystyle
    \frac{\delta E}{\delta \sigma} 
    (\alpha_t, \alpha_t) = \ell & \textrm{on the support of } \alpha_t, 
    \phantom{\Bigg|}
    \\
    \displaystyle
    \frac{\delta E}{\delta \sigma} 
    (\alpha_t, \alpha_t) \ge \ell & \textrm{otherwise},
    \end{cases}
    \]
    for some $\ell\in\R$, where $\delta E/\delta \sigma$ denotes the first variation of $E$ with respect to $\sigma$. Formally we expect
    \[
    \frac{\delta E}{\delta \sigma}(\alpha_t,\alpha_t) = \varphi + f'(\alpha_t), 
    \]
    where $\varphi:\X \to \mathbb{R}$ is an optimal $c$-concave Kantorovich potential for transporting $\alpha_t$ to itself with cost $c$; see \cite[Section 7.2]{Santambrogio:2015}.
    Then
    \[
    \begin{cases} 
    \varphi + f'(\alpha_t) = \ell & \textrm{on the support of } \alpha_t, 
    \\
    \varphi + f'(\alpha_t) \ge \ell & \textrm{otherwise}.
    \end{cases}
    \]
    If $T_{\alpha_t}=\mathrm{id}$ is the optimal transport map from $\alpha_t$ to itself , then $\nabla \varphi(\vb{x}) = \nabla_{\vb{x}}c(\vb{x},\vb{x})$ for almost all $\vb{x}$ in the support of $\alpha_t$ \cite[Proposition 1.15]{Santambrogio:2015}. Integrating this expression gives $\varphi(\vb{x}) = g \ln x_3$, up to a constant. Then solving the Euler-Lagrange equations yields
    \[
    \alpha_t(\vb{x}) = (f')^{-1}(\ell - \varphi(\vb{x})) =
    (f^*)'(\ell - g \ln x_3)
    \]
    for all $\vb{x}$ in the support of $\alpha_t$.
    This is the steady state given in Proposition \ref{prop:steady}. The constant $\ell$ is determined by the constraint that $\alpha_t$ is a probability measure. 
\end{rmk}

\begin{proof}[Proof of Proposition \ref{prop:steady}]
    First we show that there exists $\ell_*\in\R$ such that 
        $\int_{\X}\sigma_{\ell_*}(\vb{x}) \,\dd \vb{x} = 1.$
    Since $X \subset \Y = \R^2 \times (\delta,1/\delta)$ and $(f^*)'$ is non-decreasing, then for $\ell\in\R$
    \[
    (f^*)'(\ell - g \ln x_3) \ge (f^*)'(\ell - g \ln (1/\delta))
    \]
    for all $\vb{x} \in \X$. Therefore
    \[
    \int_{\X} \sigma_{\ell}(\vb{x}) \,\dd \vb{x}
    \ge (f^*)'(\ell - g \ln (1/\delta)) |\X| \to \infty \quad \textrm{as} \quad \ell \to \infty.
    \]
    Let $\ell_{\delta}=g\ln \delta$.
    Since $x_3 \ge \delta$ for all $\vb{x} \in \X$,  
    \[
    \int_{\X} \sigma_{\ell_{\delta}}(\vb{x}) \,\dd \vb{x}
    = \int_{\X} (f^*)'(\ell_{\delta} - g \ln x_3) \,\dd \vb{x} 
    \leq \int_{\X} (f^*)'(\ell_{\delta} - g \ln \delta) \,\dd \vb{x} 
    = (f^*)'(0) |\X| 
    =0.
    \]
    The map $\ell \mapsto \int_\X \sigma_{\ell}(\vb{x}) \, \dd \vb{x}$ is continuous. Therefore, by the Intermediate Value Theorem, there exists $\ell_*\in\R$ such that 
    \begin{equation*}
        \int_{\X}\sigma_{\ell_*}(\vb{x}) \,\dd \vb{x} = 1,
    \end{equation*}
    as claimed. For the remainer of this proof, we take $\ell=\ell_*$

    Next we show that the identity map is the optimal transport map for transporting any measure $\sigma \in \Pn(\X)$ to itself, and that an optimal Kantorovich potential pair is $(\varphi,\varphi^c$), where $\varphi:\X \to \R$ is given by $\varphi(\vb{x})=g \ln x_3$, and $\varphi^c:\Y \to \R$ is given by
\[
\varphi^c(\vb{y}) = \min_{\vb{x}\in\X} \, ( c(\vb{x},\vb{y})-\varphi(\vb{x}))
= \min_{\vb{x}\in\X} \, ( c(\vb{x},\vb{y})-g \ln x_3).
\]
First we compute an explicit expression for $\varphi^c(\vb{y})$ for the special case $\vb{y} \in \X \cap \Y = \X$. Given $\vb{y} \in \X$, define $\theta_{\vb{y}}:\X \to \R$ by $\theta_{\vb{y}}(\vb{x})=c(\vb{x},\vb{y})-g \ln x_3$. Then 
\[
\nabla \theta_{\vb{y}}(\vb{x})
= 
\begin{pmatrix}
\frac{\f^2}{y_3}(x_1-y_1) 
\\ 
\frac{\f^2}{y_3}(x_2-y_2) 
\phantom{\bigg|}
\\ 
\frac{g}{y_3}-\frac{g}{x_3}
\end{pmatrix},
\qquad
D^2 \theta_{\vb{y}}(\vb{x})
= 
\begin{pmatrix}
\frac{\f^2}{y_3} & 0 & 0
\\
0 & \frac{\f^2}{y_3} & 0
\phantom{\bigg|}
\\
0 & 0 & \frac{g}{x_3^2}
\end{pmatrix}.
\]
Therefore $D^2 \theta_{\vb{y}}(\vb{x})$ is positive definite for all $\vb{x} \in \X$. Hence $\theta_{\vb{y}}$ is strictly convex, its minimiser is $\vb{y}$, and its minimum value is 
\[
\varphi^c(\vb{y}) = \min_{\vb{x} \in \X} \theta_{\vb{y}}(\vb{x}) 
= \theta_{\vb{y}}(\vb{y}) = 
c(\vb{y},\vb{y}) - g \ln y_3 = 
g - g \ln y_3. 
\]
Note that $\varphi(\vb{x})+\varphi^c(\vb{x})=g$ for all $\vb{x} \in \X$ and that $\sigma$ is supported on $\X \subset \Y$.
Then 
\[
\mathcal{T}_c(\sigma,\sigma) 
\le 
\int_{\X} c(\vb{x},\vb{x}) \,
\dd \sigma(\vb{x}) 
= g
=  \int_\X \varphi \, \dd \sigma  + \int_\Y \varphi^c \, \dd \sigma 
\le 
\mathcal{T}_c(\sigma,\sigma),
\]
which proves that the identity map is optimal for transporting $\sigma \in \Pn(\X)$ to itself, and that $(\varphi,\varphi^c$) is an optimal Kantorovich potential pair, as claimed.
    
    Lastly, we show that 
    $\sigma_{\ell} = \argmin_{\sigma \in \Pac(\X)} E(\sigma,\sigma_{\ell})$.
    If $\vb{x} \in \X$ satisfies $\sigma_{\ell}(\vb{x}) = 0$, then $0 \ge \ell - g \ln x_3 = \ell - \varphi(\vb{x})$. Therefore 
    \[
    \varphi(\vb{x}) + f'(\sigma_{\ell}(\vb{x})) = \varphi(\vb{x}) \ge \ell.
    \]
    On the other hand, if $\vb{x} \in \X$ satisfies $\sigma_{\ell}(\vb{x}) > 0$, then
    \[
        \varphi(\vb{x}) + f'(\sigma_{\ell}(\vb{x}))
    = g \ln x_3 + f'((f^*)'(\ell - g \ln x_3)) = g \ln x_3 + \ell - g \ln x_3 = \ell
    \]
    because $f'((f^*)'(t))=t$ for all $t>0$. In summary,
    \begin{equation}
    \label{eq: Jeff Bezos}    
\varphi(\vb{x}) + f'(\sigma_{\ell}(\vb{x})) \ge \ell \quad \forall \; \vb{x} \in \X,
    \end{equation}
and
\begin{equation}
\label{eq: Shabaka Hutchings}    
\big( \varphi(\vb{x}) + f'(\sigma_{\ell}(\vb{x})) \big)  \sigma_{\ell}(\vb{x}) = \ell \sigma_{\ell}(\vb{x})
\quad \forall \; \vb{x} \in \X.
\end{equation}
By the Kantorovich Duality Theorem and the convexity of $f$, for all $\sigma \in\Pac(\X)$,
\begin{align}
\nonumber
        E(\sigma,\sigma_{\ell}) 
        & = 
        \T(\sigma,\sigma_{\ell})+\int_{\X} f(\sigma) \,\dd \mathcal{L}^3 
        \\
        \nonumber
        & \ge 
        \int_{\X} \varphi \,\dd \sigma + \int_{\Y} \varphi^c \, \dd \sigma_{\ell} + \int_{\X} 
        \big( f(\sigma_{\ell}) + f'(\sigma_{\ell})(\sigma - \sigma_{\ell}) \big)
         \,\dd \mathcal{L}^3
         \\
         \label{eq: 6.3}
         & = 
         E(\sigma_{\ell}, \sigma_{\ell}) + \int_{\X} \big( \varphi + f'(\sigma_{\ell}) \big) ( \sigma -\sigma_{\ell} ) \,\dd \mathcal{L}^3
\end{align}
since $\varphi$ is an optimal Kantorovich potential for the transport from $\sigma_{\ell}$ to itself.
Therefore, by \eqref{eq: Jeff Bezos}, \eqref{eq: Shabaka Hutchings} and \eqref{eq: 6.3}, 
\[
E(\sigma,\sigma_{\ell})
\ge
        E(\sigma_{\ell}, \sigma_{\ell})
        + \ell \int_\X ( \sigma -\sigma_{\ell} ) \,\dd \mathcal{L}^3
        =
        E(\sigma_{\ell}, \sigma_{\ell}),
    \]
and so
$\sigma_{\ell} = \argmin_{\sigma \in \Pac(\X)} E(\sigma,\sigma_{\ell})$, as claimed.

    In conclusion, we have shown that $\sigma_*[\sigma_{\ell}] = \sigma_{\ell}$ and $T_{\sigma_{\ell}} = \mathrm{id}$. Therefore the time-independent map $[0,\tau] \ni t \mapsto \alpha_t = \sigma_{\ell}$ is a weak solution of \eqref{eq:PDE}, as required.
\end{proof}

\subsection{Single seed example}
When there is only one seed, it is possible to construct an explicit solution for the case $\gamma=2$. This is because the single Laguerre cell is the entire domain. Note that, to be physically meaningful, $\gamma\in(1,2)$. However, for these values of $\gamma$, we are not able to derive such an explicit solution.
\begin{prop}[Elliptic orbit of a single seed]
Fix $\gamma=2$, $\kappa=\frac{1}{2}$ 
and $\tau >0$. 
Let $\X=[-a,a]\times[-b,b]\times[0,h]$, where $a,b,h>0$. 
Let $\overline{\vb{z}} \in \Y$ be the initial position of the particle. 
Define $\vb{z} : (-\infty,\infty) \to \Y$ to be the solution of the linear ODE
\begin{align*}
\dot{\vb{z}}(t) & =
\f 
\begin{pmatrix}
0 & -A & 0 
\\
B & 0 & 0
\\
0 & 0 & 0
\end{pmatrix} \vb{z}(t),
\\
\vb{z}(0) & = \overline{\vb{z}},
\end{align*} 
where 
\[
A = 1 - \frac{|\X| \f^2}{3 \overline{z}_3} b^2, \qquad B = 1 - \frac{|\X| \f^2}{3 \overline{z}_3} a^2, \qquad |\X|=4abh.
\]
Assume that $\overline{z}_3 > \frac{|\X| \f^2 }{3} \max \{a^2,b^2\}$ so that $A>0$ and $B>0$. Define the ellipse $\mathcal{E} \subset \Y$ by
\[
\mathcal{E} \coloneqq 
\left\{ \vb{z}(t) : t \in (-\infty,\infty) \right\} = 
\left\{ (y_1,y_2,\overline{z}_3) \in \R^2 \times \{ \overline{z}_3 \} \, : \, \frac{y_1^2}{A} + \frac{y_2^2}{B} = \frac{\overline{z}_1^2}{A} + \frac{\overline{z}_2^2}{B} \right\}.
\]
Assume additionally that 
$\overline{z}_3$  is sufficiently large so that
\begin{equation}
\label{eq: Tinderbox Orchestra}
c(\vb{x},\vb{y}) - \frac{1}{|\X|} \int_\X c(\cdot,\vb{y}) \, \dd \mathcal{L}^3 < \frac{1}{|\X|} \quad \forall 
\; (\vb{x},\vb{y}) \in \X \times \mathcal{E}.
\end{equation}
Define $\alpha : [0,\tau] \to \Pc(\Y)$ by $\alpha_t = \delta_{\vb{z}(t)}$. Then the following hold:
\begin{enumerate}
    \item The optimal weight map $w_* : D \to \R$ (defined in \eqref{def:opt_w}) satisfies the following:
    \[
    w_*(\vb{z}(t)) = \frac{1}{|\X|} \left( 1 + \int_\X c(\vb{x},\vb{z}(t)) \, \dd \vb{x} \right)
    \quad \forall \; t \in (-\infty,\infty).
    \]
    \item The optimal source measure $\sigma_*[\alpha_t] \in \Pac(\X)$ is given by
    \[
    \sigma_*[\alpha_t](\vb{x}) = 
    w_*(\vb{z}(t)) - c(\vb{x},\vb{z}(t)).
    \]
    \item $\alpha$ is a weak solution of the compressible SG equation \eqref{eq:PDE}. In particular, the trajectory of the seed is the ellipse $\mathcal{E}$. Moreover, if $a=b$, then the seed moves on a circle. If $\overline{z}_1=\overline{z}_2 = 0$, then the seed is stationary. 
\end{enumerate}
\end{prop}

\begin{proof}
Since $\gamma =2$ and $\kappa = \frac 12$,
\[
f^*(t)=\begin{cases}
\frac{1}{2}t^{2} & \textrm{if } t>0,
\\
0 & \textrm{if } t \le 0.
\end{cases}
\]
Define $w:(-\infty,\infty) \to \mathbb{R}$ by
\[
w(t) = \frac{1}{|\X|} \left( 1 + \int_\X c(\vb{x},\vb{z}(t)) \, \dd \vb{x} \right).
\]
Then, for all $\vb{x} \in \X$, $t \in (-\infty,\infty)$, 
\[
w(t) - c(\vb{x},\vb{z}(t)) 
= \frac{1}{|\X|} \left( 1 + \int_\X c(\tilde{\vb{x}},\vb{z}(t)) \, \dd \tilde{\vb{x}} \right)
- c(\vb{x},\vb{z}(t)) > 0
\]
by \eqref{eq: Tinderbox Orchestra}. Therefore
\[
\frac{\partial \G}{\partial w}(w(t),\vb{z}(t))
= 1 - \int_\X (f^*)'( w(t) - c(\vb{x},\vb{z}(t))  ) \, \dd \vb{x}  
= 1 - \int_\X ( w(t) - c(\vb{x},\vb{z}(t))  ) \, \dd \vb{x}
= 0,
\]
hence $w(t)$ maximises the concave function $\G(\cdot,\vb{z}(t))$.
Since the maximiser of $\G(\cdot,\vb{z}(t))$ is unique by Theorem \ref{thm:dualitytheorem}, it follows that
$w_*(\vb{z}(t))=w(t)$, as claimed. The expression for $\sigma_*[\alpha_t]$ then follows immediately from \eqref{eq:density}.

It remains to prove that $\alpha$ is a weak solution the compressible SG equation \eqref{eq:PDE}.
By Proposition \ref{lem:Correspondence}, it suffices to show that $\vb{z}$ satisfies the following ODE:
\begin{align*}
\dot{z}_1(t) & = \f (- z_2(t) + C_2(\vb{z}(t))),    
\\
\dot{z}_2(t) & = \f(z_1(t)-C_1(\vb{z}(t))),
\\
\dot{z}_3(t) & = 0.
\end{align*}
In other words, we need to show that 
\[
C_1(\vb{z}(t)) = (1-B) z_1(t), \qquad 
C_2(\vb{z}(t)) = (1-A) z_2(t).
\]
We will just verify the expression for $C_1$; the expression for $C_2$ is similar. By Definition \ref{def:centroid} and the symmetry of the domain $\X$,
\begin{align*}
C_1(\vb{z}(t)) & = \int_{\X} x_1 \, \sigma_*[\alpha_t](\vb{x}) \, \dd \vb{x}
\\
& = \int_{\X} x_1  \left( w_*(\vb{z}(t)) - c(\vb{x},\vb{z}(t)) \right) \, \dd \vb{x}
\\
& = 
- \frac{1}{\overline{z}_3} \int_{-a}^a \int_{-b}^b \int_{0}^h x_1  \left(\frac{\f^2}{2}(x_1-z_1(t))^2+\frac{\f^2}{2}(x_2-z_2(t))^2+g x_3\right) \, \dd x_3\, \dd x_2\, \dd x_1
\\
& = \frac{\f^2}{\overline{z}_3} \int_{-a}^a \int_{-b}^b \int_{0}^h x_1^2 z_1(t) \, \dd x_3\, \dd x_2\, \dd x_1
\\
& = \frac{4 a^3 b h \f^2}{3 \overline{z}_3} z_1(t)
\\
& = (1-B) z_1(t),
\end{align*}
as required.    
\end{proof}

\section*{Acknowledgements} 
We thank Mike Cullen for stimulating discussions related to this work.
DPB would like to thank the UK Engineering and Physical Sciences Research Council (EPSRC) for financial support via the grant EP/V00204X/1. 
CE gratefully acknowledges the support of the EPSRC via the grant EP/W522570/1, and of the German Research Foundation (DFG) through the CRC 1456 `Mathematics of Experiment', project A03.
TPL was supported by The Maxwell Institute Graduate School in Modelling, Analysis, and Computation, a Centre for Doctoral Training funded by the EPSRC (grant EP/S023291/1), the Scottish Funding Council, Heriot-Watt University and the University of Edinburgh. 
BP gratefully acknowledges the support of the EPSRC via the grant EP/P011543/1.

\appendix



\section{Assumptions on the fluid domain $\X$}
\label{sec:domain}
Recall that $\Phi:\R^3 \to \R^3$ is the diffeomorphism defined by
\[
    \Phi(\vb{x}) = \left( \f^{-2} \, x_1,  \,\f^{-2} \, x_2, \, g^{-1} \left( x_3 - \tfrac 12 \f^{-2} \left( x_1^2 + x_2^2 \right) \right) \right).
\]
In this appendix we discuss the assumption that the fluid domain $\X$ satisfies
\begin{equation}
\label{eq: JD Vance}
\Phi^{-1}(\X) \text{ is convex.}
\end{equation}
In particular, 
we give examples of domains $\X$ satisfying \eqref{eq: JD Vance},
we prove that \eqref{eq: JD Vance} is equivalent to $\X$ being $c$-convex \cite[Definition 1.2]{Kitagawa:2019},  and we prove that, given any $c$-Laguerre tessellation $\{ \L \}_{i=1}^N$, then 
$\{ \Phi^{-1}(\L) \}_{i=1}^N$ is a classical Laguerre tessellation (with respect to the quadratic cost). In particular, $\Phi^{-1}(\L)$ is a convex polyhedron for all $i$.   
Moreover, we will see that the cost $c$ satisfies Loeper's condition \cite[Definition 1.1]{Kitagawa:2019}.

\begin{rmk}[Interpreting $\Phi$ in terms of $c$-exponential maps]
Let $\vb{y}_0 = (0,0,1)$. It is easy to check that $\Phi$ is the $c$-exponential map $\exp_{\vb{y}_0}^c$ \cite[Remark 4.4]{Kitagawa:2019}:
\[
\Phi =  \exp_{\vb{y}_0}^c := (-D_{\vb{y}} c(\cdot,\vb{y}_0))^{-1}.  
\]
In other words,
\[
- D_{\vb{y}} c( \Phi(\vb{x}),\vb{y}_0) = \vb{x}.
\]
Note that $c$-exponential maps play an important role in regularity theory for semi-discrete optimal transport \cite{Kitagawa:2019}.
\end{rmk}

\begin{example}[Domains satisfying assumption \eqref{eq: JD Vance}]
\label{example:c-convex set}
Let $(a_1,a_2,a_3) \in \mathbb{R}^3$.
It can be shown that any paraboloid of the form
\begin{equation}\label{eq: c-convex set}
    \X_1 = 
    \left\{ 
    \vb{x} \in \mathbb{R}^3 : 
    x_3 \le - \frac{\f^2}{2g} \left( (x_1 - a_1)^2 + (x_2 - a_2)^2 \right) + a_3  
    \right\}
\end{equation}
satisfies \eqref{eq: JD Vance}.
Its inverse image $\Phi^{-1}(\X_1)$ is a half-space 
of the form 
\begin{equation*}
    \Phi^{-1}(\X_1) = \{ \vb{p} \in \mathbb{R}^3 : p_3 \le b_1 p_1 + b_2 p_2 + c\},
\end{equation*}
where $b_1,b_2,c \in \mathbb{R}$. In particular, it is convex.
Replacing the inequality $\le$
in \eqref{eq: c-convex set} by $\ge$  gives another example of a set satisfying \eqref{eq: JD Vance}.
Another example is the half-space
\begin{equation*}
    \X_2 = \{ \vb{x} \in \mathbb{R}^3 :
    \vb{x} \cdot \vb{n} \ge d
    \}
\end{equation*}
where $d \in \mathbb{R}$ and $\vb{n} \in \mathbb{R}^3$ with $n_3 \ge 0$.
Its inverse image $\Phi^{-1}(\X_2)$ has the following form
for some $\vb{b} \in \mathbb{R}^3$, $c \in \mathbb{R}$: 
\begin{equation}\label{eq: c-convex set 2}    
    \Phi^{-1}(\X_2) =
    \left\{
    \vb{p} \in \mathbb{R}^3 : 
    \frac{n_3}{g \f^2} (p_1^2 + p_2^2) + \vb{b} \cdot \vb{p} + c \le 0
    \right\}.
\end{equation}
This is convex because it has the form $\{ \vb{p} \in \mathbb{R}^3 : F(\vb{p}) \le 0\}$ with $F$ convex.
We can build other examples by using the fact that any intersection of sets satisfying \eqref{eq: JD Vance} also satisfies \eqref{eq: JD Vance}. For example, ``blister pack" sets of the form 
\begin{equation*}
    \X_3 = \left\{ \vb{x} \in \mathbb{R}^3 : 0 \le 
    x_3 \le - \frac{\f^2}{2g} \left( (x_1 - a_1)^2 + (x_2 - a_2)^2 \right) + a_3 
    \right\}
\end{equation*}
satisfy \eqref{eq: JD Vance} for all $\vb{a} \in \mathbb{R}^3$.
It can be seen from \eqref{eq: c-convex set 2}
that the half-space
$\{ \vb{x}\in\R^3:\vb{x}\cdot\vb{n} \ge d\}$ does \emph{not} satisfy \eqref{eq: JD Vance} 
if $n_3 < 0$.
\end{example}

\begin{rmk}[Rectangular domains do not satisfy \eqref{eq: JD Vance}]
Rectangular domains of the form  $\X = [a_1,b_1] \times [a_2,b_2] \times [a_3,b_3]$ do not satisfy assumption \eqref{eq: JD Vance}.
These domains are perhaps the most natural from a physical point of view. However, in this case an alternative, bespoke argument can be used to prove 
Corollary \ref{cor: regularity},
and hence Theorems \ref{thm:discretesolutions}, \ref{thm:solutions}; see \cite{Lavier:2025}. 
\end{rmk}

Recall the following from \cite[Definitions 1.1 and 1.2]{Kitagawa:2019}: The cost $c$ satisfies \emph{Loeper's condition} if, for each $\vb{y} \in \Y$, there exists a diffeomorphism $\Phi_{\vb{y}} : \R^3 \to \R^3$ such that the map
\[
\R^3 \ni \vb{x} \mapsto c(\Phi_{\vb{y}}(\vb{x}),\vb{y})
    - c(\Phi_{\vb{y}}(\vb{x}),\vb{z})
\]
is quasi-convex for all $\vb{z} \in \Y$ (meaning that its sublevel sets are convex). If in addition the set $\Phi_{\vb{y}}^{-1}(\X)$ is convex for each $\vb{y} \in \Y$, then we say that $\X$ is \emph{$c$-convex}. 

\begin{lemma}[Loeper's condition and $c$-convexity]
\label{lemma: Loeper's condition}
Let $\vb{x} \in \X$, $\vb{y}, \vb{z} \in \Y$. Then
\begin{equation}
\label{eq: Department of Government Efficiency}    c(\Phi(\vb{x}),\vb{y}) = - \frac{1}{y_3}
\begin{pmatrix}
x_1 \\ x_2 \\ x_3    
\end{pmatrix}
\cdot 
\begin{pmatrix}
y_1 \\ y_2 \\ -1    
\end{pmatrix}
+ \frac{\f^2}{2 y_3} \left( y_1^2 + y_2^2 \right).
\end{equation}
In particular, 
$u: \mathbb{R}^3 \to \mathbb{R}$ defined by
\begin{equation*}
    u(\vb{x}) = c(\Phi(\vb{x}),\vb{y})
    - c(\Phi(\vb{x}),\vb{z})
\end{equation*}
is affine and the cost $c$ satisfies Loeper's condition. 
Moreover, $\X$ satisfies assumption \eqref{eq: JD Vance} if and only if $\X$ is $c$-convex. 
\end{lemma}

\begin{proof}
Equation \eqref{eq: Department of Government Efficiency} is a direct computation. 
Since $u$ is affine, then it is quasi-convex. 
Therefore, by definition,
$c$ satisfies Loeper's condition (with $\Phi_{\vb{y}}=\Phi$ independent of $\vb{y}$), and assumption \eqref{eq: JD Vance} is equivalent to $c$-convexity of $\X$. 
\end{proof}


To prove that $c$ satisfies Loeper's condition, we could have also chosen the $c$-exponential maps $\Phi_{\vb{y}} =  \exp_{\vb{y}}^c := (-D_{\vb{y}} c(\cdot,\vb{y}))^{-1}$ depending on $\vb{y}$.
However, the fact that $\Phi$ is independent of $\vb{y}$ (and that $u$ is affine) has important consequences, as we will now see. 
Let
\[
D^N_2:=\left\{ \vb{y}= \left(\vb{y}^1,\ldots,\vb{y}^N \right) \in (\R^3)^{N}: \vb{y}^i\neq \vb{y}^j \textrm{ whenever } i \neq j \right\}.
\]
Recall the following definition:

\begin{defn}[Classical Laguerre tessellations]
    Given $(\vb*{\psi},\vb{y}) \in \R^N \times D_2^N$,  the (classical) \emph{Laguerre tessellation} of $\Phi^{-1}(\X)$ generated by $(\vb*{\psi},\vb{y})$ is the partition $\{ L^i_2 (\vb*{\psi},\vb{y}) \}_{i=1}^N$ defined by
\[
L^i_2 (\vb*{\psi},\vb{y})
= \left\{
\vb{x} \in \Phi^{-1}(\X) : 
\| \vb{x} - \vb{y}^i \|^2 - \psi^i \le \| \vb{x} - \vb{y}^j \|^2 - \psi^j \; \forall \, j \in \{ 1,\ldots,N \} 
\right\}.
\]
\end{defn}

The Laguerre cells $L^i_2 (\vb*{\psi},\vb{y})$ are the intersection of $\Phi^{-1}(\X)$ with convex polyhedra. Consequently, unlike $c$-Laguerre cells, they can be computed very efficiently. Classical Laguerre tessellations arise in semi-discrete optimal transport problems with \emph{quadratic} cost. The following lemma asserts that $c$-Laguerre tessellations are just classical Laguerre tessellations in disguise. It also gives a practical way of computing them.

\begin{lemma}[Rewriting $c$-Laguerre tessellations as classical Laguerre tessellations]
\label{lemma: Laguerre cells are polygons}
Define $\widehat{\vb{y}}:\Y^N \to (\R^3)^N$ by $\widehat{\vb{y}}(\vb{z})=(\widehat{\vb{y}}^1(\vb{z}),\ldots,\widehat{\vb{y}}^N(\vb{z}))$, where
\[
\widehat{\vb{y}}^i(\vb{z}) 
= \frac{1}{2 z_3^i} 
\begin{pmatrix}
z_1^i \\ z_2^i \\ -1    
\end{pmatrix}.
\]
Define $\widehat{\vb*{\psi}}: \R^N \times \Y^N \to \R^N$ by
\[
\widehat{\psi}^i(\vb*{w},\vb{z}) = 
 w^i + \left( \frac{z_1^i}{2z^i_3} \right)^2 + \left( \frac{z^i_2}{2z^i_3} \right)^2 + \left( \frac{1}{2z^i_3} \right)^2 - \frac{\f^2}{2z^i_3}\left(\left(z^i_1\right)^2 + \left(z^i_2\right)^2\right).
\]
Then, for all $i \in \{1,\ldots,N\}$ and $(\vb*{w},\vb{z}) \in \R^N \times D^N$,
\[
\Phi^{-1}(L^i_c (\vb*{w},\vb{z})) = L^i_2(\widehat{\vb*{\psi}}(\vb*{w},\vb{z}),\widehat{\vb{y}}(\vb{z})).
\]
\end{lemma}

\begin{proof}
By \eqref{eq: Department of Government Efficiency}, for all $\vb{x} \in \R^3$, $\vb{z} \in \Y$, $\vb*{w} \in \R^N$,
\begin{align}
\nonumber
c(\Phi(\vb{x}),\vb{z}^i) - w^i & =     
- \frac{1}{z_3^i}
\begin{pmatrix}
x_1 \\ x_2 \\ x_3    
\end{pmatrix}
\cdot 
\begin{pmatrix}
z_1^i \\ z_2^i \\ -1 
\end{pmatrix}
+ \frac{\f^2}{2 z_3^i} \left( \left(z_1^i\right)^2 + \left(z_2^i\right)^2 \right) - w^i
\\
\label{eq: Ezra Collective - Charlie's favourite band}
& = \| \vb{x} - \widehat{\vb{y}}^i(\vb{z}) \|^2 - \widehat{\psi}^i(\vb*{w},\vb{z}) - \| \vb{x} \|^2.
\end{align}
Let $\vb{x} \in \Phi^{-1}(\X)$. Then
\begin{align*}
\vb{x} \in
\Phi^{-1}(L^i_c (\vb*{w},\vb{z})) 
\quad 
& \Longleftrightarrow \quad 
c(\Phi(\vb{x}),\vb{z}^i) - w^i
\le c(\Phi(\vb{x}),\vb{z}^j) - w^j \quad \forall \; j
\\
& \Longleftrightarrow \quad 
\| \vb{x} - \widehat{\vb{y}}^i(\vb{z}) \|^2 - \widehat{\psi}^i(\vb*{w},\vb{z}) \le 
\| \vb{x} - \widehat{\vb{y}}^j(\vb{z}) \|^2 - \widehat{\psi}^j(\vb*{w},\vb{z}) \quad \forall \; j
\\
& \Longleftrightarrow \quad 
\vb{x} \in 
L^i_2(\widehat{\vb*{\psi}}(\vb*{w},\vb{z}),\widehat{\vb{y}}(\vb{z})),
\end{align*}
as required.
\end{proof}

\section{Proof of Lemma \ref{lem:dualdifferentiability}}
\label{Sec: Proof of Lemma 3.8}

In this section we prove 
Lemma \ref{lem:dualdifferentiability}. The main ingredient of the proof is the following lemma.

\begin{lemma}[Regularity theory for semi-discrete optimal transport]
\label{cor: regularity}
Let $\zeta: \X \to \mathbb{R}$ be continuous. 
Let $U \subset \R^N \times D_0^N$ be the open set
\[
U = \bigg\{ (\vb*{w},\vb{z}) \in \R^N \times D_0^N \, : \, \mathcal{L}^3 (L^i_c(\vb*{w},\vb{z})) > 0 \; \forall \; i \in \{1,\ldots,N\} \bigg\}.
\]
Define 
$\vb{\Psi}=(\Psi^1,\ldots,\Psi^N):U\to\R^N$
by
\begin{equation*}
    \Psi^i(\vb*{w},\vb{z}):=\int_{\L} \zeta(\vb{x}) \,\dd \vb{x}.
\end{equation*}
Then $\vb{\Psi}$ is continuously differentiable.
Moreover, for all $j \in \{ 1, \ldots , N\} \ \setminus \{ i \}$,
\begin{align}
\label{eq: 1}
\frac{\partial \Psi^i}{\partial w^j}(\vb*{w}, \vb{z}) 
& = 
- \int_{L^i_c(\vb*{w},\vb{z}) \cap L^j_c(\vb*{w},\vb{z})} 
\frac{\zeta(\vb{x})}{\| \nabla_{\vb{x}} c(\vb{x},\vb{z}^i) -  \nabla_{\vb{x}} c(\vb{x},\vb{z}^j) \|} \, \dd  \mathcal{H}^2(\vb{x}),
\\
\label{eq: 
2}
\frac{\partial \Psi^i}{\partial \vb{z}^j}(\vb*{w}, \vb{z}) 
& = 
\int_{L^i_c(\vb*{w},\vb{z}) \cap L^j_c(\vb*{w},\vb{z})} 
\frac{\nabla_{\vb{y}}c(\vb{x},\vb{z}^j) \, \zeta(\vb{x})}{\| \nabla_{\vb{x}} c(\vb{x},\vb{z}^i) -  \nabla_{\vb{x}} c(\vb{x},\vb{z}^j) \|} \, \dd  \mathcal{H}^2(\vb{x}),
\end{align}
and
\begin{equation}
\label{eq: 3}
\frac{\partial \Psi^i}{\partial w^i}(\vb*{w}, \vb{z}) 
= - \sum_{j \ne i} \frac{\partial \Psi^j}{\partial w^i}(\vb*{w}, \vb{z}), 
\qquad    
\frac{\partial \Psi^i}{\partial \vb{z}^i}(\vb*{w}, \vb{z}) 
= - \sum_{j \ne i} \frac{\partial \Psi^j}{\partial \vb{z}^i}(\vb*{w}, \vb{z}). 
\end{equation}
\end{lemma}

This lemma could be proved using \cite[Theorem 1]{deGournay:2019}, by proving that our transport cost $c$ satisfies assumptions (Diff-2) and (Cont-2) of \cite{deGournay:2019}. However, this is very involved. We give a shorter proof using Lemma \ref{lemma: Laguerre cells are polygons} to reduce to the case of the classical quadratic transport cost. 

\begin{proof}
Note that $\mathrm{det}(D \Phi)=\f^{-4}g^{-1}>0$ is constant.
By the change of variables 
formula and Lemma \ref{lemma: Laguerre cells are polygons}, 
\[
\Psi^i(\vb*{w},\vb{z})=
\int_{\L} \zeta \, \dd \mathcal{L}^3 
= \int_{\Phi^{-1}(\L)} \zeta \circ \Phi  \, \mathrm{det} (D \Phi ) \, \dd \mathcal{L}^3 
= \mathrm{det} (D \Phi ) \int_{L^i_2(\widehat{\vb*{\psi}}(\vb*{w},\vb{z}),\widehat{\vb{y}}(\vb{z}))}
\zeta \circ \Phi \,\dd \mathcal{L}^3.
\]
Let $U_2 \subset \R^N \times D_2^N$ be the open set
\[
U_2 = \bigg\{ (\vb*{\psi},\vb{y}) \in \R^N \times D_2^N \, : \, \mathcal{L}^3 (L^i_2(\vb*{\psi},\vb{y})) > 0 \; \forall \; i \in \{1,\ldots,N\} \bigg\}.
\]
Define 
$\vb{\Psi}_2=(\Psi^1_2,\ldots,\Psi^N_2):U_2\to\R^N$
by
\begin{equation*}
    \Psi^i_2(\vb*{\psi},\vb{y}):=\int_{L^i_2(\vb*{\psi},\vb{y})} \zeta(\Phi(\vb{x})) \,\dd \vb{x}.
\end{equation*}
Then
\begin{equation}
\label{eq: Jack Quaid}    
\vb{\Psi}(\vb*{w},\vb{z})
= \mathrm{det} (D \Phi ) 
\vb{\Psi}_2(\widehat{\vb*{\psi}}(\vb*{w},\vb{z}),\widehat{\vb{y}}(\vb{z})).
\end{equation}
By assumption \eqref{eq: JD Vance} and \cite[Proposition 2]{deGournay:2019}, $\vb{\Psi}_2$ is continuously differentiable and, for all $j \in \{ 1, \ldots , N\} \ \setminus \{ i \}$,
\begin{align}
\label{eq: John Wick 1}
\frac{\partial \Psi_2^i}{\partial \psi^j}(\vb*{\psi}, \vb{y}) 
& = 
- \frac{1}{2 \| \vb{y}^i - \vb{y}^j \|}\int_{L^i_2(\vb*{\psi},\vb{y}) \cap L^j_2(\vb*{\psi},\vb{y})} 
\zeta(\Phi(\vb{x})) \, \dd  \mathcal{H}^2 (\vb{x}),
\\
\label{eq: John Wick 2}
\frac{\partial \Psi_2^i}{\partial \vb{y}^j}(\vb*{\psi}, \vb{y}) 
& = 
\frac{1}{\| \vb{y}^i - \vb{y}^j \|} \int_{L^i_2(\vb*{\psi},\vb{y}) \cap L^j_2(\vb*{\psi},\vb{y})} 
(\vb{y}^j - \vb{x}) \zeta(\Phi(\vb{x})) \, \dd  \mathcal{H}^2(\vb{x}),
\end{align}
and
\begin{equation}
\label{eq: John Wick 3}    
\frac{\partial \Psi_2^i}{\partial \psi^i}(\vb*{\psi}, \vb{y}) 
= - \sum_{j \ne i} \frac{\partial \Psi_2^j}{\partial \psi^i}(\vb*{\psi}, \vb{y}), 
\qquad    
\frac{\partial \Psi_2^i}{\partial \vb{y}^i}(\vb*{\psi}, \vb{y}) 
= - \sum_{j \ne i} \frac{\partial \Psi^j_2}{\partial \vb{y}^i}(\vb*{\psi}, \vb{y}). 
\end{equation}
Note that $(\widehat{\vb*{\psi}}(\vb*{w},\vb{z}),\widehat{\vb{y}}(\vb{z})) \in U_2$ for all $(\vb*{w},\vb{z}) \in U$.
Therefore, by the Chain Rule and the smoothness of $\widehat{\vb*{\psi}}$ and $\widehat{\vb{y}}$, it follows that $\vb{\Psi}$ is also continuously differentiable. 

Next we compute $\partial \Psi^i / \partial w^j$ for $i \ne j$.
For brevity let
\[
L_2^{ij} = L^i_2(\widehat{\vb*{\psi}}(\vb*{w},\vb{z}),\widehat{\vb{y}}(\vb{z})) \cap L^j_2(\widehat{\vb*{\psi}}(\vb*{w},\vb{z}),\widehat{\vb{y}}(\vb{z})).
\]
Then
\begin{align}
\nonumber
\frac{\partial \Psi^i}{\partial w^j}(\vb*{w},\vb{z}) 
& = 
\mathrm{det}(D \Phi) \sum_{k=1}^N \frac{\partial \Psi^i_2}{\partial \psi^k} (\widehat{\vb*{\psi}}(\vb*{w},\vb{z}),\widehat{\vb{y}}(\vb{z})) \frac{\partial \widehat{\psi}^k}{\partial w^j}(\vb*{w},\vb{z}) 
\\
\nonumber
& = \mathrm{det}(D \Phi) \frac{\partial \Psi^i_2}{\partial \psi^j} (\widehat{\vb*{\psi}}(\vb*{w},\vb{z}),\widehat{\vb{y}}(\vb{z}))
\\
\label{eq: a}
& = - \frac{\mathrm{det}(D \Phi)}{2 \| \widehat{\vb{y}}^i(\vb{z}) - \widehat{\vb{y}}^j(\vb{z}) \|}\int_{L^{ij}
_2} 
\zeta \circ \Phi \, \dd  \mathcal{H}^2.
\end{align}
By the generalised area formula (see, e.g., \cite[Theorem 2.91]{AbrosioFuscoPallara:2000} or \cite[Equation (2.7)] {Cagnetti2008}), 
\begin{equation}
\label{eq: b}
\int_{L^{ij}_2} \zeta \circ \Phi \, \dd \mathcal{H}^2
= \int_{\Phi(L^{ij}_2)} \frac{\zeta}{J_{\Phi} \circ \Phi^{-1}} \, \dd \mathcal{H}^2
=
\int_{L^i_c(\vb*{w},\vb{z}) \cap L^j_c(\vb*{w},\vb{z})} \frac{\zeta}{J_{\Phi} \circ \Phi^{-1}} \, \dd \mathcal{H}^2,
\end{equation}
where $J_\Phi:\R^3 \to \R$ is the Jacobian
\[
J_\Phi = \| (D \Phi)^{-T} \vb{n} \| \, \mathrm{det} (D\Phi), 
\qquad 
\vb{n} = \frac{\widehat{\vb{y}}^j(\vb{z}) - \widehat{\vb{y}}^i(\vb{z})}{ \| \widehat{\vb{y}}^j(\vb{z}) - \widehat{\vb{y}}^i(\vb{z}) \|}.
\]
Note that $\vb{n}$ is the outer unit normal to $L^i_2(\widehat{\vb*{\psi}}(\vb*{w},\vb{z}),\widehat{\vb{y}}(\vb{z}))$ on the face $L^{ij}_2$.
We will prove below that $J_\Phi \circ \Phi^{-1}$ is uniformly bounded from below by a positive constant, hence the integrals in \eqref{eq: b} are well defined and the generalised area formula is valid.

By equation \eqref{eq: Ezra Collective - Charlie's favourite band},
\begin{equation}
\label{eq: Thunderbolts*}    
c(\vb{x},\vb{z}^i)-\vb*{w}^i =
- 2 \Phi^{-1}(\vb{x}) \cdot \widehat{\vb{y}}^i(\vb{z}) + \| \widehat{\vb{y}}^i(\vb{z}) \|^2  
 - \widehat{\psi}^i(\vb*{w},\vb{z}) .
\end{equation}
Differentiating this with respect to $\vb{x}$ gives
\[
\nabla_{\vb{x}}c(\vb{x},\vb{z}^i)   =  - 2 (D\Phi^{-1}(\vb{x}))^T \widehat{\vb{y}}^i(\vb{z})
=  - 2 (D \Phi(\Phi^{-1}(\vb{x}))^{-T}
\widehat{\vb{y}}^i(\vb{z}).
\]
Therefore
\begin{equation}
\label{eq: c}    
\| \nabla_{\vb{x}}c(\vb{x},\vb{z}^i)   - \nabla_{\vb{x}}c(\vb{x},\vb{z}^j) \| = \frac{2 J_\Phi(\Phi^{-1}(\vb{x})) \| \widehat{\vb{y}}^j(\vb{z}) - \widehat{\vb{y}}^i(\vb{z}) \|}{\mathrm{det}(D \Phi)}.   
\end{equation}
Combining equations \eqref{eq: a}, \eqref{eq: b} and \eqref{eq: c} proves \eqref{eq: 1}, as required.

Next we prove a lower bound on $J_\Phi \circ \Phi^{-1}$.
Let $\vb{z} \in D_0^N$ and $\eta = \min_{i \ne j} |z_3^i - z_3^j|>0$.
Then 
\begin{equation}
\label{eq:lb}
    \|
\nabla_{\vb{x}}c(\vb{x},\vb{z}^i)-\nabla_{\vb{x}}c(\vb{x},\vb{z}^j)
    \|
    \ge | \partial_{x_3}  c(\vb{x},\vb{z}^i)-\partial_{x_3} c(\vb{x},\vb{z}^j) | 
    = g \frac{|z_3^i-z_3^j|}{z_3^i z_3^j} \ge g \eta \delta^2
\end{equation}
since $\vb{z}^i,\vb{z}^j \in \Y = \R^2 \times (\delta,1/\delta)$.
By \eqref{eq: c}, for all $\vb{x} \in \X$,
\[
J_\Phi(\Phi^{-1}(\vb{x}))
= \frac{\mathrm{det}(D \Phi) \| \nabla_{\vb{x}}c(\vb{x},\vb{z}^i)   - \nabla_{\vb{x}}c(\vb{x},\vb{z}^j) \|}{2 \| \widehat{\vb{y}}^j(\vb{z}) - \widehat{\vb{y}}^i(\vb{z}) \|}
\ge \frac{\f^{-4} \, \eta \, \delta^2}{2 \| \widehat{\vb{y}}^j(\vb{z}) - \widehat{\vb{y}}^i(\vb{z}) \|}, 
\]
which is the required uniform lower bound.

Next we prove \eqref{eq: 2}. 
Differentiating \eqref{eq: Thunderbolts*} with respect to $\vb{z}^i$ gives
\begin{equation}
\label{eq: Andor Season 2}   \nabla_{\vb{y}} c(\vb{x},\vb{z}^i) =
2 \left( \frac{\partial \widehat{\vb{y}}^i}{\partial \vb{z}^i} (\vb{z}) \right)^T ( - \Phi^{-1}(\vb{x}) + \widehat{\vb{y}}^i(\vb{z}) ) - \frac{\partial \widehat{\psi}^i}{\partial \vb{z}^i}(\vb*{w},\vb{z}).
\end{equation}
By \eqref{eq: Jack Quaid}, for all $i \ne j$,
\begin{align*}
& \frac{\partial \Psi^i}{\partial \vb{z}^j}(\vb*{w},\vb{z}) =
\\
& = 
\mathrm{det}(D \Phi) 
\left( 
\frac{\partial \Psi^i_2}{\partial \psi^j}
(\widehat{\vb*{\psi}}(\vb*{w},\vb{z}),\widehat{\vb{y}}(\vb{z}))
\frac{\partial \widehat{\psi}^j}{\partial \vb{z}^j} (\vb*{w},\vb{z}) +
\left( \frac{\partial \widehat{\vb{y}}^j}{\partial \vb{z}^j} (\vb{z}) \right)^T
\frac{\partial \Psi^i_2}{\partial \vb{y}^j}
(\widehat{\vb*{\psi}}(\vb*{w},\vb{z}),\widehat{\vb{y}}(\vb{z}))
\right)
\\
& = 
- \frac{\mathrm{det}(D \Phi)}{2 \| \widehat{\vb{y}}^i(\vb{z}) - \widehat{\vb{y}}^j(\vb{z}) \|}\int_{L_2^{ij}} 
\zeta(\Phi(\vb{x})) \, \dd  \mathcal{H}^2 (\vb{x}) \, 
\frac{\partial \widehat{\psi}^j}{\partial \vb{z}^j} (\vb*{w},\vb{z}) 
\\
& \quad +
\left( \frac{\partial \widehat{\vb{y}}^j}{\partial \vb{z}^j} (\vb{z}) \right)^T
\frac{\mathrm{det}(D \Phi)}{\| \widehat{\vb{y}}^i(\vb{z}) - \widehat{\vb{y}}^j(\vb{z}) \|} \int_{L_2^{ij}} 
(\widehat{\vb{y}}^j(\vb{z}) - \vb{x}) \zeta(\Phi(\vb{x})) \, \dd  \mathcal{H}^2(\vb{x})
\\
& = 
- \frac{\mathrm{det}(D \Phi)}{2 \| \widehat{\vb{y}}^i(\vb{z}) - \widehat{\vb{y}}^j(\vb{z}) \|}\int_{\Phi(L_2^{ij})} 
\frac{\zeta(\vb{x})}{J_\Phi(\Phi^{-1}(\vb{x}))} \, \dd  \mathcal{H}^2 (\vb{x}) \, 
\frac{\partial \widehat{\psi}^j}{\partial \vb{z}^j} (\vb*{w},\vb{z}) 
\\
& \quad +
\left( \frac{\partial \widehat{\vb{y}}^j}{\partial \vb{z}^j} (\vb{z}) \right)^T
\frac{\mathrm{det}(D \Phi)}{\| \widehat{\vb{y}}^i(\vb{z}) - \widehat{\vb{y}}^j(\vb{z}) \|} \int_{\Phi(L_2^{ij})} 
\frac{(\widehat{\vb{y}}^j(\vb{z})- \Phi^{-1}(\vb{x})) \zeta(\vb{x})}{J_\Phi(\Phi^{-1}(\vb{x}))} \, \dd  \mathcal{H}^2(\vb{x})
\\
& = \frac{\mathrm{det}(D \Phi)}{2 \| \widehat{\vb{y}}^i(\vb{z}) - \widehat{\vb{y}}^j(\vb{z}) \|}
\int_{\Phi(L^{ij}_2)}
\frac{\zeta}{J_\Phi \circ \Phi^{-1}}
\left( 
2 \left( 
\frac{\partial 
\widehat{\vb{y}}^j}{\partial \vb{z}^j} (\vb{z}) \right)^T (  \widehat{\vb{y}}^j(\vb{z})
- \Phi^{-1}
) 
-\frac{\partial \widehat{\psi}^j}{\partial \vb{z}^j} (\vb*{w},\vb{z}) 
\right) \, \dd \mathcal{H}^2.
\end{align*}
Combining this with 
\eqref{eq: c} and
\eqref{eq: Andor Season 2} proves \eqref{eq: 2}, as desired.

Finally, \eqref{eq: 3} follows by differentiating
\[
\sum_{j=1}^N \Psi^j (\vb*{w},\vb{z}) = \int_{\X} \zeta(\vb{x}) \, \dd \vb{x}
\]
with respect to $w^i$ and $\vb{z}^i$.
\end{proof}

\begin{rmk}[The well-preparedness assumption]
\label{remark: B.2}
The estimate \eqref{eq:lb}, and hence Lemma \ref{cor: regularity}, could be proved under the weaker assumption that the seeds are merely distinct, $\vb{z} \in D^N$, without the well-preparedness assumption $\vb{z} \in D_0^N$. This is done by keeping the $x_1$ and $x_2$ components of $\nabla_{\vb{x}}c(\vb{x},\vb{z}^i)-\nabla_{\vb{x}}c(\vb{x},\vb{z}^j)$ in the inequality \eqref{eq:lb}. It follows that Lemma \ref{lem:dualdifferentiability} also holds without the well-preparedness assumption, i.e., we could replace $D_0^N$ by $D^N$ in the definition of $U$ in both Lemma \ref{cor: regularity} and Lemma \ref{lem:dualdifferentiability}. However, since the well-preparedness assumption is crucial for the proof of Proposition \ref{prop:ODESolutions} (to ensure that seeds do not collide), we keep it in Lemmas \ref{cor: regularity} and \ref{lem:dualdifferentiability} to simplify the proof of \eqref{eq:lb}.
\end{rmk}

Finally we are in a position to prove Lemma \ref{lem:dualdifferentiability}.




\begin{proof}[Proof of Lemma \ref{lem:dualdifferentiability}]
For each $i \in \{1,\ldots,N\}$, 
define the auxiliary function $h^i: U \times U \to \mathbb{R}$ by 
\[
h^i(\vb*{w}_1, \vb{z}_1, \vb*{w}_2,  \vb{z}_2)
=
\int_{L^i_c(\vb*{w}_1, \vb{z}_1)} \zeta(\vb{x},\vb*{w}_2,\vb{z}_2) \, \dd \vb{x}.
\]
Then $\Psi^i(\vb*{w},\vb{z}) = h^i(\vb*{w},\vb{z},\vb*{w},\vb{z})$.
Fix an arbitrary point $(\vb*{w}_0,\vb{z}_0) \in U$ and let $K \subset U$ be a compact neighbourhood of 
$(\vb*{w}_0,\vb{z}_0)$ such that $|z^i_3 - z^j_3| \ge \eta$ for all $i,j \in \{1,\ldots,N\}$, $i \ne j$, 
for some $\eta > 0$.
We will prove that $h^i$ is continuously differentiable at $(\vb*{w}_0,\vb{z}_0,\vb*{w}_0,\vb{z}_0)$, and hence $\Psi^i$ is continuously differentiable at $(\vb*{w}_0,\vb{z}_0)$. 
To be precise, firstly we will prove that
\begin{gather}
\label{eq: cont1}
\forall \; (\vb*{w}_2,\vb{z}_2) \in U, \;
\frac{\partial h^i}{\partial \vb*{w}_1}(\cdot,\cdot,\vb*{w}_2,\vb{z}_2) \text{ and }
\frac{\partial h^i}{\partial \vb{z}_1}(\cdot,\cdot,\vb*{w}_2,\vb{z}_2)
\text{ are continuous in $U$}, 
\\
\label{eq: eqcont1}
\left\{
\frac{\partial h^i}{\partial \vb*{w}_1}(\vb*{w}_1,\vb{z}_1,\cdot,\cdot) \right\}_{(\vb*{w}_1,\vb{z}_1) \in K}
\text{ is uniformly Lipschitz continuous in $K$}, 
\\
\label{eq: eqcont1 z}
\left\{ \frac{\partial h^i}{\partial \vb{z}_1}(\vb*{w}_1,\vb{z}_1,\cdot,\cdot) \right\}_{(\vb*{w}_1,\vb{z}_1) \in K}
\text{ is uniformly Lipschitz continuous in $K$}.
\end{gather}
Here \emph{uniformly Lipschitz} means that the Lipschitz constants are independent of $(\vb*{w}_1,\vb{z}_1) \in K$.
It follows from \eqref{eq: cont1}--\eqref{eq: eqcont1 z} that the partial derivatives $\partial h^i/\partial \vb*{w}_1$ and $\partial h^i/\partial \vb{z}_1$ are continuous at $(\vb*{w}_0,\vb{z}_0,\vb*{w}_0,\vb{z}_0)$.
Secondly we will prove that 
\begin{gather}
\label{eq: cont2}
\forall \; (\vb*{w}_1,\vb{z}_1) \in U, \; \frac{\partial h^i}{\partial \vb*{w}_2}(\vb*{w}_1,\vb{z}_1,\cdot,\cdot) \text{ and }
\frac{\partial h^i}{\partial \vb{z}_2}(\vb*{w}_1,\vb{z}_1,\cdot,\cdot)
\text{ are continuous in $U$},
\\
\label{eq: eqcont2}
\left\{ \frac{\partial h^i}{\partial \vb*{w}_2}(\cdot,\cdot,\vb*{w}_2,\vb{z}_2)
\right\}_{(\vb*{w}_2,\vb{z}_2) \in K}
\text{ and }
\left\{ \frac{\partial h^i}{\partial \vb{z}_2}(\cdot,\cdot,\vb*{w}_2,\vb{z}_2)
\right\}_{(\vb*{w}_2,\vb{z}_2) \in K}
\text{ are equicontinuous in } K.
\end{gather}
Therefore the partial derivatives $\partial h^i/\partial \vb*{w}_2$ and $\partial h^i/\partial \vb{z}_2$ are also continuous at $(\vb*{w}_0,\vb{z}_0,\vb*{w}_0,\vb{z}_0)$.

\textit{Step 1: proof of \eqref{eq: cont1}--\eqref{eq: eqcont1 z}.}
By Lemma \ref{cor: regularity}, the function $h^i(\cdot,\cdot,\vb*{w}_2,\vb{z}_2)$ is continuously differentiable in $U$ for all 
$(\vb*{w}_2,\vb{z}_2) \in U$. (Here we have used the assumption in Lemma \ref{lem:dualdifferentiability} that $\zeta(\cdot,\vb*{w}_2,\vb{z}_2)$ is continuous.)
This proves \eqref{eq: cont1}.
Moreover, we read off from Lemma \ref{cor: regularity} that, for all $j \in \{ 1, \ldots , N\} \ \setminus \{ i \}$,
\begin{align}
\label{eq: Mikey Madison}
\frac{\partial h^i}{\partial w_1^j}(\vb*{w}_1, \vb{z}_1, \vb*{w}_2, \vb{z}_2) 
& = 
- \int_{L^i_c(\vb*{w}_1,\vb{z}_1) \cap L^j_c(\vb*{w}_1,\vb{z}_1)} 
\frac{\zeta(\vb{x},\vb*{w}_2,\vb{z}_2)}{\| \nabla_{\vb{x}} c(\vb{x},\vb{z}^i_1) -  \nabla_{\vb{x}} c(\vb{x},\vb{z}^j_1) \|} \, \dd  \mathcal{H}^2(\vb{x}),
\\
\frac{\partial h^i}{\partial \vb{z}_1^j}(\vb*{w}_1, \vb{z}_1, \vb*{w}_2, \vb{z}_2) 
& = 
\int_{L^i_c(\vb*{w}_1,\vb{z}_1) \cap L^j_c(\vb*{w}_1,\vb{z}_1)} 
\frac{\nabla_{\vb{y}}c(\vb{x},\vb{z}^j_1) \, \zeta(\vb{x},\vb*{w}_2,\vb{z}_2)}{\| \nabla_{\vb{x}} c(\vb{x},\vb{z}^i_1) -  \nabla_{\vb{x}} c(\vb{x},\vb{z}^j_1) \|} \, \dd  \mathcal{H}^2(\vb{x}),
\end{align}
and
\begin{align}
\label{eq: Anora}
\frac{\partial h^i}{\partial w_1^i}(\vb*{w}_1, \vb{z}_1, \vb*{w}_2, \vb{z}_2) 
& = - \sum_{j \ne i} \frac{\partial h^j}{\partial w_1^i}(\vb*{w}_1, \vb{z}_1, \vb*{w}_2, \vb{z}_2), 
\\    
\label{eq: The Substance}
\frac{\partial h^i}{\partial \vb{z}_1^i}(\vb*{w}_1, \vb{z}_1, \vb*{w}_2, \vb{z}_2) 
& = - \sum_{j \ne i} \frac{\partial h^j}{\partial \vb{z}_1^i}(\vb*{w}_1, \vb{z}_1, \vb*{w}_2, \vb{z}_2). 
\end{align}
Next we prove \eqref{eq: eqcont1}.
For all $\vb{x} \in \X$, $(\vb*{w}_1,\vb{z}_1) \in K$,
\[
\|
\nabla_{\vb{x}}c(\vb{x},\vb{z}^i_1)-\nabla_{\vb{x}}c(\vb{x},\vb{z}^k_1)
    \|
\ge g \eta \delta^2
\]
by \eqref{eq:lb}.
For all $(\vb*{w}_1,\vb{z}_1), \, (\vb*{w}_2,\vb{z}_2), \, (\widetilde{\vb*{w}}_2,\widetilde{\vb{z}}_2) \in K$ and all $i,j \in \{1,\ldots,N\}$, $i \ne j$,
\begin{align}
\nonumber
    & \left| 
    \frac{\partial h^i}{\partial w^j_1}
    \left(\vb*{w}_1,\vb{z}_1,\vb*{w}_2,\vb{z}_2\right) - \frac{\partial h^i}{\partial w^j_1}
    \left(\vb*{w}_1,\vb{z}_1,\widetilde{\vb*{w}}_2,\widetilde{\vb{z}}_2\right) \right| 
    \\
    \nonumber
    & \quad 
    \leq \int_{L^i_c (\vb*{w}_1,\vb{z}_1) \cap L^j_c (\vb*{w}_1,\vb{z}_1)} \frac{|\zeta (\vb{x},\vb*{w}_2,\vb{z}_2)-\zeta (\vb{x},\widetilde{\vb*{w}}_2,\widetilde{\vb{z}}_2)|}{\| \nabla_{\vb{x}} c(\vb{x},\vb{z}^i_1) -  \nabla_{\vb{x}} c(\vb{x},\vb{z}^j_1) \|} \,\dd \mathcal{H}^2(\vb{x}) 
    \\
    \label{eq: trade tarrifs}
    & \quad \leq \mathcal{H}^2 (L^i_c (\vb*{w}_1,\vb{z}_1) \cap L^j_c (\vb*{w}_1,\vb{z}_1))
    \frac{L(K)}{g\eta\delta^2} \| (\vb*{w}_2,\vb{z}_2)-(\widetilde{\vb*{w}}_2,\widetilde{\vb{z}}_2) \|,
\end{align}
where $L(K)$ was defined in the statement of Lemma \ref{lem:dualdifferentiability}. 
Define
\[
f^{ij} = \left\{\vb{x}\in\R^3 : c(\vb{x},\vb{z}_1^i) - w_1^i = c(\vb{x},\vb{z}_1^j) - w^j_1\right\}.
\]
Note that
$L^i_c (\vb*{w}_1,\vb{z}_1) \cap L^j_c (\vb*{w}_1,\vb{z}_1) \subseteq f^{ij}$.
By \cite[Theorem 2.8]{EvansGariepy:2015},
\begin{align}
\nonumber
    \mathcal{H}^2 (L^i_c (\vb*{w}_1,\vb{z}_1) \cap L^j_c (\vb*{w}_1,\vb{z}_1)) 
    & \leq \mathcal{H}^2 (\X \cap f^{ij})
    \\
    \nonumber
    & \leq \mathrm{Lip}(\Phi |_{\X}) ^2 \mathcal{H}^2 \left( \Phi^{-1}(\X \cap f^{ij}) \right) 
    \\
    \nonumber
    & = \mathrm{Lip}(\Phi |_{\X}) ^2 \mathcal{H}^2 \left( \Phi^{-1}(\X) \cap \Phi^{-1} (f^{ij}) \right) 
    \\
    \label{eq: Emilia Perez}
    & \leq \mathrm{Lip}(\Phi |_{\X}) ^2 \, \pi \left( \frac{\mathrm{diam}(\Phi^{-1}(\X))}{2} \right)^2,
\end{align}
where in the last line we used the fact that $\Phi^{-1}(f^{ij})$ is a plane (by Lemma \ref{lemma: Loeper's condition})
and the isodiametric inequality.
Combining \eqref{eq: trade tarrifs} and \eqref{eq: Emilia Perez} proves that $\partial h^i/\partial w_1^j(\vb*{w}_1,\vb{z}_1,\cdot,\cdot)$ is Lipschitz continuous in $K$ with Lipschitz constant independent of $(\vb*{w}_1,\vb{z}_1) \in K$.
By \eqref{eq: Anora}, the same holds for $\partial h^i/\partial w_1^i(\vb*{w}_1,\vb{z}_1,\cdot,\cdot)$.
This proves \eqref{eq: eqcont1}.
The proof of \eqref{eq: eqcont1 z} is very similar.

\textit{Step 2: proof of \eqref{eq: cont2}, \eqref{eq: eqcont2}.}
First we show that the partial derivative $\partial h^i / \partial \vb*{w}_2$ exists.
Fix $(\vb*{w}_1,\vb{z}_1,\vb*{w}_2,\vb{z}_2) \in U \times U$ and let $k \in \{1,\ldots,N\}$ and $\varepsilon \in (0,1)$. 
To show that $h^i$ is differentiable with respect to $\vb*{w}_2$, consider the difference quotient
\[
    \frac{h^i(\vb*{w}_1,\vb{z}_1,\vb*{w}_2+\varepsilon\vb{e}_k, \vb{z}_2)-h^i(\vb*{w}_1,\vb{z}_1,\vb*{w}_2,\vb{z}_2)}{\varepsilon}
    = \int_{L^i_c(\vb*{w}_1,\vb{z}_1)}\frac{\zeta(\vb{x}, \vb*{w}_2+\varepsilon\vb{e}_k, \vb{z}_2)-\zeta(\vb{x}, \vb*{w}_2, \vb{z}_2)}{\varepsilon} \,\dd \vb{x},
\]
where $\{\vb{e}_j\}_{j=1}^N$ is the standard basis of $\R^N$. Note that
\[
\lim_{\varepsilon \to 0}
\frac{\zeta(\vb{x}, \vb*{w}_2+\varepsilon\vb{e}_k, \vb{z}_2)-\zeta(\vb{x}, \vb*{w}_2, \vb{z}_2)}{\varepsilon}
= 
\frac{\partial \zeta}{\partial \vb*{w}} (\vb{x},\vb*{w}_2,\vb{z}_2)
\]
for all $\vb{x} \in S(\vb*{w}_2,\vb{z}_2)$, in particular, for $\mathcal{L}^3$-almost every $\vb{x} \in L^i_c(\vb*{w}_1,\vb{z}_1)$. Moreover,
\[
\left|
\frac{\zeta(\vb{x}, \vb*{w}_2+\varepsilon\vb{e}_k, \vb{z}_2)-\zeta(\vb{x}, \vb*{w}_2, \vb{z}_2)}{\varepsilon}
\right|
\le L(\widetilde{K}),
\]
where $\widetilde{K}= \{ (\vb*{w}_2+\rho \vb{e}_k, \vb{z}_2) : \rho \in [0,1] \}$ and $L$ was defined in the statement of Lemma \ref{lem:dualdifferentiability}.
Therefore, by the Lebesgue Dominated Convergence Theorem,
\begin{equation}
\label{eq: The Brutalist}
\frac{\partial h^i}{\partial \vb*{w}_2}(\vb*{w}_1,\vb{z}_1,\vb*{w}_2,\vb{z}_2) = 
\int_{L^i_c(\vb*{w}_1,\vb{z}_1)} \frac{\partial \zeta}{\partial \vb*{w}} (\vb{x},\vb*{w}_2,\vb{z}_2) \,\dd \vb{x}.
\end{equation}
Similarly, 
\begin{equation}
\label{eq: BAFTA}
\frac{\partial h^i}{\partial \vb{z}_2}(\vb*{w}_1,\vb{z}_1,\vb*{w}_2,\vb{z}_2) = 
\int_{L^i_c(\vb*{w}_1,\vb{z}_1)} \frac{\partial \zeta}{\partial \vb{z}} (\vb{x},\vb*{w}_2,\vb{z}_2) \,\dd \vb{x}.
\end{equation}
Then \eqref{eq: cont2} follows from the assumption that the partial derivatives 
$\partial \zeta/\partial \vb*{w}(\vb{x},\cdot,\cdot)$, $\partial \zeta/\partial \vb{z}(\vb{x},\cdot,\cdot)$ are continuous at $(\vb*{w}_2,\vb{z}_2)$ for all $\vb{x} \in S(\vb*{w}_2,\vb{z}_2)$, along with another application of the Lebesgue Dominated Convergence Theorem.


Finally we prove \eqref{eq: eqcont2}. 
For all $(\vb*{w}_1,\vb{z}_1), \, (\widetilde{\vb*{w}}_1,\widetilde{\vb{z}}_1), \, (\vb*{w}_2,\vb{z}_2) \in K$,
\begin{align}
\nonumber
& \left\| \frac{\partial h^i}{\partial \vb*{w}_2}(\widetilde{\vb*{w}}_1,\widetilde{\vb{z}}_1,\vb*{w}_2,\vb{z}_2) - \frac{\partial h^i}{\partial \vb*{w}_2}(\vb*{w}_1,\vb{z}_1,\vb*{w}_2,\vb{z}_2) \right\|
\\
\nonumber
& \quad = \bigg\| \int_{L^i_c(\widetilde{\vb*{w}}_1,\widetilde{\vb{z}}_1)} \frac{\partial \zeta}{\partial \vb*{w}} (\vb{x},\vb*{w}_2,\vb{z}_2) \,\dd\vb{x} - \int_{L^i_c(\vb*{w}_1,\vb{z}_1)} \frac{\partial \zeta}{\partial \vb*{w}} (\vb{x},\vb*{w}_2,\vb{z}_2) \, \dd \vb{x} \, \bigg \|
\\
\nonumber
& \quad \le \max_{\vb{x} \in \X} \left\| \frac{\partial \zeta}{\partial \vb*{w}} (\vb{x},\vb*{w}_2,\vb{z}_2) \right\| 
\left\|
\chi_{L^i_c(\widetilde{\vb*{w}}_1,\widetilde{\vb{z}}_1)}-\chi_{L^i_c(\vb*{w}_1,\vb{z}_1)}
\right\|_{L^1(\X)}
\\
\label{eq: Twitter}
& \quad \le 
L(K)
\left\|
\chi_{L^i_c(\widetilde{\vb*{w}}_1,\widetilde{\vb{z}}_1)}-\chi_{L^i_c(\vb*{w}_1,\vb{z}_1)}
\right\|_{L^1(\X)}.
\end{align}
It can be read off from the proof of \cite[Proposition 38(vii)]{Merigot:2020} that 
\[
\|
\chi_{L^i_c(\widetilde{\vb*{w}}_1,\widetilde{\vb{z}}_1)}-\chi_{L^i_c(\vb*{w}_1,\vb{z}_1)}
\|_{L^1(\X)} \to 0 \quad \text{as} \quad
(\widetilde{\vb*{w}}_1,\widetilde{\vb{z}}_1) \to (\vb*{w}_1,\vb{z}_1).
\]
Then, by \eqref{eq: Twitter}, $\partial h^i/\partial \vb*{w}_2 (\cdot,\cdot,\vb*{w}_2,\vb{z}_2)$ is continuous on the compact set $K$, hence uniformly continuous, and moreover its modulus of continuity is independent of $(\vb*{w}_2,\vb{z}_2)$. 
The same holds for $\partial h^i/\partial \vb{z}_2 (\cdot,\cdot,\vb*{w}_2,\vb{z}_2)$. This proves \eqref{eq: eqcont2}.


\emph{Step 3.} The expressions for the partial derivatives of $\Psi^i$ in the statement of Lemma \ref{lem:dualdifferentiability} can be read off immediately from \eqref{eq: Mikey Madison}--\eqref{eq: The Substance},  
\eqref{eq: The Brutalist} and \eqref{eq: BAFTA}.
\end{proof}

\bibliographystyle{abbrv}
\bibliography{references_final}

\end{document}